\documentclass[final]{siamltex}
\usepackage{amsmath,amsfonts,amssymb}
\usepackage{latexsym}
\usepackage{graphicx,overpic}
\usepackage{subfig}
\usepackage[table,xcdraw]{xcolor}
\usepackage{pgfplots}
\pgfplotsset{compat=newest}
\usepackage{color}
\usepackage[hidelinks]{hyperref}
\hypersetup{
    colorlinks=true,
    linkcolor={red!50!black},
    citecolor={blue!50!black},
    urlcolor={blue!80!black}
}

\newcommand{\R}{\mathbb R}

\newcommand{\bH}{\mathbf H}

\newcommand{\bP}{\mathbf P}

\newcommand{\bj}{\mathbf j}
\newcommand{\blf}{\mathbf f}
\newcommand{\bn}{\mathbf n}
\newcommand{\bm}{\mathbf m}

\newcommand{\bp}{\mathbf p}

\newcommand{\bu}{\mathbf u}

\newcommand{\bv}{\mathbf v}
\newcommand{\bw}{\mathbf w}

\newcommand{\bx}{\mathbf x}
\newcommand{\by}{\mathbf y}

\newcommand{\T}{\mathcal T}

\newcommand{\divG}{{\mathop{\,\rm div}}_{\Gamma}}
\newcommand{\gradG}{\nabla_{\Gamma}}
\newcommand{\nablaG}{\nabla_{\Gamma}}
\newcommand{\laplG}{\Delta_{\Gamma}}

\newcommand{\cO}{\mathcal O}

\newcommand{\Gs}{\mathcal{S}} 
\newcommand{\OGamma}{\Omega^\Gamma_h}

\renewcommand{\div}{\textrm{div}\ \!}

\newcommand{\tr}{{\rm tr}}
\DeclareGraphicsExtensions{.pdf,.eps,.ps,.eps.gz,.ps.gz,.eps.Y}

\def\cl {\nonumber \\}
\def\el {\nonumber }

\newtheorem{remark}{Remark}[section]

\begin{document}
\title{Numerical modelling of phase separation on dynamic surfaces}
\author{
Vladimir Yushutin\thanks{Department of Mathematics, University of Houston, Houston, Texas 77204 (yushutin@math.uh.edu).}
\and Annalisa Quaini\thanks{Department of Mathematics, University of Houston, Houston, Texas 77204 (quaini@math.uh.edu); Partially supported by NSF through grant  DMS-1620384 and by the  Bridge Funding Program of University of Houston.}
\and
Maxim Olshanskii\thanks{Department of Mathematics, University of Houston, Houston, Texas 77204 (molshan@math.uh.edu); Partially supported by NSF through grant  DMS-1717516.}
}
\maketitle
%
\begin{abstract} The paper presents a model of lateral phase separation in a two component material surface.
The resulting fourth order nonlinear PDE can be seen as a Cahn-Hilliard equation posed on a time-dependent surface.
Only elementary tangential calculus and  the embedding of the surface in $\mathbb{R}^3$ are used to formulate the model,
thereby facilitating the development of a fully Eulerian discretization method to solve the problem numerically.
A  hybrid method, finite difference in time and trace finite element in space,  is introduced
and stability of its semi-discrete version is proved. The method avoids any triangulation of
the surface and uses a surface-independent background mesh to discretize the equation. Thus, the method
is capable of solving the Cahn-Hilliard equation numerically on implicitly defined surfaces
and surfaces undergoing strong deformations and topological transitions. We assess the approach
on a set of test problems and apply it to model spinodal decomposition and pattern formation on colliding surfaces.
Finally, we consider the phase separation on a sphere splitting into two droplets.
\end{abstract}
\begin{keywords}
Surface Cahn--Hilliard equation;  Evolving interfaces; TraceFEM; Membrane fusion.
 \end{keywords}

\section{Introduction}\label{sec:intro}

The Cahn--Hilliard (CH) equation on a stationary, flat domain has been introduced
in the late 50s to model phase separation in binary alloy systems \cite{Cahn_Hilliard1958,CAHN1961}.
As a prototype model for segregation of two components in a mixture, over the years it has been
used in many areas beyond materials science. Applications in the biology are particularly numerous.
For example, CH type equations have been used to model and simulate
tumor growth \cite{WISE2008524,hillhorst2015},  dynamics of  plasma membranes and  multicomponent vesicles
\cite{Baumgart_et_al2003,Lowengrub2009,sohn2010dynamics,Li_et_al2012,mercker2012multiscale,barrett2014stable,barrett2017finite}, and
lipid rafts formation \cite{garcke2016coupled}. 
In some applications, such as sorting in biological membranes,
phase separation and coarsening happen in a thin, evolving layer of self-organizing molecules,
which in continuum-based approach can be modeled as a material surface.
This motivates the recent interest in the CH equation posed on time-dependent surfaces. 

The Cahn--Hilliard equation is challenging to solve numerically due to non-linearity,
stiffness, and the presence of a fourth order derivative in space.
For some recent publications on the CH equation in planar and volumetric domains, we refer to
\cite{guillen2014second,tierra2015numerical,liu2015stabilized, cai2017error} and references therein.
The numerical solution of the CH equation posed on surfaces is further complicated
by the need to discretize tangential differential operators and to approximately recover complex shapes.
Several authors have opted for a finite difference method.
For example, the closest point finite difference method was applied to solve the
CH equation on a stationary torus  in \cite{gera2017cahn} and on more general
stationary domains in \cite{jeong2015microphase}. A finite difference method for a diffuse volumetric representation of the surface
CH equation was introduced in \cite{greer2006fourth}.
However, a finite element method (FEM) is often considered to be the most flexible numerical approach
to handle complex geometries.
Concerning the CH equation posed on a stationary surface,
the convergence of a FEM was studied in \cite{du2011finite}, where numerical examples
for a sphere and saddle surface are provided; results obtained by FEM on more general surfaces
are presented in \cite{garcke2016coupled,li2017unconditionally}.
In \cite{nitschke_voigt_wensch_2012}, solutions to the  surface
Cahn--Hilliard--Navier--Stokes equation were computed  on a sphere and torus.
All of the above references use a sharp surface representation and
a discretization mesh \emph{fitted} to the surface.
In \cite{Yushutin_IJNMBE2019}, we studied for the first time
a \emph{geometrically unfitted} finite element method for the surface Allen--Cahn and
Chan--Hilliard equations.
In the unfitted  FEM, the surface is not triangulated in the common sense and
may overlap a background computational (bulk) mesh in an arbitrary way. Moreover, stability and accuracy
of the discretization do not depend on the position of the surface relative to the bulk mesh.
Thus, the next natural development is to allow the surface to evolve through the time-independent
bulk mesh that is used to discretize the equation posed on the dynamic surface itself.
Such development is presented in this paper,
which extends and studies the unfitted FEM for CH equations on time-dependent surfaces.

Very few works deal  with a CH equation posed on an evolving surface.
In \cite{elliott2015evolving}, the authors show a rigorous well-posedness result for a continuous CH-type equation, which
is a simplification of the model for surface dissolution
set out in \cite{EILKS20089727}.
The FEM used for the space discretization in \cite{elliott2015evolving}, known as evolving surface finite element method,
relies on evolving an initial surface triangulation by moving the nodes according to a prescribed velocity.
The asymptotic limit (as the interfacial width parameter tends to zero) of the CH equation on a surface
evolving with prescribed velocity
is studied theoretically and numerically in \cite{o2016cahn}.
The authors of \cite{barrett2017finite} studied a finite element method for the bulk Navier--Stokes equations
coupled to the surface CH model, where the surface evolution
is driven by the bulk fluid dynamics and a curvature energy. 
Again, in \cite{EILKS20089727,elliott2015evolving,o2016cahn,barrett2017finite}, 
the discretization mesh is \emph{fitted} to the computational surface and evolves with it. Although this approach offers many
advantages, problems arise when strong deformations and topological changes of the surface occur.
The numerical approximation of the surface CH equation using an isogeometric approach was studied in
\cite{Zimmermann2019}.
The point of using isogeometric analysis is that spline bases
with high-order and high-continuity allow to treat the fourth order problem, without resorting to a mixed formulation.
In \cite{Zimmermann2019}, an isogeometric finite element formulation was
introduced for the CH equation written in intrinsic surface variables on a surface evolved by
the PDE derived according to Kirchhoff--Love shell theory.
Finally, we mention a different approach to study surface phase distribution \cite{Wang2008,Lowengrub2009,VALIZADEH2019599}:
a pair of phase-field variables is introduced such that one variable characterizes
the surface (in a diffuse manner) while the other describes the distribution of the surface phases.

In this paper, we use elementary tangential calculus to derive a CH equation posed on an evolving material surface.
There is a (non-essential) difference between the equation we get
compared to the one in \cite{elliott2015evolving,o2016cahn}, which is explained in Remark~\ref{Rem1}.
The more substantial difference with the previous works 
is that we study a \emph{geometrically unfitted} finite element method.
This method builds upon earlier work on an unfitted FEM, known as trace finite element method (TraceFEM), for elliptic PDEs on stationary surfaces in~\cite{ORG09}
and evolving surfaces in~\cite{lehrenfeld2018stabilized}.
Unlike some other geometrically unfitted methods for surface PDEs,
TraceFEM employs a sharp surface representation. The surface can be defined implicitly, e.g.~as the
zero of a level-set function, and
no knowledge of the surface parametrization is required.  We shall see that the developed method is very well suited
for CH equation posed on evolving surfaces. 

The work presented in here is a first step in the direction of understanding the evolution
of multicomponent lipid vesicles used as drug carriers and their fusion with the target cell.
Membrane fusion is recognized as a potentially efficient mechanism for the delivery of
macromolecular therapeutics to the cellular cytoplasm. The process of lipid membrane
phase separation to concentrate binding lipids within distinct regions of the membrane surface
has been shown to enhance membrane fusion \cite{Imam2017}.
In addition, recent developments on targeted lipid vesicles suggest that the formation of reversible
phase-separated patterns on the vesicle surface increase target selectivity, cell uptake and overall efficacy \cite{Bandekar_et_al2013,KARVE20104409}.
The numerical experiments presented in Sec.~\ref{sec:coll_dyn} mimic phase separation occurring on a two component
lipid vesicle, leading to fusion with another vesicle. These and other numerical results in
Sec.~\ref{sec:num_res} showcase the ease with which the numerical method handles changes in the surface topology,
making it a  perfect computational tool to support and complement experimental practice in the design
of drug carriers.

The outline of the paper is as follows. In Section~\ref{s_cont},
we will derive from  conservation laws a Cahn--Hilliard equation on an evolving surface.
The hybrid finite difference in time -- finite element in space variant of the
TraceFEM for the CH equation is introduced in Section~\ref{sec:num_method}.
After the assessment of the numerical method for a set of model problems, in Section~\ref{sec:num_res}
we will study phase separation modeled by our surface Cahn--Hilliard equation
on evolving surfaces with topological transitions: two colliding spheres
and a sphere splitting into two droplets.
Section~\ref{sec:concl} closes the paper with a few concluding remarks.

\section{Mathematical model}\label{s_cont}

\subsection{Preliminaries}

Let $\Gamma(t)$ be a closed, smooth, evolving surface in $\mathbb{R}^3$
for $t \in [0, T]$, where $T > 0$ is a final time.
Let the 3D bulk domain $\Omega$ be such that $\Gamma(t)\subset\Omega$, for all $t \in [0, T]$.
We assume surface $\Gamma(t)$ evolves according to a given smooth velocity field $\bu$,  i.e., $\bu(\bx,t)$ is the velocity of point $\bx\in\Gamma(t)$.  Consider the decomposition of $\bu$ into normal $(\bu_N=u_N\bn)$ and tangential $(\bu_T)$ components:
\begin{equation}
\bu = \bu_T + u_N\bn,\quad \bu_T\cdot\bn=0, \el
\end{equation}
where $\bn$ is the outward normal vector on $\Gamma$.
While the normal velocity $u_N$ completely defines the geometric evolution of the surface, we are interested in the motion
of $\Gamma(t)$ as a material surface, which also includes tangential deformations. The material deformation  can be defined through the Lagrangian mapping $\Psi(t, \cdot)$ from $\Gamma(0)$ to $\Gamma(t)$, i.e. for $\by\in\Gamma(0)$, $\Psi(t,\by)$ solves the ODE system
\begin{equation}\label{Lagrange}
\Psi(0,\by)=\by,\quad  \frac{\partial \Psi(t,\by)}{\partial t}=\bu(t,\Psi(t,\by)),\quad t\in[0,T].
\end{equation}

For any sufficiently smooth function $f$ in a neighborhood of $\Gamma(t)$ its tangential gradient is
defined as $\nablaG f=\nabla f-(\bn\cdot\nabla f)\bn$. The tangential (surface) gradient $\nablaG f$
then depends only on values of $f$ restricted to $\Gamma(t)$ and $\bn\cdot\nablaG f=0$ holds.
For a vector field $\blf$ on $\Gamma(t)$ we define $\nablaG \blf$ componentwise.
The surface divergence operator for $\blf$ and the Laplace--Beltrami operator for $f$
are given by:
\[
\divG \blf  := \tr (\gradG \blf),  \quad \laplG f := \divG (\gradG f),
 \]
where $\tr(\cdot)$ is the trace of a matrix.
For a fixed $t\in[0,T]$, $L^2(\Gamma(t))$ is the Lebesgue space of square-integrable functions on $\Gamma(t)$ and
$H^1(\Gamma(t))$ is the Sobolev space of all functions $f\in L^2(\Gamma(t))$ such that $\nabla_\Gamma f\in L^2(\Gamma(t))^3$.
For a subdomain $S(t)\subset\Gamma(t)$, we recall the integration by parts identity:
 \begin{equation}\label{int_parts}
 \int_{S(t)} f\div_\Gamma\bv \, ds=\int_{\partial S(t)} f\bv\cdot\bm\, d\gamma -\int_{S(t)} \bv\cdot\nabla_\Gamma f\, ds+\int_{S(t)} \kappa f\bv\cdot\bn \, ds,
 \end{equation}
where $\kappa$ is the sum of principle curvatures and $\bm$ is the normal vector on $\partial S(t)$ that is tangential to $\Gamma(t)$ and outward for $S(t)$;  vector field $\bv$ and scalar function $f$ are such that all
the quantities in \eqref{int_parts} exist.

Let  $S(t)\subset\Gamma(t)$ be a material area evolving according to \eqref{Lagrange}. Then, the following surface analogue of the Reynolds transport theorem holds (see, e.g., Lemma 2.1 in
\cite{DziukElliot2013a}):
\begin{align}\label{eq:transport_rhoxi}
\frac{{\rm d}}{{\rm d}t} \int_{S(t)} f  \, ds =  \int_{S(t)} \left(\dot{f} + f \divG \bu \right) \, ds,
\end{align}
where $\dot{f}$ is the material derivative of  $f$.
One can write the material derivative in terms of partial derivatives
\begin{align}\label{eq:material_der}
\dot{f} = \frac{\partial f}{\partial t} + \bu \cdot \nabla f.
\end{align}
The terms on the right-hand side of \eqref{eq:material_der} are well defined if one identifies $f$ with its arbitrary smooth extension
from the surface to  a neighborhood of the space-time manifold $\Gs\subset \mathbb{R}^4$,
\[
\Gs=\bigcup_{t\in[0,T]}\Gamma(t)\times\{t\}.
\]
Of course, $\dot{f}$ is an intrinsic surface quantity and the quantity on right-hand side of \eqref{eq:material_der} is independent of the choice of smooth extension.

Now, we are prepared to set up the mathematical model to describe the separation of two conserved phases on $\Gamma$.

\subsection{Model setup}
On $\Gamma(t)$ we consider a heterogeneous mixture of two species with densities $\rho_i$, {\small $i = 1, 2$}.
Let $\rho=\rho_1+\rho_2$ be the total density. The conservation of mass  written for an arbitrary \emph{material} area  $S(t)\subset\Gamma(t)$ and \eqref{eq:transport_rhoxi} yield
\begin{align*}
0=\frac{{\rm d}}{{\rm d}t} \int_{S(t)} \rho  \, ds =\int_{S(t)} \left(\dot{\rho} + \rho \divG \bu \right) \, ds.
\end{align*}
Since the above identity holds for  arbitrary   $S(t)\subset\Gamma(t)$, it implies
\begin{align}\label{eq:evol_rho}
\dot{\rho} + \rho \divG \bu=0\quad\text{on}~\Gamma(t).
\end{align}
Following \cite{Cahn_Hilliard1958,CAHN1961,Lowengrub1998},
to describe the dynamics of phases
we introduce specific mass concentrations $c_i = m_i/m$, {\small$i = 1, 2$}, where $m_i$
are the masses of the components and $m$ is the total mass.  Since $m = m_1 + m_2$,
it holds $c_1 + c_2 = 1$. Let $c_1$ be the representative concentration $c$ (order parameter), that is $c = c_1$
and $c\in [0,1]$. 
Mass conservation  for one component on $S(t)\subset\Gamma(t)$ takes the form
\begin{align}\label{eq:evol_rhoxi}
\frac{{\rm d}}{{\rm d}t} \int_{S(t)} \rho c  \, ds =    -\int_{\partial S(t)} \bj\cdot\bm \, d\gamma.
\end{align}
where $\bj$ is a mass flux.

By applying the transport formula  \eqref{eq:transport_rhoxi} to the left-hand side of \eqref{eq:evol_rhoxi}  and the integration by
parts formula \eqref{int_parts} to the right-hand side of \eqref{eq:evol_rhoxi}, we obtain
\begin{align*}
  \int_{S(t)} (\dot{\rho} c+ \dot{c} \rho  +{c} \rho\divG \bu) \, ds = - \int_{S(t)}  \divG \bj  \, ds
\end{align*}
Thanks to \eqref{eq:evol_rho}, the above equality simplifies to
\begin{align*}
  \int_{S(t)} \dot{c} \rho \, ds = - \int_{S(t)}  \divG \bj  \, ds
\end{align*}
Since $S(t)$ can be taken arbitrary, we get
\begin{align}\label{eq:evol_xi}
 \dot{c}   +  \frac{1}{\rho}\divG \bj = 0 \quad \text{on}~\Gamma(t),
\end{align}

The classical assumption for the flux $\bj$ is  Fick's law:
\begin{align}\label{eq:flux}
\bj = - M \gradG \mu \quad \text{on}~\Gamma, \quad \mu = \frac{\delta f}{\delta c},
\end{align}
where $M$ is the so-called mobility coefficient (see \cite{Landau_Lifshitz_1958}) and
$\mu$ is the chemical potential, which is defined as the functional derivative of the total
specific free energy $f$ with respect to the concentration $c$.
The choice of $f$ defines a specific model of mixing. Following \cite{Cahn_Hilliard1958}, we set
\begin{align}\label{eq:total_free_e}
f(c) = \frac1{\epsilon}f_0(c) + \frac{\epsilon}{2}  | \gradG c |^2.
\end{align}
In eq.~\eqref{eq:total_free_e}, $f_0(c)$ is the free energy per unit surface, a non-convex function of $c$, while the second term represents
the interfacial free energy based on the concentration gradient.
 In model \eqref{eq:evol_xi}--\eqref{eq:total_free_e},
the interface between the two components is a layer of size ${\epsilon}$.

By combining eqs.~\eqref{eq:evol_rho}, \eqref{eq:evol_xi}--\eqref{eq:total_free_e}, we obtain
the Cahn--Hilliard equation posed on evolving material surface $\Gamma(t)$:
\begin{align}
\dot{\rho} + \rho \divG \bu&=0 \quad \text{on}~\Gamma(t), \label{eq:CH_1} \\
\dot{c}  -  \rho^{-1}\divG \left(M \gradG \left(\frac1{\epsilon}f_0' - \epsilon\Delta_\Gamma c\right)\right) &= 0\quad \text{on}~\Gamma(t), \label{eq:CH_2}
\end{align}
for all $t \in (0, T)$.
We close system \eqref{eq:CH_1}--\eqref{eq:CH_2} with initial conditions
\begin{align}
\rho(\cdot, t) = \rho_0,\quad\text{and}\quad c(\cdot, t) = c_0  \quad \text{on}~\Gamma(0).
\end{align}
There are no boundary conditions since the boundary of $\Gamma(t)$ is empty.

Note that eq.~\eqref{eq:CH_1} is decoupled from eq.~\eqref{eq:CH_2},
in the sense that the evolution of the total density distribution depends only on the given surface motion, but not on the order
parameter $c$. In the particular case of only normal velocities, i.e. $\bu_T=0$, eq.~\eqref{eq:CH_1}  simplifies to
 \[
 \dot{\rho} = -\kappa u_N\rho.
 \]
For inextensible membranes, it holds $\divG \bu=0$ and eq.~\eqref{eq:CH_1} further simplifies  to $\dot\rho=0$.

\begin{remark}\rm\label{Rem1}
The system \eqref{eq:CH_1}--\eqref{eq:CH_2} is similar to the Cahn-Hilliard equation on evolving surface found in \cite{elliott2015evolving,o2016cahn},
but different. In \cite{elliott2015evolving,o2016cahn}, the equation is written for the conserved variable $u=\rho\, c$ and under the assumption that the free energy functional depends on $u$ rather than on the order parameter $c$, i.e. $f=f(u)$ and  $\mu = \frac{\delta f}{\delta u}$.
Formally, system \eqref{eq:CH_1}--\eqref{eq:CH_2} and the surface Cahn--Hilliard equation in \cite{elliott2015evolving,o2016cahn} are
the same problem for inextensible membranes and if one assumes $\rho_0=const$. Here the situation partially resembles the coupling of a two-component compressible fluid flow with dissipative Ginzburg--Landau interface dynamics, where, depending on the choice of the variables to define the energy functional, the effects of compressibility   has to be considered in the definition of
the reactive stress tensor and the chemical potential; see, e.g.,~\cite{Lowengrub1998} for further discussion and references.
\end{remark}

\subsection{Weak form of the equations}
Eq.~\eqref{eq:CH_1} is a transport problem for the total density.
It can be integrated independently of eq.~\eqref{eq:CH_2} as long as the evolution of the surface
is given. Since $\rho$ is completely determined by the initial data and $\bu$, but not by the phase composition of $\Gamma(t)$,
for the purposes of this paper we shall assume that $\rho$ is given.
Eq.~\eqref{eq:CH_2} is a fourth-order nonlinear evolutionary equation, which is more challenging
numerically. In particular, casting it in a weak form would lead to second-order spatial derivatives.
From the numerical point of view,
it is beneficial to avoid higher order spatial derivatives. Hence, following the  common practice we rewrite eq.~\eqref{eq:CH_2} in mixed form, i.e. as two coupled second-order equations:
\begin{align}
&\dot{c} -  \rho^{-1}\divG \left(M \gradG \mu \right)  = 0 \quad \text{on}~\Gamma(t) \label{eq:sys_CH1}, \\
&\mu = \frac1{\epsilon}f_0' - \epsilon \Delta_\Gamma c \quad \text{on}~\Gamma(t). \label{eq:sys_CH2}
\end{align}

System \eqref{eq:sys_CH1}--\eqref{eq:sys_CH2} needs to be supplemented with the
definitions of mobility $M$ and free energy per unit surface $f_0$. A possible choice for $M$
is given by
\begin{align}\label{eq:M}
M = M(c) = \sigma{} c(1 - c),
\end{align}
with $\sigma>0$.
This mobility is referred to as a degenerate mobility, since it is not strictly positive.
As for the free energy per unit surface, a common choice 
is given by the double-well potential
\begin{align}\label{eq:f0}
f_0(c) = \frac{1}{4} c^2(1 - c)^2.
\end{align}

We are interested in a finite element numerical method for the evolving surface Cahn--Hilliard problem
\eqref{eq:sys_CH1}--\eqref{eq:f0}. We first  need a weak  formulation of system
\eqref{eq:sys_CH1}--\eqref{eq:sys_CH2}.
To devise it, we start with multiplying \eqref{eq:sys_CH1} by $\rho$
and a smooth $v:\Gs\to\R$, while we multiply \eqref{eq:sys_CH2}
by a smooth $q:\Gs\to\R$. Next, we integrate  over $\Gamma(t)$ and employ
the integration by parts identity \eqref{int_parts} with $S(t)=\Gamma(t)$ (this implies $\partial S(t)=\emptyset$). The integration by parts is applied to the diffusion terms in \eqref{eq:sys_CH1} and \eqref{eq:sys_CH2}.
The curvature terms vanish  since the fluxes $M \nabla_\Gamma \mu$ and $\nabla_\Gamma c$ are tangential. We obtain the following equalities:
\begin{align}
&\int_{\Gamma(t)}\rho \dot{c} \,v \, ds +  \int_{\Gamma(t)} M \gradG \mu \, \gradG v \, ds = 0, \label{eq:sys_CH1_weak} \\
&\int_{\Gamma(t)}  \mu \,q \, ds - \frac1{\epsilon}\int_{\Gamma(t)} f_0'(c) \,q \, ds -  \epsilon\int_{\Gamma(t)}\gradG c \, \gradG q \, ds = 0.\label{eq:sys_CH2_weak}
\end{align}

A rigorous week formulation requires the definition of the test and trial functional spaces. Suitable  spaces
are formulated for functions defined on the space-time manifold $\Gs$: $L^q(\Gs)$, $1\le q\le\infty$ and $H^1(\Gs)$ are standard Lebesgue and Sobolev spaces on $\Gs$. We also need
\begin{align*}
  H^{1,\Gamma} &=\{v\in L^2(\Gs)\,:\, \|\nablaG v\|_{L^2(\Gs)}\le\infty\},  \\
  L^{\infty}_{1,\Gamma} &=\{v\in L^\infty(\Gs)\,:\, \mbox{ess}\sup\limits_{t\in[0,T]}\|\nablaG v\|_{L^2(\Gs)}\le\infty\}.
\end{align*}
Note that  $\nabla_\Gs$  and $\int_{\Gs}$ can be written in terms of tangential calculus on $\Gamma$ and the geometric motion of $\Gamma$, i.e. $u_N$ and $\bn$. See  \cite{olshanskii2014eulerian,elliott2015evolving} for details and properties of the above spaces.

We introduce the week formulation of \eqref{eq:sys_CH1}--\eqref{eq:sys_CH2}:\\
\textit{
 Find $c\in L^{\infty}_{1,\Gamma}\cap H^1(\Gs)$ and $\mu\in H^{1,\Gamma}$ satisfying \eqref{eq:sys_CH1_weak}--\eqref{eq:sys_CH2_weak} for almost every $t\in(0,T)$ and for all $ (v,q) \in H^{1,\Gamma} \times H^{1,\Gamma}$.}

For a closely related problem (which coincides with \eqref{eq:sys_CH1_weak}--\eqref{eq:sys_CH2_weak} if $\div_\Gamma \bu=0$, $\rho=1$, $M=\mbox{const}$), the above weak formulation was shown in \cite{elliott2015evolving} to have the unique solution.  Extending the well-posedness analysis for \eqref{eq:sys_CH1_weak}--\eqref{eq:sys_CH2_weak} is out of the scope for this paper.

\section{A hybrid finite difference\,/\,finite element numerical method}\label{sec:num_method}\label{s_FEM}

In this section, we present a completely Eulerian numerical method for solving problem \eqref{eq:CH_1}--\eqref{eq:CH_2}.
In our approach, the mesh does not follow the evolution of the surface, as is typical for
Lagrangian methods. Examples of finite element methods for surface PDEs based on the
Lagrangian description can be found, e.g., in \cite{Dziuk07,barrett2015stable,elliott2015error,elliott2015evolving,sokolov2015afc,macdonald2016computational}.
We rather  allow the surface to travel through a background mesh without  restrictions.
Furthermore, here we look for a numerical scheme that extends, in some sense, the method of lines,
which is (arguably) the most popular computational approach for parabolic problems in stationary domains.
We recal that the method of lines provides the convenience of a separated numerical treatment
of spatial and temporal variables. This is by no means straightforward for problem
\eqref{eq:CH_1}--\eqref{eq:CH_2}, since the equations are defined on a time-dependent
surface and there is no evident way of separating variables in the Eulerian setting.
To circumvent these difficulties, we build on ideas from \cite{olshanskii2017btrace,lehrenfeld2018stabilized}
and suggest a  hybrid,  finite-difference in time and finite element in space discretization of
problem \eqref{eq:CH_1}--\eqref{eq:CH_2}.

We start with several important observations about the differential problem  \eqref{eq:CH_1}--\eqref{eq:CH_2} and its alternative integral formulation.

\subsection{Extending $c$ and $\mu$ off the surface} We are  able to decompose the material derivative into the sum of Eulerian terms,
\begin{equation}\label{material}
  \dot{c}= \frac{\partial c}{\partial t} + \bu \cdot \nabla c,
\end{equation}
if we assume  an arbitrary smooth extension of $c$ to  a neighborhood of $\Gs$, denoted further by $\cO(\Gs)$. Let us define $\cO(\Gs)=\bigcup\limits_{t\in(0,T)}\cO_\delta(\Gamma(t))\times\{t\}$, where
\[
\cO_\delta(\Gamma(t))=\{\bx\in\R^3\,:\,\mbox{dist}(\bx,\Gamma(t))<\delta\}.
\]
We let $\delta$ be sufficiently small such that $\cO_\delta(\Gamma(t))\subset\Omega$ for all $t\in(0,T)$.

Among infinitely many  possible smooth extensions of $c$ to $\cO(\Gs)$,
it will be convenient to assume the one given by the close point projection on $\Gamma(t)$:
Fix $t\in(0,T)$, then for  $\bx\in \cO_\delta(\Gamma(t))$ we denote its close point projection on $\Gamma(t)$ by $\bp(\bx)$. For smoothly evolving $\Gamma(t)$,
the projection $\bp(\bx)$ is well defined in $\cO(\Gs)$, and we consider the extensions given by
 $c^e(\bx,t)=c(\bp(\bx),t)$ and $\mu^e(\bx,t)=\mu(\bp(\bx),t)$ for $(\bx,t)\in\cO(\Gs)$.
We also introduce the extension of the spatial normal vector field, $\bn^e(\bx,t)=\bn(\bp(\bx),t)$.
Recall that $\bn(\cdot,t)$ is the outward normal for $\Gamma(t)$.
Now $c^e$ and $\mu^e$ can be characterized as extensions along the directions given by the  normal field,
i.e.~the so called normal extensions:
  \begin{equation}\label{ExtNorm}
 \frac{\partial c^e}{\partial \bn^e}= \frac{\partial \mu^e}{\partial \bn^e}=0\quad\text{in}~\cO(\Gs).
 \end{equation}
If $\Gamma(t)$ is a $C^2$ surface, then its normal field is $C^1$-smooth and $c\in H^1(\Gamma(t))$
implies $c^e\in H^1(\cO_\delta(\Gamma(t)))$. In turn, the regularity of $\Gamma(t)$ together
with assumptions on its evolution (e.g., in terms of smoothness of the mapping $\Psi$)
yield the space-time regularity of $c^e$ (if $c$ is smooth); see more details in, e.g.,
\cite[Sec.~6.1]{olshanskii2014error}.

We shall identify functions defined on $\Gs$ with their normal extensions and skip the upper index $e$. Once we do this,  the weak formulation  \eqref{eq:sys_CH1_weak}--\eqref{eq:sys_CH2_weak} yields the following identities for $c$ and $\mu$:
\begin{align}
\int_{\Gamma(t)}\rho (\frac{\partial c}{\partial t} + \bu \cdot \nabla c) \,v \, ds +  \int_{\Gamma(t)} M \gradG \mu \, \gradG v \, ds
&= 0, \label{eq:sys_CH1_weak2} \\
\int_{\Gamma(t)}  \mu \,q \, ds - \frac1{\epsilon}\int_{\Gamma(t)} f_0'(c) \,q \, ds -  \epsilon\int_{\Gamma(t)}\gradG c \, \gradG q \, ds
&= 0,\label{eq:sys_CH2_weak2} \\
\frac{\partial \mu}{\partial \bn} = \frac{\partial c}{\partial \bn} &=0\quad \text{in}~\cO_\delta(\Gamma(t)), \label{eq:sys_CH3_weak2}
\end{align}
for all smooth $v$ and $q$ defined in $\Gamma(t)$. 
We shall base our discretization method on equalities  \eqref{eq:sys_CH1_weak2}--\eqref{eq:sys_CH3_weak2}. The crucial idea is that our numerical solution should approximate the concentration and chemical potential on $\Gamma(t)$ \emph{together with their normal extensions} in a suitable neighborhood of the discrete surface. This will allow for a completely Eulerian method with a separate treatment of spatial and temporal derivatives.

Next, we turn to describing the discrete formulation.

\subsection{Time discretization}

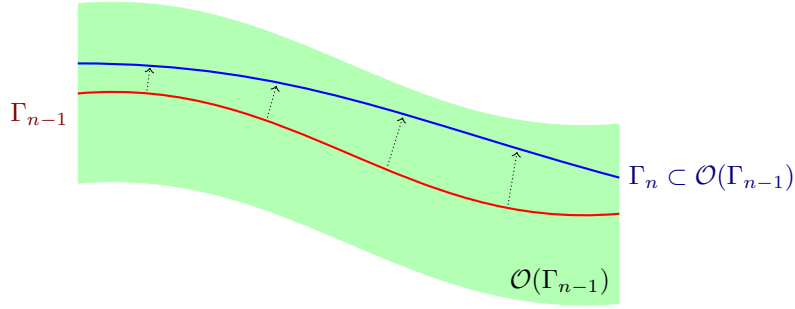
\begin{figure}[ht]
  \vspace*{-0.4cm}
  \begin{center}
    \begin{tikzpicture}[scale=0.8]
      \draw[draw=none, fill=green, opacity=0.3] (1,5)
      to[out=5,in=185] (10,3)
      -- (10,0)
      to[out=185,in=5] (1,2)
      -- (1,5)
      ;
      \draw[red,thick] (1,3.5) to[out=5,in=185] (10,1.5);
      \draw[blue,thick] (1,4) to[out=0,in=165] (10,2.1);
      \node[above left] at (10,0) {\color{black} $\cO(\Gamma_{n-1})$};
      \node[below left] at (1,3.5) {\color{red!50!black} $\Gamma_{n-1}$};
      \node[right] at (10,2.1) {\color{blue!50!black} $\Gamma_n \subset \cO(\Gamma_{n-1})$};
      \draw[->,densely dotted,black] (2.15,3.55) -- (2.2,3.925);
      \draw[->,densely dotted,black] (4.15,3.1) -- (4.3,3.625);
      \draw[->,densely dotted,black] (6.15,2.3) -- (6.4,3.1);
      \draw[->,densely dotted,black] (8.15,1.65) -- (8.3,2.525);
    \end{tikzpicture}
  \end{center}
  \vspace*{-0.4cm}
  \caption{Illustration of interface positions at different time instances satisfying condition~\eqref{ass1}.}
  \label{fig:Gammaneighborhood}
\end{figure}

In the classical method of lines, one first discretizes spatial variables and leaves time continuous; next the resulting Cauchy problem is integrated numerically. Here, we follow the reverse order  by first treating time derivatives, which is a popular alternative~\cite{Thomee97}  and often leads at the end to the same fully discrete scheme.
To this end, consider the uniform time step $\Delta t=T/N$, and let  $t_n=n\Delta t$ and $I_n=[t_{n-1},t_n)$.
Denote by $c^n$, $\mu^n$ an approximation of $c(t_n)$ and $\mu(t_n)$, define $\Gamma_n := \Gamma(t_n)$, $n=0,\dots,N$.
We assume that $\cO(\Gs)$ is a sufficiently large neighborhood of $\Gs$ such that
\begin{equation}\label{ass1}
  \Gamma_n\subset\cO_\delta(\Gamma_{n-1})\quad\text{for}~n=1,\dots,N; 
\end{equation}
see Fig.~\ref{fig:Gammaneighborhood}.
In this case, $c^{n-1}$ is well-defined on $\Gamma_n$ and we can approximate the time derivative by simple finite difference
\[
\frac{\partial c}{\partial t}(t_n)\simeq\left[{c}\right]_{t}^{n}:=\frac{c^{n}-c^{n-1}}{\Delta t}.
\]

This brings us to the implicit Euler method for problem \eqref{eq:sys_CH1_weak2}--\eqref{eq:sys_CH3_weak2}:
\begin{align}
\int_{\Gamma_n}\rho^n (\left[{c}\right]_{t}^{n} + \bu^n \cdot \nabla c^n) \,v \, ds +  \int_{\Gamma^n} M^{n} \gradG \mu^n \, \gradG v \, ds
&= 0, \label{eq:sys_CH1_im} \\
\int_{\Gamma_n}  \mu^n   \,q \, ds  -  \epsilon\int_{\Gamma_n}\gradG c^n \, \gradG q \, ds-
\frac1{\epsilon}\int_{\Gamma_n} f_0'(c^{n}) \,q \, ds &= 0,\label{eq:sys_CH2_im}\\
\frac{\partial \mu^n}{\partial \bn} = \frac{\partial c^n}{\partial \bn} =0\quad \text{in}~&\cO_\delta(\Gamma_n), \label{eq:sys_CH3_im}
\end{align}
It is important that all terms in \eqref{eq:sys_CH1_im}--\eqref{eq:sys_CH2_im} are defined on the
`current' surface $\Gamma_n$ and its neighborhood.
In particular, $c^{n-1}$ is well-defined on $\Gamma_{n}$ thanks to condition \eqref{ass1},
which gives the meaning to $\left[{c}\right]_{t}^{n}$ in \eqref{eq:sys_CH1_im}--\eqref{eq:sys_CH2_im},

In practice, one can use a semi-implicit Euler method, which compromises
the stability of the implicit method in favor of computational time savings:
\begin{align}
\int_{\Gamma_n}\rho^n (\left[{c}\right]_{t}^{n} + \bu^n \cdot \nabla c^n) \,v \, ds +  \int_{\Gamma^n} M^{n-1} \gradG \mu^n \, \gradG v \, ds 
&= 0, \label{eq:sys_CH1_num} \\
\int_{\Gamma_n}  (\mu^n - \beta_s \Delta t\left[{c}\right]_{t}^{n})  \,q \, ds  -  \epsilon\int_{\Gamma_n}\gradG c^n \, \gradG q \, ds
&= \frac1{\epsilon}\int_{\Gamma_n} f_0'(c^{n-1}) \,q \, ds,\label{eq:sys_CH2_num}\\
\frac{\partial \mu^n}{\partial \bn} = \frac{\partial c^n}{\partial \bn} &=0\quad \text{in}~\cO_\delta(\Gamma_n).\label{eq:sys_CH3_num}
\end{align}
Following the idea of first order stabilization in~\cite{Shen_Yang2010}, we included
term $ -\beta_s \Delta t\left[{c}\right]_{t}^{n}$ in eq.~\eqref{eq:sys_CH2_im}.
 On a stationary domain, for constant mobility and a slightly modified double-well potential $f_0$
 than \eqref{eq:f0}, the scheme above was shown in~\cite{Shen_Yang2010} to be
 stable for large enough $\beta_s=O(\epsilon^{-1})$.
In the semi-implicit Euler method both nonlinear terms,
i.e.~the mobility coefficient in eq.~\eqref{eq:sys_CH1_num} and the potential
in eq.~\eqref{eq:sys_CH2_num}, are extrapolated from the previous time step.
So, technically at each time step the method requires solving a linear system
\eqref{eq:sys_CH1_num}--\eqref{eq:sys_CH2_num} and further extending the solution along normals,
i.e. solving \eqref{eq:sys_CH3_num} with given data on $\Gamma_n$.
We will show in Sec.~\ref{sec:fully_discrete}
that \emph{both steps} (solving for $\mu^n$ and $c^n$ on $\Gamma_n$ and extension)
\emph{can be naturally combined in one linear solve on the finite element level}.

\subsection{Stability of the semi-discrete scheme}

It is well known that in a steady domain the Cahn-Hilliard problem defines the $H^{-1}$-gradient
flow of an energy functional. However, we are not aware of a minimization property for the
Cahn-Hilliard problem in time-dependent domains. In Section~\ref{sec:num_res},
we shall illustrate by numerical examples that the evolution
of $\Gamma(t)$ can actually produce a (local) increase of the system energy.

It is natural to say that \eqref{eq:sys_CH1_im}--\eqref{eq:sys_CH3_im} or
\eqref{eq:sys_CH1_num}--\eqref{eq:sys_CH3_num} are numerically stable
if the solution satisfy an energy bound  uniform in the discretization parameter $\Delta t$.
In this Section, we demonstrate such a bound for \eqref{eq:sys_CH1_im}--\eqref{eq:sys_CH3_im}
subject to the following simplifications. First, we assume constant density and mobility: $\rho^n=1$, $M^{n}=1$.
We also assume a tangential incompressibility condition:
 \begin{equation}\label{Incompr}
 \div_\Gamma \bu(\bx,t)=0,\quad \forall~\bx\in\Gamma(t),
 \end{equation}
 which corresponds to the motion of an \emph{inextensible} membrane.

We start by taking $v=\Delta t\, \mu^n$ in \eqref{eq:sys_CH1_num} and $q=-\Delta t(\left[{c}\right]_{t}^{n}+ \bu^n \cdot \nabla c^{n})$ in \eqref{eq:sys_CH2_num} and sum the resulting equalities.
After few cancellations, we obtain:
\begin{align}\label{aux510}
 \Delta t\|\gradG \mu^n\|^2_{\Gamma_n} &+\beta_s \|c^{n}-c^{n-1}\|^2_{\Gamma_n}
 + \beta_s\Delta t \int_{\Gamma_n} (c^{n}-c^{n-1})(\bu^n \cdot \nabla c^{n})\, ds \cl
  &+ \frac{\epsilon}2\left(\|\gradG c^{n}\|^2_{\Gamma_n} -\|\gradG c^{n-1}\|^2_{\Gamma_n}+ \|\gradG(c^{n}-c^{n-1})\|^2_{\Gamma_n}\right) \cl
   &+\epsilon\Delta t\int_{\Gamma_n}\gradG c^n \, \gradG (\bu^n \cdot \nabla c^{n}) \, ds\cl
  &= -\frac1{\epsilon}\int_{\Gamma_n} {f}_0'(c^{n}) \,(c^{n}-c^{n-1}) \, ds -\frac{\Delta t}{\epsilon}\int_{\Gamma_n} {f}_0'(c^{n}) \,(\bu^n \cdot \nabla c^{n}) \, ds,
\end{align}
where $\|\cdot\|_{\Gamma_n}$ stands for the $L^2(\Gamma_n)$ norm.
First, we treat the terms with ${f}_0'$ on the right-hand side.
Here and throughout the rest of this section, we use the fact that \eqref{eq:sys_CH3_im}
yields $\nabla c^n=\nabla_\Gamma c^n$, $\nabla {f}_0(c^{n})=\gradG {f}_0(c^{n})$, etc.
Using integration by parts, formula \eqref{int_parts}, and condition \eqref{Incompr}, we have:
\begin{equation}\label{aux517}
\int_{\Gamma_n} {f}_0'(c^{n}) \,(\bu^n \cdot \nabla c^{n}) \, ds=\int_{\Gamma_n} \bu^n \cdot\gradG {f}_0(c^{n})  \, ds = -\int_{\Gamma_n} (\div_\Gamma\bu^n){f}_0(c^{n})\, ds=0.
\end{equation}
With the help of trunceted Taylor expansion and ${f}_0''\ge-\frac14$, we get
\begin{align}\label{aux522}
-\int_{\Gamma_n} {f}_0'(c^n) \,(c^{n}-c^{n-1}) \, ds&= -\int_{\Gamma_n} ({f}_0(c^{n}) - {f}_0(c^{n-1})) \, ds \cl
&\quad -\frac12\int_{\Gamma_n} {f}_0''(\xi^{n}) \,|c^{n}-c^{n-1}|^2 \, ds\cl
&\le  -\int_{\Gamma_n} ({f}_0(c^{n}) - {f}_0(c^{n-1})) \, ds + \frac14\|c^{n}-c^{n-1}\|^2_{\Gamma_n}.
\end{align}
Next, we consider the terms on the left-hand side of \eqref{aux510}.
Using Cauchy--Schwarz inequality and estimate
$ab\le \frac1{2\delta}a^2 + \frac{\delta}2b^2$, $\forall~a,b\in \mathbb{R}$ with $\delta\in(0,+\infty)$,
we obtain:
\begin{equation}\label{aux536}
 \beta_s\Delta t \int_{\Gamma_n} (c^{n}-c^{n-1})(\bu^n \cdot \nabla c^{n})\, ds \le \frac{\beta_s}{2} \|c^{n}-c^{n-1}\|^2_{\Gamma_n} + C \frac{|\Delta t|^2}{\epsilon}\|\gradG c^{n}\|^2_{\Gamma_n},
 \end{equation}
 where  constant $C$ depends on $\mathbf{u}(\bx,t)$.
The last term on the left-hand side of  \eqref{aux510} is handled by the integration by parts
\begin{align}\label{aux554}
&\epsilon\Delta t\int_{\Gamma_n}\gradG c^{n} \, \gradG (\bu^n \cdot \nabla c^{n})\, ds
=\epsilon\Delta t\int_{\Gamma_n}\gradG c^{n} \cdot [\gradG \bu^n]\gradG c^{n}\, ds \cl
&\quad + \epsilon\Delta t\int_{\Gamma_n}\gradG c^{n} \,(\bu^n \cdot  \gradG ) \gradG c^{n}\, ds
= \epsilon\Delta t\int_{\Gamma_n}\gradG c^{n} \cdot [\gradG \bu^n]\gradG c^{n}\, ds \cl
& \quad \le  C\epsilon\Delta t\|\gradG c^{n}\|_{\Gamma^n}^2.
\end{align}
Combining \eqref{aux510}--\eqref{aux554}, for $\beta_s>\frac1{2\epsilon}$ we get
\begin{multline}\label{aux576}
 \frac{\epsilon}2\|\gradG c^{n}\|^2_{\Gamma_n} + \frac1{\epsilon}\int_{\Gamma_n} {f}_0(c^{n}) \, ds +\Delta t\|\gradG \mu^n\|^2_{\Gamma_n} \le \frac{\epsilon}2\|\gradG c^{n-1}\|^2_{\Gamma_n}+\frac1{\epsilon}\int_{\Gamma_n} {f}_0(c^{n-1})\, ds\\
     +C\Delta t(\epsilon+\frac{|\Delta t|}{\epsilon})\|\gradG c^{n}\|_{\Gamma^n}^2,
\end{multline}
where the constant $C$ depends on $\mathbf{u}(\bx,t)$.

Before applying the discrete Gronwall inequality to pass from \eqref{aux576} to a priori stability estimate, we need to relate
the quantities $\|\gradG c^{n-1}\|^2_{\Gamma_n}$ and $\int_{\Gamma_n} {f}_0(c^{n-1})\, ds$ to $\|\gradG c^{n-1}\|^2_{\Gamma_{n-1}}$ and $\int_{\Gamma_{n-1}} {f}_0(c^{n-1})\, ds$, respectively.
To this end, consider the closest point projection $\bp^{n-1}:
\cO_\delta(\Gamma_{n-1})\to \Gamma_{n-1}$. Since $c^{n-1}$ satisfies \eqref{eq:sys_CH3_num}, we can write $c^{n-1}(\bx)=c^{n-1}(\bp^{n-1}(\bx))$ in $\cO_\delta(\Gamma_{n-1})$. In particular, due to \eqref{ass1}
this representation holds  for $\bx\in\Gamma_n$.
For the surface measures on $\Gamma_{n-1}$ and $\Gamma_{n}$, we have \cite[Lemma 1]{lehrenfeld2018stabilized}:
\begin{equation}\label{aux594}
 \mu^n(\bx)ds^n(\bx)= ds^{n-1}(\bp^{n-1}(\bx)),\quad |1-\mu^n(\bx)|\le C \Delta t, \quad\bx\in\Gamma_{n},
\end{equation}
with some $C>0$ depending only on surface velocity $\bu$. Therefore, for $\Delta t$ small enough
we obtain
\begin{equation}\label{aux598}
\int_{\Gamma_n} {f}_0(c^{n-1})\, ds^n = \int_{\Gamma_{n-1}} {f}_0(c^{n-1})\, [\mu^n]^{-1} ds^{n-1}\le (1+C \Delta t)\int_{\Gamma_{n-1}} {f}_0(c^{n-1})\, ds^{n-1}.
\end{equation}
In turn, the surface gradients on $\Gamma_n$ and $\Gamma_{n-1}$ are related through (cf., e.g., \cite[section 2.3]{demlow2007adaptive})
\[
\nabla_{\Gamma_n} c^{n-1}(\bx)= \bP_{\Gamma_n}(\bx)(1-d(\bx)\bH(\bx))\nabla_{\Gamma_{n-1}} c^{n-1}(\bp^{n-1}(\bx)),\quad\bx\in\Gamma_{n},
\]
where $\bP_{\Gamma_n}$ is the orthogonal projector on the tangential space on $\Gamma_n$,
$d(\bx)$ a signed distance function for $\Gamma_{n-1}$, and $\bH(\bx)=\nabla^2 d(\bx)$ is the shape operator.  Our assumptions on $\Gamma(t)$ and its smooth evolution imply the uniform in $n$ bounds $\|\bH\|_{L^\infty(\cO_\delta(\Gamma_{n-1}))}\le C$ and
$\|d\|_{L^\infty(\Gamma_{n})}\le C\Delta t$. Together with \eqref{aux594}, this leads to estimate
\begin{equation}\label{aux608}
\|\gradG c^{n-1}\|^2_{\Gamma_n}\le (1+C \Delta t) \|\gradG c^{n-1}\|^2_{\Gamma_{n-1}}.
\end{equation}

By employing \eqref{aux598} and \eqref{aux608} in \eqref{aux576}, we are led to the following bound
\begin{align}\label{aux652}
& \frac{\epsilon}2\|\gradG c^{n}\|^2_{\Gamma_n} + \frac1{\epsilon}\int_{\Gamma_n} {f}_0(c^{n}) \, ds + \Delta t\|\gradG \mu^n\|^2_{\Gamma_n} \cl
& \quad  \le \left(1+C \Delta t\right) \left(\frac{\epsilon}2 \|\gradG c^{n-1}\|^2_{\Gamma_{n-1}}+\frac1{\epsilon}\int_{\Gamma_{n-1}} {f}_0(c^{n-1})\, ds\right) \cl
& \quad    +C\Delta t(\epsilon+\frac{\Delta t}{\epsilon})\|\gradG c^{n}\|_{\Gamma^n}^2.
\end{align}
We sum the inequalities \eqref{aux652} over $n=1,\dots,N$  to get
\begin{align}\label{aux659}
& \frac{\epsilon}2\|\gradG c^{N}\|^2_{\Gamma_N} + \frac1{\epsilon}\int_{\Gamma_N} {f}_0(c^{N}) \, ds +\Delta t\sum_{n=1}^{N}\|\gradG \mu^n\|^2_{\Gamma_n} \cl
& \quad  \le  \widetilde{C}(1+\frac{|\Delta t|}{\epsilon^{2}}) \sum_{n=0}^{N}\Delta t \left(  \frac{\epsilon}2 \|\gradG c^{n}\|^2_{\Gamma_{n}}+\frac1{\epsilon}\int_{\Gamma_{n}} {f}_0(c^{n})\, ds\right) \cl
& \quad   +\frac{\epsilon}2 \|\gradG c^{0}\|^2_{\Gamma_{0}}+\frac1{\epsilon}\int_{\Gamma_{0}} {f}_0(c^{0})\, ds.
\end{align}
Now assume that $\Delta t$ is small enough such that $1-\Delta t\widetilde{C}(1+\frac{|\Delta t|}{\epsilon^{2}})=\alpha>0$
and apply the discrete Gronwall inequality  to obtain
\begin{equation}\label{StabEstFinal}
 \frac{\epsilon}2\|\gradG c^{n}\|^2_{\Gamma_n} + \frac1{\epsilon}\int_{\Gamma_n} {f}_0(c^{n}) \, ds +\Delta t\sum_{k=1}^{n}\|\gradG \mu^k\|^2_{\Gamma_k} \le C({\rm data}),\quad \text{for all}~n=1,\dots,N,
\end{equation}
where the constant on the right-hand side depends on given problem data, i.e., on $\bu$, $\Gamma(0)$, $c^0$, but does not depend on $n$ and $\Delta t$. The constant $C({\rm data})$ is also uniformly bounded with respect to $\epsilon$ if $|\Delta t|\lesssim \epsilon$.

\subsection{Fully discrete method}\label{sec:fully_discrete}

The (quasi)-stationary system \eqref{eq:sys_CH1_num}--\eqref{eq:sys_CH3_num} involves only integrals
over $\Gamma^n$ and $\cO_\delta(\Gamma^n)$. This enables us to apply a surface finite element
method developed for a steady surface. Below we consider an unfitted finite element method,
known as the TraceFEM~\cite{olshanskii2017trace}, to discretize problem
\eqref{eq:sys_CH1_num}--\eqref{eq:sys_CH3_num} in space.

We start by assuming a shape regular triangulation  $\T_h$
of the bulk computational domain $\Omega\subset\R^3$.
At this point, we consider the same surface-independent mesh for all $t\in(0,T)$.
However, in practice dynamic refinement is possible and the
numerical examples in Section~\ref{sec:num_res} demonstrate its
feasibility and use. We will address this option later.
The domain formed by all the elements in $\T_h^\Gamma$ is denoted by $\OGamma$.
On $\T_h$ we consider a standard finite element space of continuous functions that are piecewise-polynomials
of degree $1$.
This bulk (volumetric) finite element space is denoted by $V_h^{\rm bulk}$:
\[
V_h^{\rm bulk}=\{v\in C(\Omega)\,:\, v\in P_1(T)~\text{for any}~T\in\T_h\}.
\]

\begin{figure}[h]
\centering
\begin{overpic}[width=.4\textwidth,grid=false]{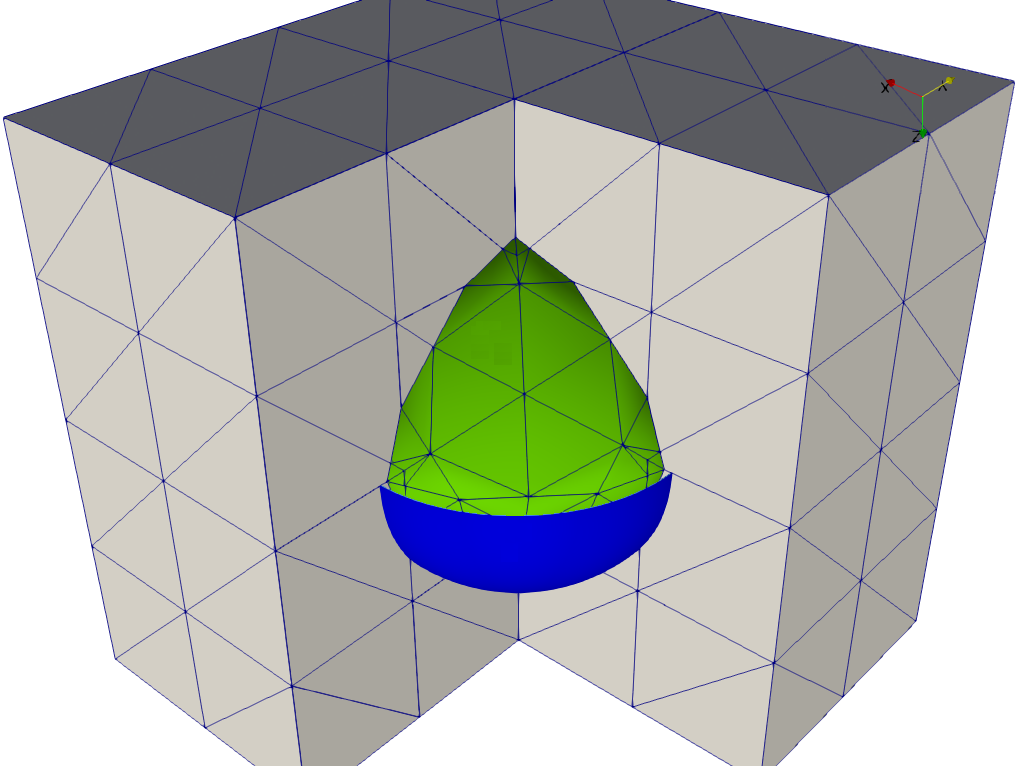}
         \put(26,43){\small{$\Omega$}}
         \put(50,35){\small{$\Gamma_h$}}
         \put(50,19){\small{\textcolor{white}{$\Gamma$}}}
\end{overpic}
 \caption{Example of bulk domain $\Omega$ with a given triangulation and sphere surface $\Gamma$ with the
  corresponding ``discrete'' surface $\Gamma_h$.}
  \label{fig:Gamma_h}
\end{figure}

Let the surface be given implicitly as the zero level of the level set function $\phi: \Omega\times[0,T]\to\R$:
\[
\Gamma(t)=\{\bx\in\R^3\,:\, \phi(\bx,t)=0\},
\]
with $|\nabla\phi|>0$ in $\cO(\Gs)$. Note that $\phi$ can be defined as the solution to the level-set equation~\cite{Sethian96AN}:
\begin{equation}\label{levelset}
\phi_t+\bu^{\rm bulk}\cdot\nabla\phi=0,
\end{equation}
where $\bu^{\rm bulk}$ is a given velocity field in $\Omega$ such that $\bu^{\rm bulk}\cdot\bn=\bu\cdot\bn$ on $\Gamma(t)$.
For the purpose of numerical integration, we approximate $\Gamma$ with a ``discrete'' surface $\Gamma_h$, which is defined as the zero level set of a $P_1$ Lagrangian interpolant for $\phi$ on the one time refined mesh:
\[
\Gamma_h=\{\bx\in\Omega\,:\, \phi_h(\bx)=0,\quad \phi_h:=I_{h/2}(\phi(\bx))\in V^{\rm bulk}_{h/2}\}, \quad \text{with}~\bn_h=
\frac{\nabla I^2_{h/2}(\phi)}{|\nabla I^2_{h/2}(\phi)|},
\]
where $I^2_{h/2}(\phi)$ is a $P_2$ nodal interpolant of the level set function. See Fig.~\ref{fig:Gamma_h}.
Note that $\bn_h$ is defined in a neighborhood of $\Gamma_h$.
The discrete surface $\Gamma_h$ satisfies the following geometric approximation properties:
\begin{equation}\label{geomApprox}
\mbox{dist}(\Gamma,\Gamma_h)=O(h^2)\quad\text{and}\quad\|\bn^e-\bn_h\|_{L^\infty(\Gamma_h)}=O(h).
\end{equation}
The analysis of the TraceFEM for scalar elliptic problem on a steady surface~\cite{olshanskii2017trace}
and parabolic problems on evolving surfaces~\cite{lehrenfeld2018stabilized} shows that the geometric consistency error \eqref{geomApprox} allows to prove optimal order convergence for $P^1$ elements (both in energy and $L^2$ norms). The full convergence analysis of TraceFEM for fourth order problems as
the surface Cahn-Hilliard equation \eqref{eq:CH_2} is an open problem that we shall address elsewhere.
When only the initial position of the surface and the bulk advection field $\bu^{\rm bulk}$ are given,
the evolution of $\Gamma$ can be recovered solving eq.~\eqref{levelset} numerically.
Then, the calculated solution $\phi_h$ defines $\Gamma_h$.

The numerical method approximates the solution to the surface
Cahn-Hilliard problem and its extension to a neighborhood of the surface.
In practice, at time  $t^n$ we extend the current solution to a narrow band of
$\Gamma^n_h$ such that $\Gamma^{n+1}_h$ stays inside this narrow band
and so all discrete quantities at $t^{n+1}$ are computable.
This narrow band  $\cO(\Gamma^n_h)$ is defined as the union  of tetrahedra within distance $\delta_n$ from
$\Gamma^n_h$:
\begin{equation} \label{defGamma}
 \overline{\cO(\Gamma^n_h)}=
\bigcup\left\{ \overline{T}\,:\,T\in\mathcal{T}_h\,: \mbox{dist}(T,\Gamma_h)\le\delta_n\right\}.
\end{equation}
To ensure
\begin{equation} \label{cond1}
\Gamma^{n+1}_h\subset \cO(\Gamma^n_h),
\end{equation}
 we set the minimum thickness of the extension layer to be
\begin{equation} \label{e:delta}
  \delta_n := c_\delta \Delta t \sup_{t\in(t_n,t_{n+1})}\| \bu\cdot\bn_h\|_{L^\infty(\Gamma_h)}
\end{equation}
where $c_\delta\geq 1 $ is an $O(1)$ mesh-independent constant.
If $\phi_h$ is an approximate distance function, then the condition $\mbox{dist}(S,\Gamma_h)\le\delta_n$ in \eqref{defGamma} can be replaced by $| \phi_h^n(\bx) |\le\delta_n$  for any vertex $\bx$ of $T$.

We define finite element spaces
\begin{equation}
V_h^n=\{v \in C(\cO(\Gamma^n_h))\,:\, v\in P_1(T), \forall T\subset \cO(\Gamma^n_h)\}.
\end{equation}
These spaces are the restrictions of the time-independent bulk space $V_h^{\rm bulk}$ on all tetrahedra from $\cO(\Gamma^n_h)$.

Our hybrid finite difference in time\,/\,finite element in space method is based on the semi-discrete formulation \eqref{eq:sys_CH1_weak2}--\eqref{eq:sys_CH3_weak2}. It reads:
Given $c_h^0 \in V_h^0$, for $n=1,\dots,N$ find  $c_h^n\in V_h^n$ and $\mu_h^n\in V_h^n$ satisfying
\begin{align}
&\int_{\Gamma^n_h}\rho^n \left(\frac{c_h^n-c_h^{n-1}}{\Delta t} + \bu^e \cdot \nabla c^n_h \right) \,v_h \, ds +  \int_{\Gamma^n_h} M(c^{n-1}_h) \nabla_{\Gamma_h^n} \mu^n_h \, \nabla_{\Gamma_h^n} v_h \, ds \cl
& \quad +  \rho_\mu \int_{\cO(\Gamma^n_h)}\frac{\partial \mu^n_h}{\partial \bn_h} \frac{\partial v_h}{\partial \bn_h} \, dx= 0, \label{eq:sys_CH1_h} \\
&\int_{\Gamma^n_h}  (\mu^n_h - \beta_s \left[c_h^n-c_h^{n-1}\right]  \,q_h \, ds  -  \epsilon\int_{\Gamma^n_h}\nabla_{\Gamma_h^n} c^n_h \, \nabla_{\Gamma_h^n} q_h \, ds \cl
& \quad -\rho_c \int_{\cO(\Gamma^n_h)}\frac{\partial c^n_h}{\partial \bn_h} \frac{\partial q_h}{\partial \bn_h} \, dx = \frac1{\epsilon}\int_{\Gamma^n_h} f_0'(c^{n-1}_h) \,q_h \, ds.\label{eq:sys_CH2_h}
\end{align}
for all $v_h,\,q_h\in V_h^n$. Here  $\bu^e(\bx)=\bu(\bp^n(\bx))$, i.e. the lifted data on $\Gamma^n_h$ from $\Gamma^n$. The first term in \eqref{eq:sys_CH1_h} is well-defined thanks to condition \eqref{cond1} with the index shifted $n\to n-1$.
In accordance to the analysis of the scalar advection-diffusion problems~\cite{lehrenfeld2018stabilized},
parameters $\rho_\mu$ and $\rho_c$ are set to be
\begin{equation}\label{param}
  \rho_\mu=\rho_c=h^{-1}.
\end{equation}

The role of the  term $\rho_c \int_{\cO(\Gamma^n_h)}\frac{\partial c^n_h}{\partial \bn_h} \frac{\partial q_h}{\partial \bn_h} \, dx$
in eq.~\eqref{eq:sys_CH2_h} is twofold. First, form
\[
 a( c^n_h,  q_h)=\epsilon\int_{\Gamma^n_h}\nabla_{\Gamma_h^n} c^n_h \, \nabla_{\Gamma_h^n} q_h \, ds+\rho_c \int_{\cO(\Gamma^n_h)}\frac{\partial c^n_h}{\partial \bn_h} \frac{\partial q_h}{\partial \bn_h} \, dx
\]
is positive definite on the narrow-band finite element space $V_h^n$, rather than only on the space of traces.
Therefore, at each time step we obtain a finite element solution $c^n_h$ defined in  $\cO(\Gamma^n_h)$.
This can be seen as an \emph{implicit extension procedure} of finite element solution from $\Gamma^n_h$ to the narrow band  $\cO(\Gamma^n_h)$. The same observation holds for $\mu^n_h$
and the corresponding term in eq.~\eqref{eq:sys_CH1_h}.
Furthermore, adding two volume terms (one for $c^n_h$ and another for $\mu^n_h$) 
makes the problem algebraically stable, i.e. the condition numbers of the resulting matrices are independent on how the surface cuts through the background mesh. Actually, the algebraic stabilization was the original motivation of introducing such volumetric terms in~\cite{burman2016cutb,grande2018analysis} for unfitted surface finite element methods.

\begin{remark}[Implementation]\rm For the realization of the method,
one uses the standard nodal basis functions for the bulk volumetric mesh $\T_h$.
At time step $n$, only the degrees of freedom of the tetrahedra in the narrow band $\cO(\Gamma^n_h)$
are active. Since $\Gamma_h^n$ is the zero level of the $P_1$ finite element
function $\phi_h^n$, the tetrahedra intersected by $\Gamma_h^n$ are those $T\in \cO(\Gamma^n_h)$
for which $\phi_h^n$ has a change in sign at different vertices.
Then, $\Gamma_h^n\cap T$ is either a triangle or a flat rectangle, which can be further divided in two triangles.
The vertexes of these triangles are immediately available from the nodal values of $\phi_h^n$ in  $T$.
So, it is straightforward to apply standard quadrature rules to compute the surface
integrals in \eqref{eq:sys_CH1_h}--\eqref{eq:sys_CH2_h}.
We note that on hyperrectangles the implicit representation of $\Gamma_h$ by the bilinear level-set function
can be treated with a marching cube \cite{lorensen1987marching} approximation.
\end{remark}

\begin{remark}[Higher order FEM]\label{rem:integration}\rm
Although this paper discusses only the finite element method based on  polynomials of degree 1,
higher order trace finite elements are possible. To gain full benefits of using higher order elements, one should  
ensure that the geometric and finite element interpolation errors are of the same order.
One way to achieve this, is to use  higher order polynomials to define the discrete level
set function $\phi_h$ that implicitly defines $\Gamma_h$.
However, it is a non-trivial task to obtain a parametrized representation of the zero level
of $\phi_h$ (for polynomial degree $\ge2$) which allows for a straightforward application of numerical quadrature rules.
The higher order case  requires special approaches for the construction of quadrature rules as discussed, for example, in  \cite{fries2015,muller2013highly,olshanskii2016numerical,saye2015hoquad,sudhakar2013quadrature}. Isoparametric unfitted finite element method is an another elegant and efficient way to simultaneously leverage the geometric and interpolation accuracy of TraceFEM; see  \cite{lehrenfeld2015cmame,grande2018analysis}.
\end{remark}

\begin{remark}[Methods based on surface PDEs extension]\rm
The finite element method \eqref{eq:sys_CH1_h}--\eqref{eq:sys_CH2_h} is different from a
surface unfitted FEM based on the extension of  a PDE from the surface to an ambient bulk
domain~\cite{BCOS01,Greer,deckelnick2009h,olshanskii2016narrow}.
Both methods avoid surface triangulation and remeshing (if the surface evolves).
However, in the  method  originally suggested in \cite{BCOS01} a PDE is first extended and
then solved in a one-dimension higher domain (the challenges
of such approach are discussed and partially addressed in  \cite{Greer,olshanskii2016narrow}), while the present
method extends only the solution, i.e. a function rather than a PDE is extended to the narrow band.
In the TraceFEM, the PDE is formulated on the original  surface.
\end{remark}

\section{Numerical examples}\label{sec:num_res}

We assess the model introduced in Sec.~\ref{s_cont} and the method presented in Sec.~\ref{s_FEM}
on a series of benchmark tests in Sec.~\ref{sec:validation}.
In Sec.~\ref{sec:merging}, we apply our approach to study spinodal decomposition and pattern formation on colliding surfaces.
Finally, in Sec.~\ref{sec:splitting} we consider the phase separation on a sphere splitting into two droplets.

For  the numerical tests, we will consider mobility $M$ as in \eqref{eq:M} with $\sigma=1$
(the only exception is test 2c, where we set $\sigma=16$) and free energy \eqref{eq:f0}.
Furthermore, we always use constant density $\rho=1$.
Although, there are no reasons to assume $\rho=\mbox{const}$ for
motions that violate the inextensibility condition (cf. \eqref{eq:evol_rho}),
we  expect that moderate smooth variations of $\rho$ would lead to minor changes
in lateral phase dynamics. In addition, $\rho=1$ is a fair assumption for the
assessment of the numerical method and the model.
We set $c_\delta=1$ in \eqref{e:delta} to define the extension strip.
Implementation of the method was been done using the FE package DROPS~\cite{DROPS}.

\subsection{Validation}\label{sec:validation}

To study the accuracy of the finite element method described in Sec.~\ref{s_FEM},
we consider two sets of benchmark tests. The first set features a sphere undergoing rigid body motion
and is presented in Sec.~\ref{sec:benchmark1}. The second set, illustrated in Sec.~\ref{sec:benchmark2},
involves an ellipsoid that dilates and shrinks.

\subsubsection{Sphere undergoing rigid body motion}\label{sec:benchmark1}
We consider a unit sphere, initially centered at the origin. During the time interval of interest $[0, 1]$,
the sphere is translated with translational velocity $\bv$ and rotated with angular
velocity $\boldsymbol{\omega}$.  We will report the numerical results for two tests:
\begin{itemize}
\item[-] \emph{Test 1a}: $\bv=(1,0,0)$ and $\boldsymbol{\omega}=(1,0,0)$. Its aim is to illustrate
the methodology.
\item[-] \emph{Test 1b}: $\bv=(10,0,0)$ and $\boldsymbol{\omega}=(10,0,0)$. Its aim is to
check the spatial and temporal accuracy.
\end{itemize}
In both test we take $\epsilon = 0.1$.

We  embed the sphere undergoing rigid body motion in outer domain
\[
\Omega= [-5/3,10/3]\times[-5/3,5/3]\times[-5/3,5/3].
\]
The initial triangulation $\T_{h_0}$ of $\Omega$ consists of 12 sub-cubes,
where each of the sub-cubes is further subdivided into 6 tetrahedra.
In addition, the mesh is refined towards the surface, and $\ell\in\mathbb{N}$ denotes the level of refinement, with the associated mesh size $h_\ell= 2^{-\ell-2}10/3$. 

\noindent{\bf Test 1a}.
To better illustrate the methodology, we first consider a simple test case.
We take the homogeneous equation Cahn--Hilliard equation.
We set $c_{x_2}$ as the initial solution, where
\begin{align}
\label{steady_CH}
c_{x_i}= \frac12\left(1+\tanh{}\left(\delta+\frac{x_i}{2\sqrt{2}\epsilon}\right)\right)\,.
\end{align}
Here and in the following, $\bx = (x_1, x_2, x_3)^T$ denotes a point in $\mathbb{R}^3$.
For this test we set $\delta=0$ (we will use $c_{x_i}$ with $\delta\neq0$ and $i\neq2$ later).

\begin{figure}[h]
  \centering
  \href{https://youtu.be/UylcNt9RnAw}{
    \subfloat[$t = 0$]{\includegraphics[width=.32\textwidth]{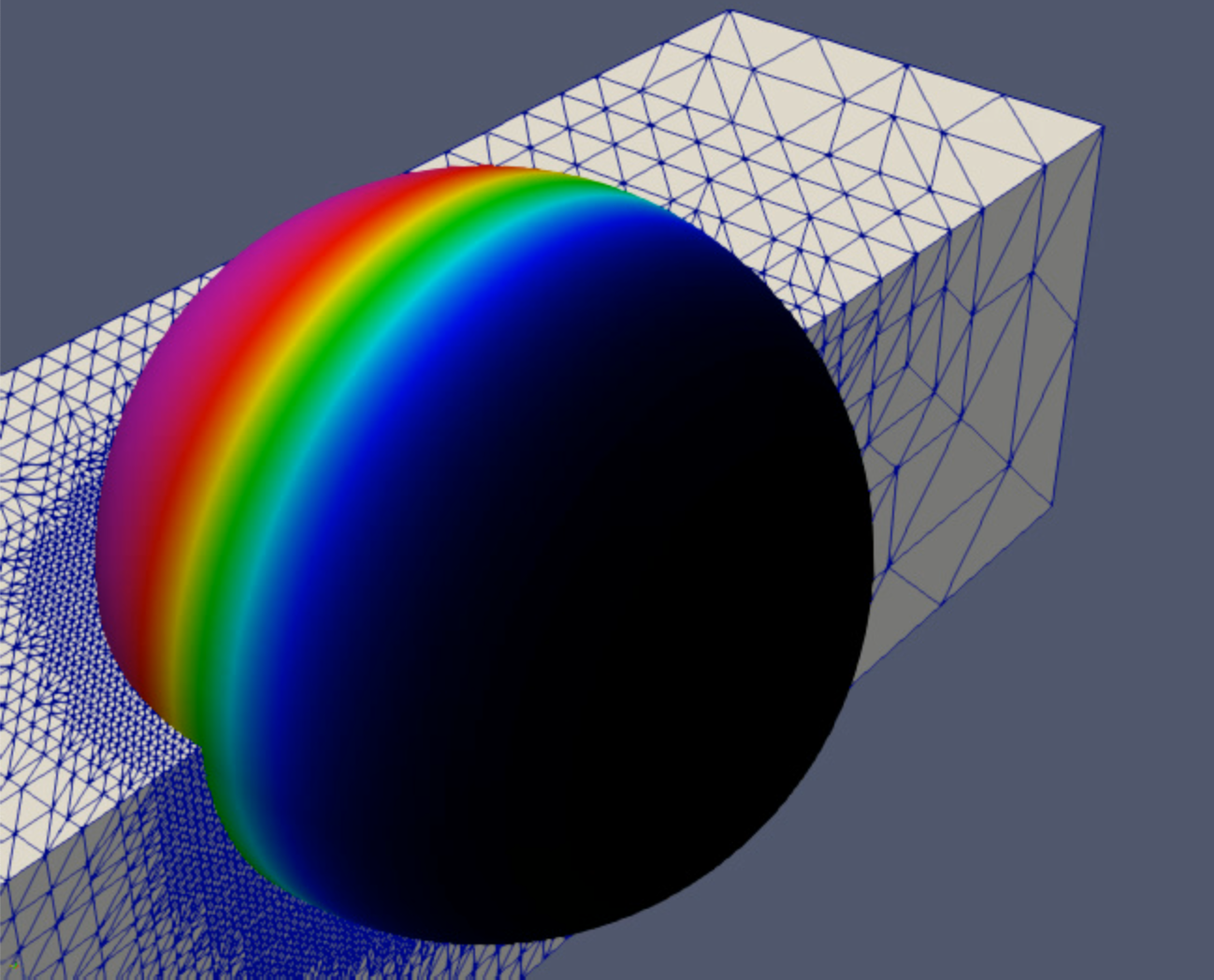}}~
  \subfloat[$t = 0.5$]{\includegraphics[width=.32\textwidth]{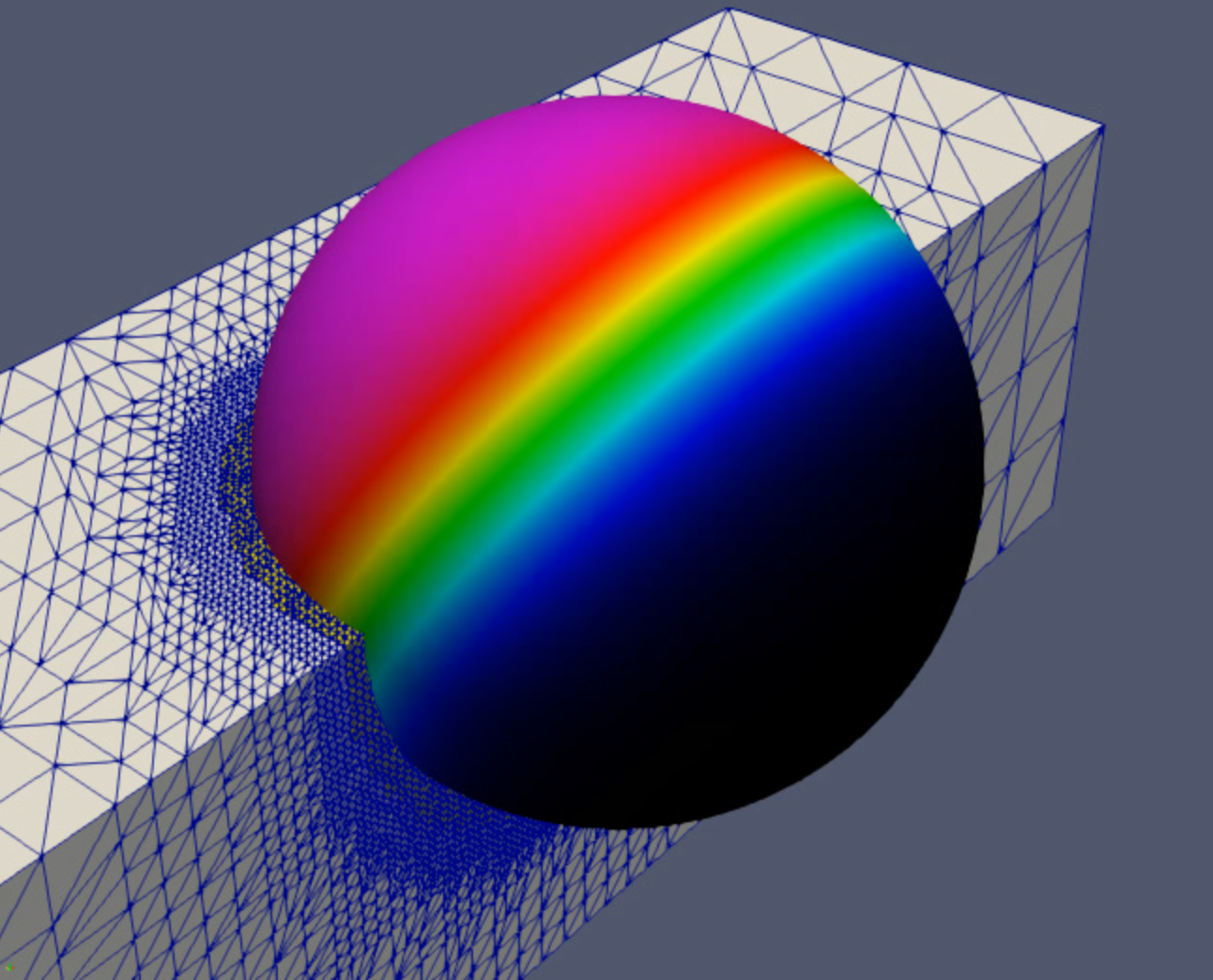}}~
   \subfloat[$t = 1$]{\includegraphics[width=.32\textwidth]{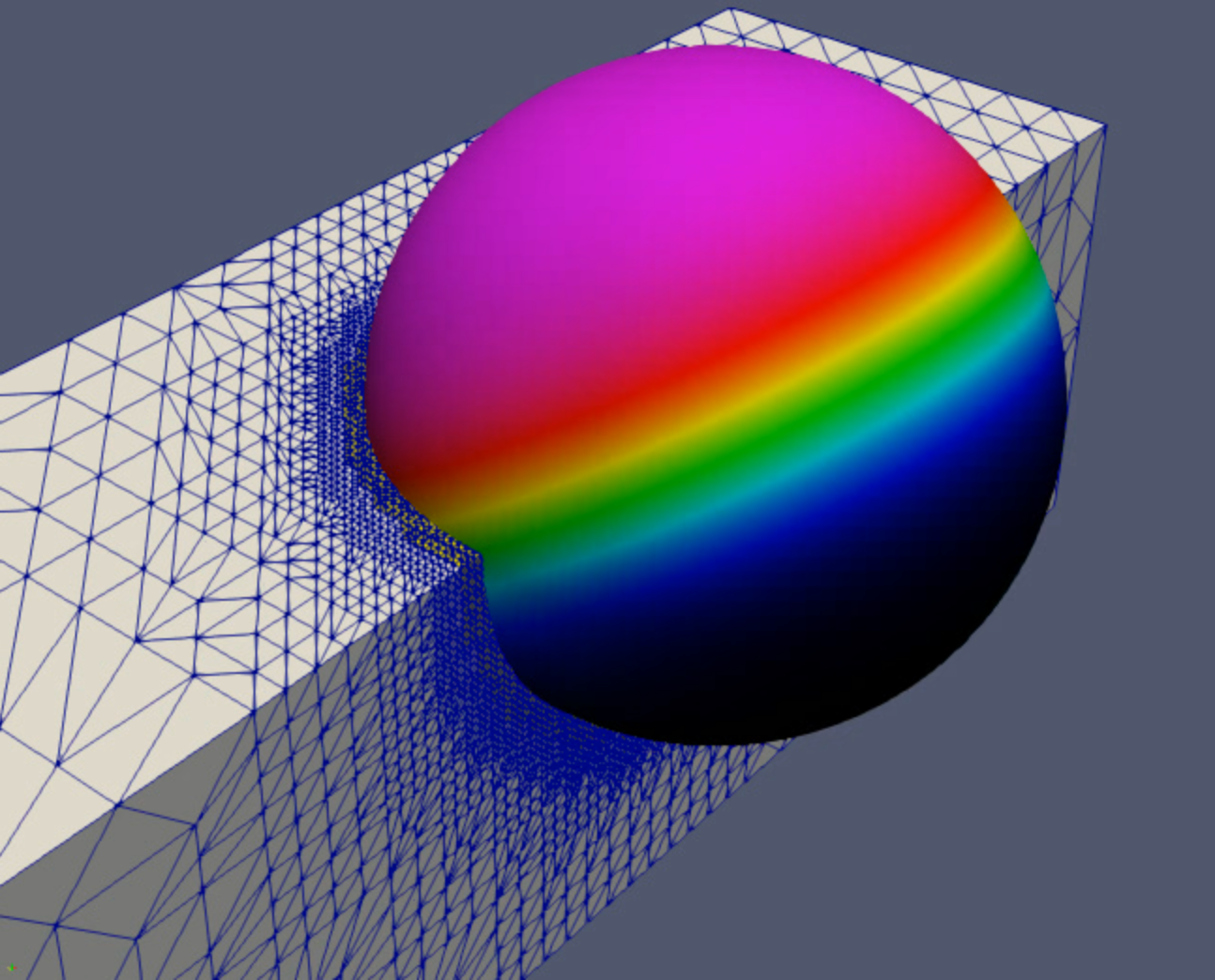}}
                                    }
  \hfill
  \caption{Test 1a: Solution computed with mesh $\ell = 5$ at times (a) $t = 0$, (b) $t = 0.5$,
  and (c) $t = 1$, together with a view of the bulk mesh. Click any picture above to run the full animation.}
  \label{fig:spiral_sphere}
\end{figure}

Fig.~\ref{fig:spiral_sphere} displays the solution computed with mesh $\ell = 5$
at times $t = 0, 0.5, 1$, 
together with a view of the bulk mesh. Notice the mesh refinement in the neighborhood of the surface.
This ensures that there are always enough elements to resolve the interface thickness $\epsilon$.

\noindent{\bf Test 1b}. We consider the following synthetic solution to the Cahn--Hilliard equation:
\begin{align*}
c^*(x_1,x_2)=\frac{1}{2}\left(Y(x_1,x_2)+ 1\right), \quad Y(x_1,x_2)=x_1 x_2, \quad t\in[0,0.1].
\end{align*}
The above $c^*$ is the exact solution to the non-homogeneous equation
$$
\dot{c}  -  \rho^{-1}\divG \left(M \gradG \left(\frac1{\epsilon}f_0' - \epsilon\Delta_\Gamma c\right)\right) = g.
$$
The exact chemical potential $\mu^*$ can be readily computed from eq.~\eqref{eq:sys_CH2}.
The non-zero right-hand side $g$ is calculated
so that the solution follows the motion of the sphere without any change.
To compute the right-hand side $g$ it is helpful to notice that $Y$ is a real spherical harmonic function.

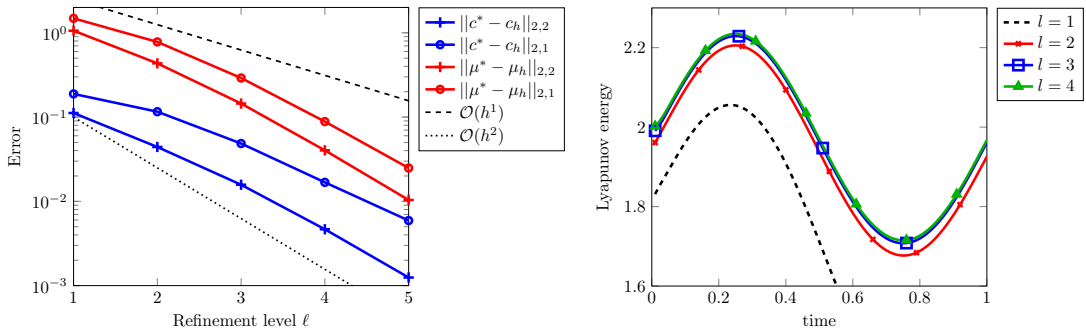
\begin{figure}[h]
	\begin{center}
\begin{minipage}{1.2\textwidth}
		\begin{tikzpicture}[scale=0.65]
		\def\vara{5}
		\def\varb{0.4}
		\begin{semilogyaxis}[ xlabel={Refinement level $\ell$},xmax=5,xmin=1, ylabel={Error}, ymin=1E-3
		, ymax=2, legend style={ cells={anchor=west}, legend pos=outer north east} ]

			\addplot+[blue,mark=+, solid,mark size=3pt,mark options={solid},line width=1.5pt] table[x=level, y=L2t_L2x_chi] {sec_spiral_h-tt.dat};
		
		\addplot+[blue,mark=o, solid,mark size=2pt,mark options={solid},line width=1.5pt] table[x=level, y=L2t_H1x_chi] {sec_spiral_h-tt.dat};
		
		\addplot+[red,mark=+, solid,mark size=3pt,mark options={solid},line width=1.5pt] table[x=level, y=L2t_L2x_omega] {sec_spiral_h-tt.dat};
		
		\addplot+[red,mark=o, solid,mark size=2pt,mark options={solid},line width=1.5pt] table[x=level, y=L2t_H1x_omega] {sec_spiral_h-tt.dat};
		
		\addplot[dashed,line width=1pt] coordinates { 
			(0,\vara) (1,\vara*0.5) (2,\vara*0.25) (3,\vara*0.125) (4,\vara*0.0625) (5,\vara*0.03125) (6,\vara*0.03125/2)
		};
		\addplot[dotted,line width=1pt] coordinates { 
			(0,\varb) (1,\varb*0.5*0.5) (2,\varb*0.25*0.25) (3,\varb*0.125*0.125) (4,\varb*0.0625*0.0625) (5,\varb*0.03125*0.03125) (6,\varb*0.03125*0.03125*0.5*0.5)
		};
		\legend{
			$|| c^* - c_h ||_{2,2}$,
			$|| c^* - c_h ||_{2,1}$,
			$|| \mu^* - \mu_h ||_{2,2}$,
			$|| \mu^* - \mu_h ||_{2,1}$,
			$\mathcal{O}(h^1)$ ,
			$\mathcal{O}(h^{2})$}
		\end{semilogyaxis}
		\end{tikzpicture}
		\def\vara{10}
		\def\varb{5}
		\begin{tikzpicture}[scale=0.65]
		\begin{axis}[ xlabel={time}, ylabel={Lyapunov energy}, ymin=1.6, xmax=1,xmin=0,
		, ymax=2.3, legend style={ cells={anchor=west}, legend pos=outer north east} ]
		
			\addplot+[black,mark=t,mark repeat={13},dashed,mark size=2pt,line width=1.5pt] table[x=Time, y=Lyapunov_energy] {l=1_N=100_eps=1e-1.txt};

			\addplot+[red,mark=x,mark repeat={13},solid,mark size=2pt,line width=1.5pt] table[x=Time, y=Lyapunov_energy] {l=2_N=100_eps=1e-1.txt};
		
		\addplot+[blue,mark=square,mark repeat={25},solid,mark size=3pt,line width=1.5pt] table[x=Time, y=Lyapunov_energy] {l=3_N=100_eps=1e-1.txt};
		
		\addplot+[black!30!green,mark=triangle,mark repeat={15},solid,mark size=2.5pt,line width=1.5pt] table[x=Time, y=Lyapunov_energy] {l=4_N=100_eps=1e-1.txt};

		\legend{
			$l=1$,
				$l=2$,
					$l=3$,
			$l=4$
		}
		\end{axis}
		\end{tikzpicture}
\end{minipage}
		\caption{Left plot: Discrete $L_2(0,T; L^2(\Gamma_h))$ norm \eqref{eq:22norm} and discrete $L_2(0,T; H^1(\Gamma_h))$ norm \eqref{eq:21norm} of the error for order parameter $c$ (blue lines) and chemical potential $\mu$ (red lines) in Test 1b plotted against the refinement level $\ell$ with time step $\Delta t=4^{1-\ell}/10$. Right plot: Convergence of the discrete Lyapunov energy defined in \eqref{eq:Lyapunov_E} with mesh refining for Test 2a.
		}
		\label{spiral_secharm_convergence} \label{fig:E_CH1}
	\end{center}
\end{figure}

We report in Fig.~\ref{spiral_secharm_convergence} (left) the discrete  $L_2(0,T; L^2(\Gamma_h))$ norm:
\begin{align}\label{eq:22norm}
\| c^* - c_h \|_{2,2}= \left(\frac{1}{T}\sum_k\, \Delta{}t\|c^*(t_k) -c_h(t_k)\|^2_{L^2(\Gamma_h)}\right)^{1/2}
\end{align}
and the discrete  $L_2(0,T; H^1(\Gamma_h))$ norm:
\begin{align}\label{eq:21norm}
\| c^* - c_h \|_{2,1}= \left(\frac{1}{T}\sum_k\, \Delta{}t\|\gradG(c^*(t_k) -c_h(t_k))\|^2_{L^2(\Gamma_h)}\right)^{1/2}
\end{align}
of the error for order parameter $c$ (blue lines) and chemical potential $\mu$ (red lines)
plotted against the refinement level $\ell$. The time step was refined together with the mesh size according to $\Delta t=4^{1-\ell}/10$.
All the norms reported in Fig.~\ref{spiral_secharm_convergence} (left) are
computed on the approximate surface $\Gamma_h$, where $c^*$ and $\mu^*$
were defined through their normal extensions from $\Gamma$.
The observed second order convergence in the $L_2(0,T; L^2(\Gamma_h))$ for both concentration and chemical potential as well as the first order convergence in the $L_2(0,T; H^1(\Gamma_h))$ norm for the concentration are consistent with the well known error analysis of finite element methods for the Cahn--Hilliard equation~\cite{elliott1992error,elliott2015evolving}.  For the error in chemical potential we also observe the almost second order of convergence in $L_2(0,T; H^1(\Gamma_h))$ norm. We do not have an explanation for this apparent super-convergence.

\subsubsection{Oscillating ellipsoid}\label{sec:benchmark2}

The second series of tests is inspired by a numerical experiment in \cite{elliott2015evolving}.
We consider time-dependent surface $\Gamma(t) = \{ \bx \in \mathbb{R}^3: \phi(\bx, t) = 0 \}$ with
$$
\phi(\bx, t) = \left(\frac{x_1}{a(t)}\right)^2+x_2^2+x_3^2 - 1\,,\quad a(t)=1+0.2\sin(2\pi{}kt),
$$
which is an ellipsoid centered at origin with variable length for the principal axis aligned
with the $x_1$-axis. We embed this evolving surface in domain $\Omega = [-5/3,5/3]^3$.
As initial solution to the Cahn--Hilliard equation, we take a small perturbation
about 0.5 given by:
$$
c_0=0.5+0.05\cos(2\pi{}x)\cos(2\pi{}y)\cos(2\pi{}z).
$$
We study the evolution of the discrete Lyapunov energy:
\begin{equation}\label{eq:Lyapunov_E}
E^L_h(c_h) = \int_{\Gamma_h} f(c_h) ds = \int_{\Gamma_h} \left( \frac1{\epsilon}f_0(c_h) + \frac{1}{2} \epsilon | \gradG c_h |^2 \right)ds
\end{equation}
for three datasets:
\begin{itemize}
\item[-] \emph{Test 2a}: $k=1$, $\sigma= 1$, $\epsilon=0.1$, $\Delta{}t=0.01$.
\item[-] \emph{Test 2b}: $k=1$, $\sigma=1$, $\epsilon=0.01$, $\Delta{}t=0.01$.
\item[-]  \emph{Test 2c}: $k=5$, $\sigma=16$, $\epsilon=0.1$, $\Delta{}t=0.001$.
\end{itemize}
For all the tests, we set $T=1$.
To mesh the domain we follow the same procedure used for the first sets of tests, i.e.~we divide
$\Omega$ into sub-cubes that are further divided into tetrahedra.
We consider different levels of refinement $\ell$ associated with
mesh size $h=2^{-\ell-2}10/3$.

\begin{figure}[h]
\begin{center}
\begin{minipage}{1.2\textwidth}
		\def\vara{10}
		\def\varb{5}
\begin{tikzpicture}[scale=0.7]
		\begin{axis}[ xlabel={time}, ylabel={Lyapunov energy}, ymin=0.1, xmax=1,xmin=0,
		, ymax=10, legend style={ cells={anchor=west}, legend pos=outer north east} ]
		
		\addplot+[black,mark=t,mark repeat={13},dotted,mark size=2pt,line width=1.5pt] table[x=Time, y=Lyapunov_energy] {l=1_N=100_eps=1e-2.txt};

		\addplot+[black,mark=t,mark repeat={23},dashed,mark size=2pt,line width=1.5pt] table[x=Time, y=Lyapunov_energy] {l=2_N=100_eps=1e-2.txt};
		
		\addplot+[red,mark=x,mark repeat={13},solid,mark size=3pt,line width=1.5pt] table[x=Time, y=Lyapunov_energy] {l=3_N=100_eps=1e-2.txt};
		
		\addplot+[black!30!blue,mark=o,mark repeat={10},solid,mark size=2pt,line width=1.5pt] table[x=Time, y=Lyapunov_energy] {l=4_N=100_eps=1e-2.txt};
		
			\addplot+[black!30!green,mark=square,mark repeat={27},solid,mark size=3pt,line width=1.5pt] table[x=Time, y=Lyapunov_energy] {l=5_N=100_eps=1e-2.txt};
		
		\addplot+[black!30!green,mark=triangle,mark repeat={10},solid,mark size=2pt,line width=1.5pt] table[x=Time, y=Lyapunov_energy] {l=6_N=100_eps=1e-2.txt};
		\legend{
			$l=1$,
			$l=2$,
			$l=3$,
			$l=4$,
			$l=5$,
			$l=6$
		}
		\end{axis}
		\end{tikzpicture}
	\begin{tikzpicture}[scale=0.7]
		\begin{axis}[ xlabel={time}, ylabel={Lyapunov energy}, ymin=0.5, xmax=1,xmin=0,
		, ymax=2.5, legend style={ cells={anchor=west}, legend pos=outer north east} ]
		
		\addplot+[black,mark=t,mark repeat={13},dotted,mark size=2pt,line width=1.5pt] table[x=Time, y=Lyapunov_energy] {l=1_N=1000_eps=1e-1.txt};
		
		\addplot+[black,mark=t,mark repeat={85},dashed,mark size=2pt,line width=1.5pt] table[x=Time, y=Lyapunov_energy] {l=2_N=1000_eps=1e-1.txt};
		
		\addplot+[red,mark=x,mark repeat={55},solid,mark size=3pt,line width=1.5pt] table[x=Time, y=Lyapunov_energy] {l=3_N=1000_eps=1e-1.txt};
		
		\addplot+[black!30!blue,mark=o,mark repeat={70},solid,mark size=2pt,line width=1.5pt] table[x=Time, y=Lyapunov_energy] {l=4_N=1000_eps=1e-1.txt};
		
		\addplot+[black!30!green,mark=square,mark repeat={100},solid,mark size=2.5pt,line width=1.5pt] table[x=Time, y=Lyapunov_energy] {l=5_N=1000_eps=1e-1.txt};
%
		\legend{
			$l=1$,
			$l=2$,
			$l=3$,
			$l=4$,
			$l=5$
				}
		\end{axis}
		\end{tikzpicture}
\end{minipage}
		\caption{Evolution of the discrete Lyapunov energy functional \eqref{eq:Lyapunov_E}. Left: Test 2b; Right: Test 2c.
		For test 2b (resp., 2c), we see mesh convergence for levels 5 and 6 (resp., 4 and 5).}
		\label{fig:E_CH2} \label{fig:E_CH3}
	\end{center}
\end{figure}
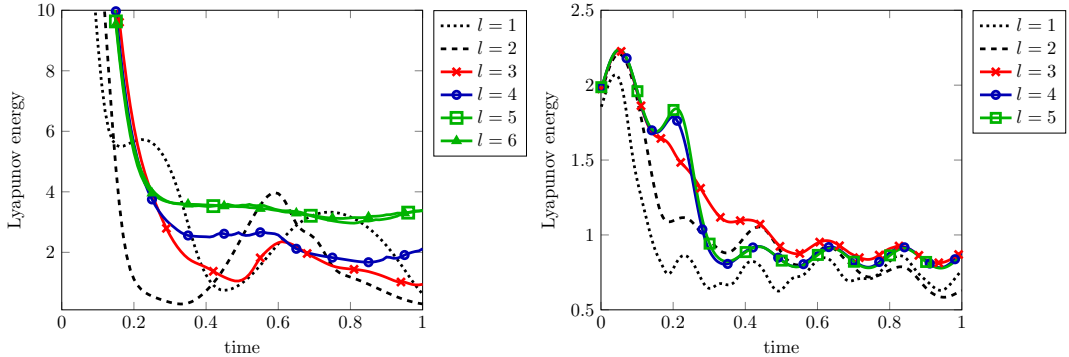

\noindent{\bf Test 2a.} We present the evolution of the discrete Lyapunov energy in Fig.~\ref{fig:E_CH1} (right).
We see that mesh $\ell = 1$ is too coarse to correctly capture the evolution of $E^L_h(c_h)$.
The results obtained with mesh $\ell = 2$ are very close to results obtained with meshes $\ell = 3,4$,
which are almost superimposed for the entire time interval under consideration.
Fig.~\ref{fig:E_CH1} (right) seems to suggest that that the solution converges to a time
periodic solution. In order to verify that, we would need to run the simulation
for a longer period of time. Instead, we prefer to decrease the value of $\epsilon$
in the next test. We recall that $\epsilon$
is a crucial modeling parameter, since it defines the thickness of the layer where
phase transition takes place and also defines the intrinsic time scale. Thus,
the case of smaller $\epsilon$ is numerically challenging.

\begin{figure}[ht]
\begin{center}
\href{https://youtu.be/-g0Y0UY9nm8}{
\begin{overpic}[height=.26\textwidth, viewport=105 130 555 520, clip,grid=false]{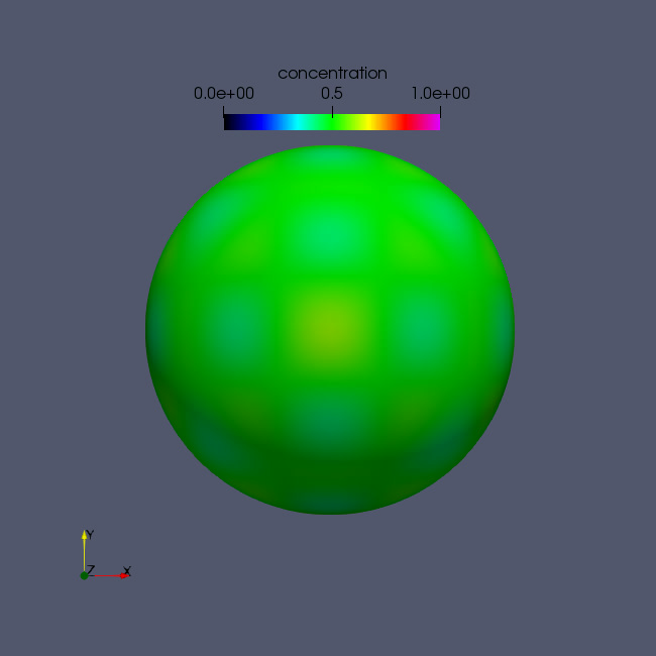}
         \put(40,90){\small{$t = 0$}}
\end{overpic}
\begin{overpic}[height=.26\textwidth, viewport=105 130 555 520, clip,grid=false]{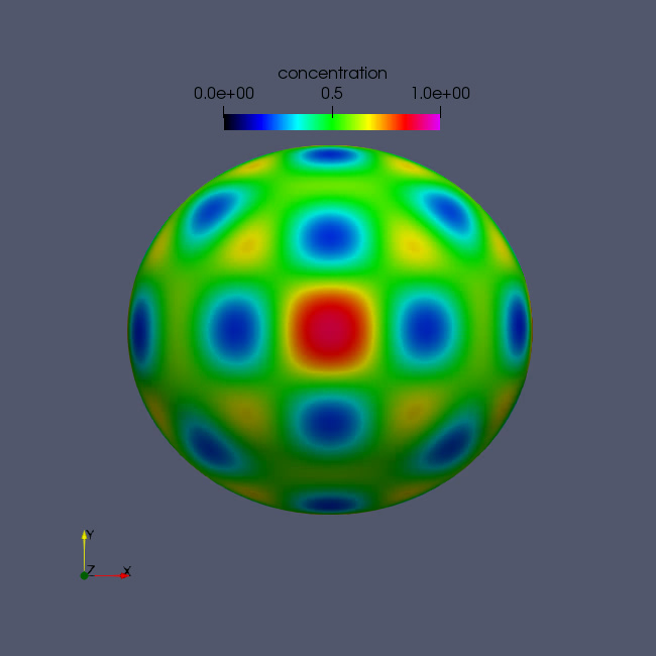}
         \put(40,90){\small{$t = 0.1$}}
\end{overpic}
\begin{overpic}[height=.26\textwidth, viewport=105 130 555 520, clip,grid=false]{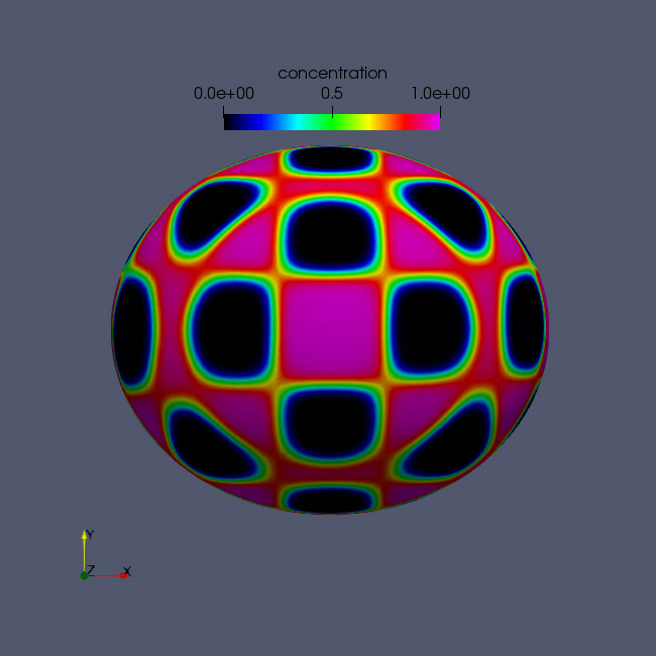}
         \put(40,90){\small{$t = 0.2$}}
\end{overpic}
}
\vskip .4cm
\href{https://youtu.be/-g0Y0UY9nm8}{
\begin{overpic}[height=.26\textwidth, viewport=105 130 555 520, clip,grid=false]{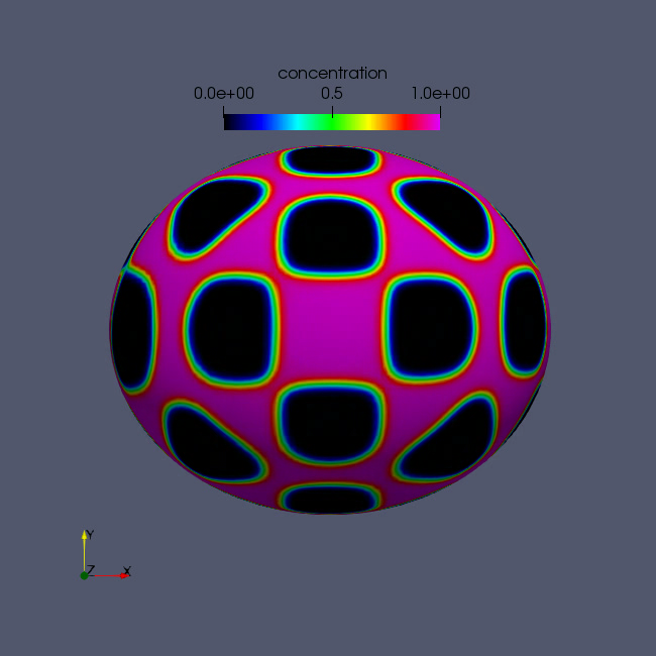}
         \put(40,90){\small{$t = 0.3$}}
\end{overpic}
\begin{overpic}[height=.26\textwidth, viewport=105 130 555 520, clip,grid=false]{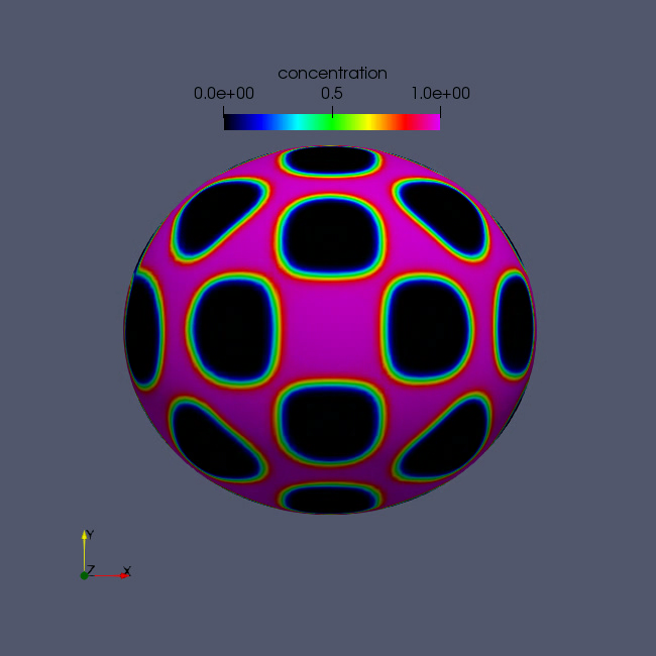}
         \put(40,90){\small{$t = 0.4$}}
\end{overpic}
\begin{overpic}[height=.26\textwidth, viewport=105 130 555 520, clip,grid=false]{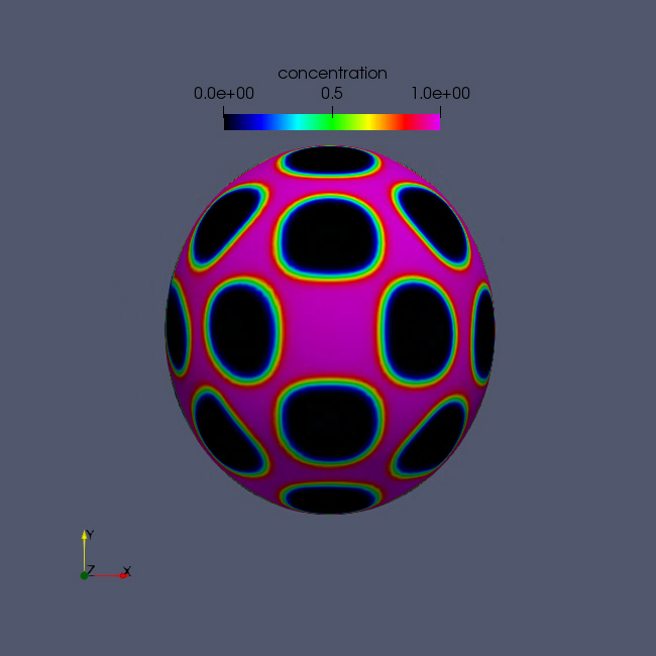}
         \put(40,90){\small{$t = 0.6$}}
\end{overpic}
}
\vskip .4cm
\hskip -1.3cm
\href{https://youtu.be/-g0Y0UY9nm8}{
\begin{overpic}[height=.26\textwidth, viewport=105 130 555 520, clip,grid=false]{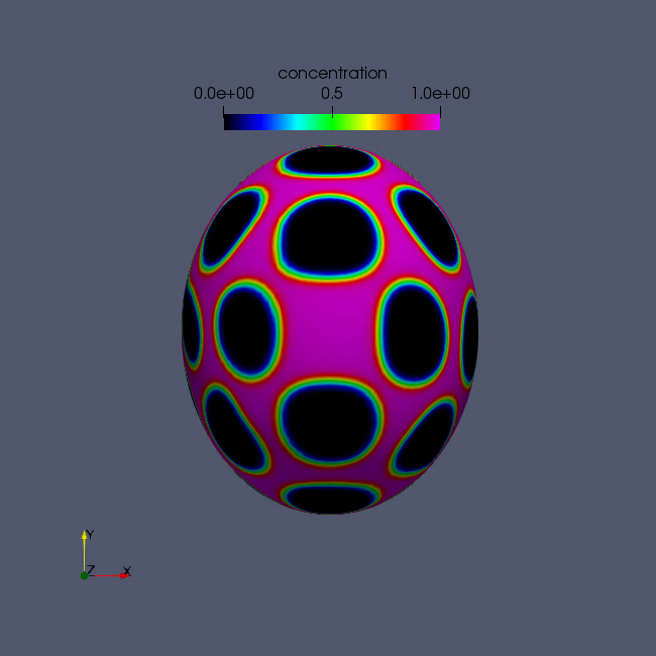}
         \put(40,90){\small{$t = 0.8$}}
\end{overpic}
\begin{overpic}[height=.26\textwidth, viewport=105 130 555 520, clip,grid=false]{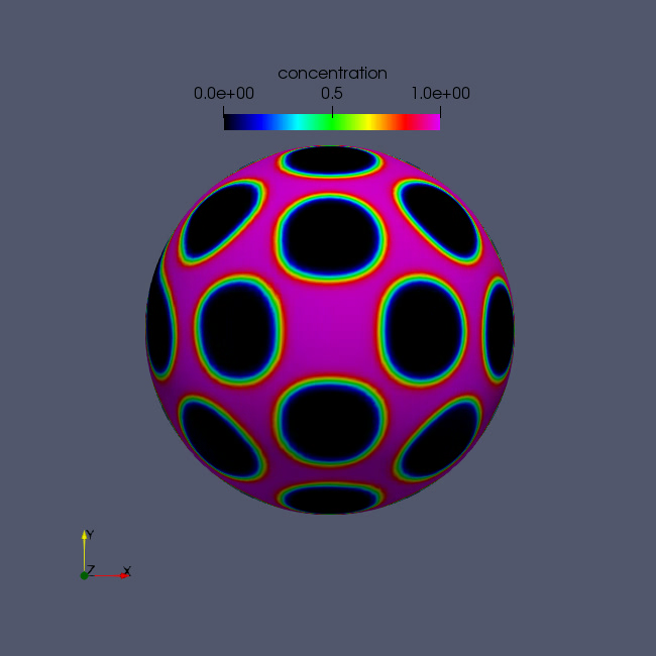}
         \put(40,90){\small{$t = 1$}}
\end{overpic}
\hspace{1cm}
\begin{overpic}[height=.26\textwidth,grid=false]{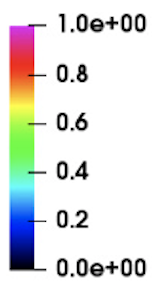}
\end{overpic}
}
\end{center}
\caption{\label{fig:test2b}
Test 2b:
Evolution of the numerical solution of the Cahn--Hilliard equation computed with mesh $\ell =6$ for $t \in [0, 1]$.
View: $x_1 x_2$-plane. Click any picture above to run the full animation.}
\end{figure}

\noindent{\bf Test 2b.}
We keep all the model and discretization parameters as in test 2a, but we reduce the value
of $\epsilon$ to $0.01$.
Fig.~\ref{fig:E_CH2} (left)  shows the evolution of $E^L_h(c_h)$ over time.
Since this test is more challenging, to observe convergence for the discrete Lyapunov energy
it takes a higher level of refinement than test 2a.
In fact, from Fig.~\ref{fig:E_CH2} (left) we see that the results computed with mesh $\ell = 4$
are still quite far from the results obtained with mesh $\ell = 5$, which are almost
superimposed to the results for mesh $\ell = 6$. This indicates that both meshes
$\ell = 5$ and $\ell = 6$ are sufficiently refined for  $\epsilon = 0.01$.
Notice that the discrete Lyapunov energy computed with meshes $\ell = 5, 6$
decreases rapidly till around $t = 0.3$ and then it gently decreases, although not
monotonically. Moreover, we note that the energy may have locally increase driven by the evolution of the domain. This phenomenon will be even more pronounced in test 2c, which is presented next .

%
%
%
%
%
%
%

We show the evolution of the solution computed with mesh $\ell = 6$
in Fig.~\ref{fig:test2b}. Around $t = 0.3$, we observe one big pink domain (i.e., $c = 1$)
and several small black domains (i.e., $c = 0$). For $t > 0.3$, the aspect of the solution
(one big pink domain and several  smaller black domains) remains almost unchanged
as the surface gets dilated and shrunk.

\begin{figure}[ht]
\begin{center}
\begin{overpic}[height=.27\textwidth, viewport=105 130 555 577, clip,grid=false]{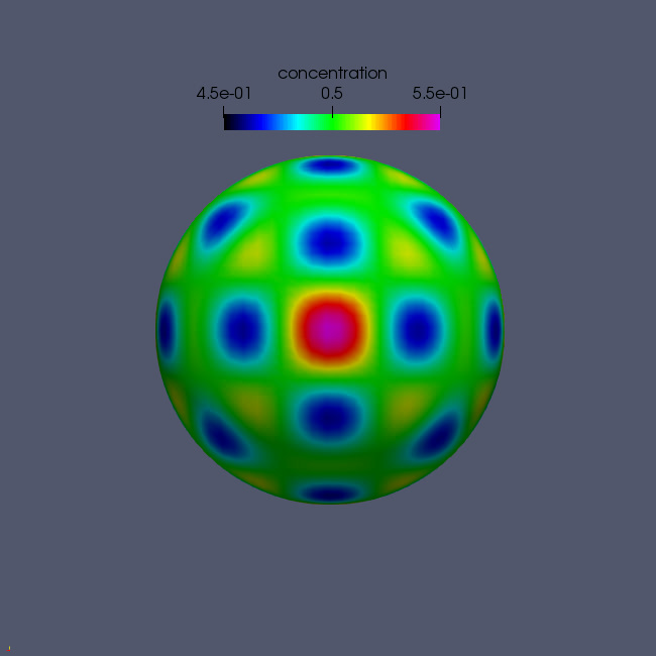}
         \put(40,102){\small{$t = 0$}}
\end{overpic}
\begin{overpic}[height=.27\textwidth, viewport=105 130 555 577, clip,grid=false]{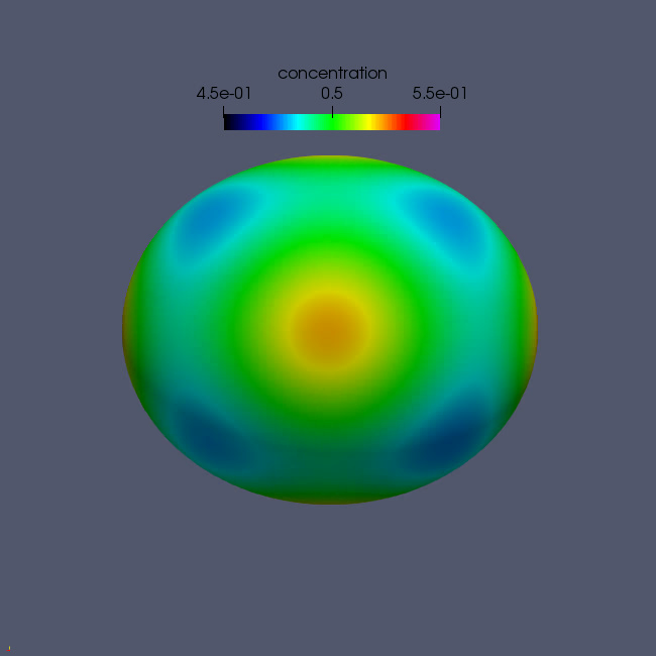}
         \put(40,102){\small{$t = 0.05$}}
\end{overpic}
\begin{overpic}[height=.27\textwidth, viewport=105 130 555 577, clip,grid=false]{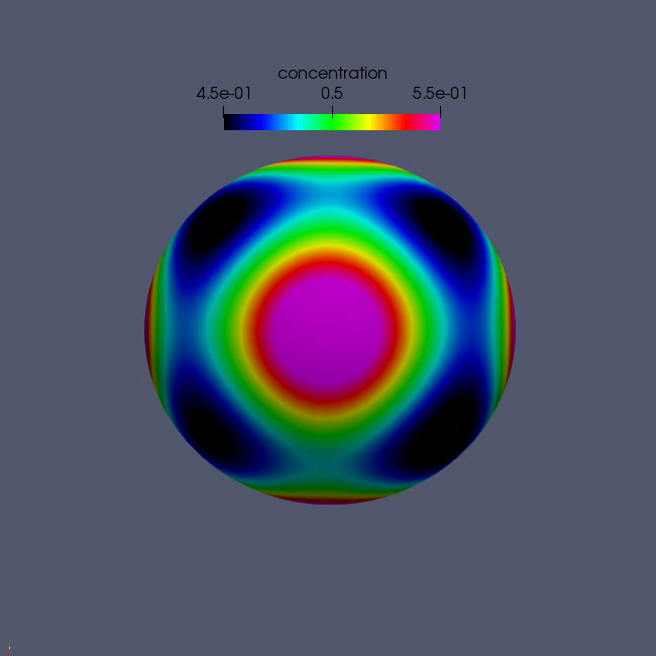}
         \put(40,102){\small{$t = 0.1$}}
\end{overpic}
\vskip .4cm
\begin{overpic}[height=.27\textwidth, viewport=105 130 555 577, clip,grid=false]{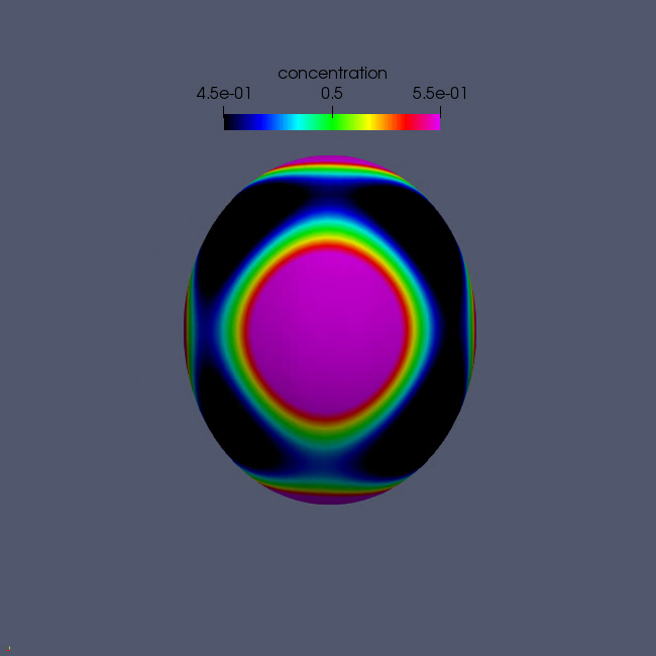}
         \put(40,102){\small{$t = 0.15$}}
\end{overpic}
\begin{overpic}[height=.27\textwidth, viewport=105 130 555 577, clip,grid=false]{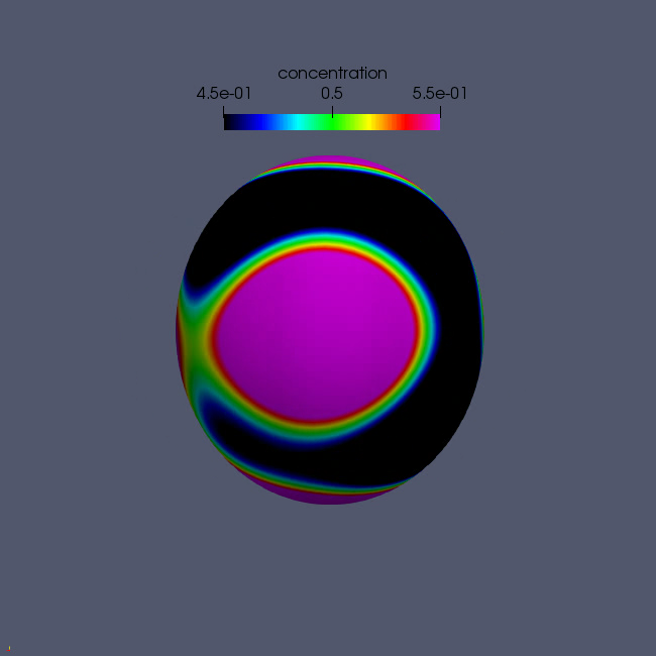}
         \put(40,102){\small{$t = 0.2$}}
\end{overpic}
\begin{overpic}[height=.27\textwidth, viewport=105 130 555 577, clip,grid=false]{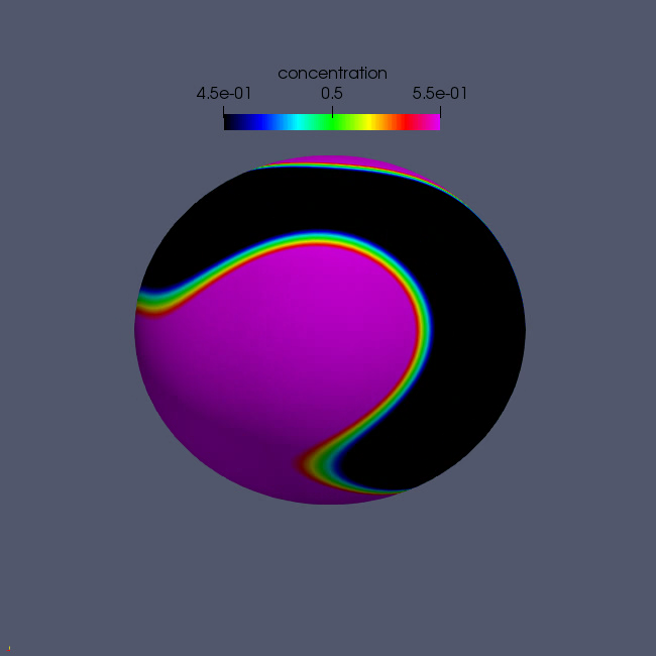}
        \put(40,102){\small{$t = 0.25$}}
\end{overpic}
\vskip .4cm
\begin{overpic}[height=.27\textwidth, viewport=105 130 555 577, clip,grid=false]{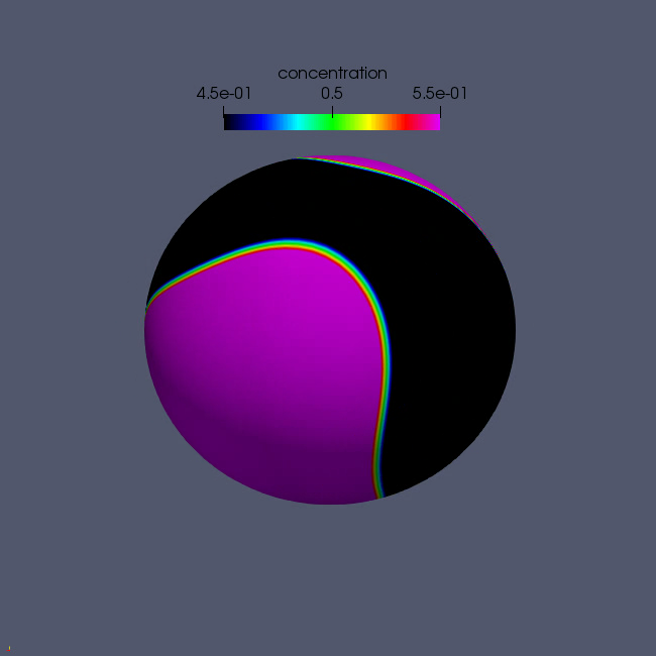}
         \put(40,102){\small{$t = 0.3$}}
\end{overpic}
\begin{overpic}[height=.27\textwidth, viewport=105 130 555 577, clip,grid=false]{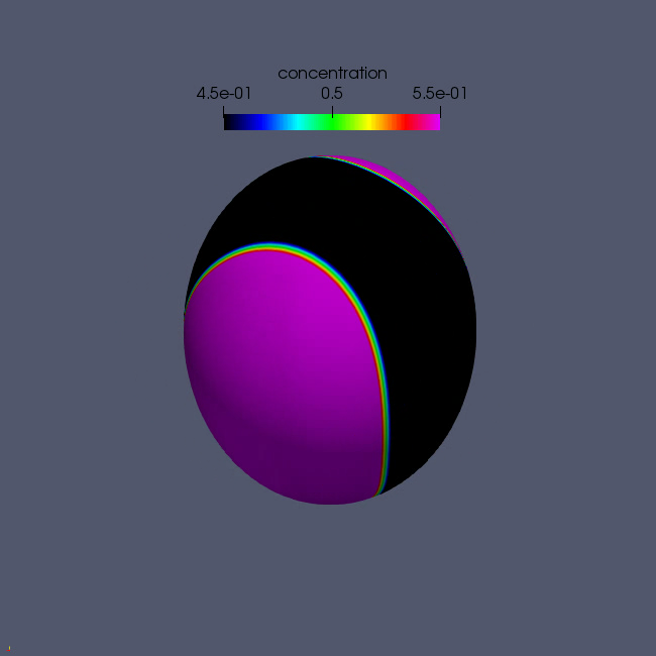}
         \put(40,102){\small{$t = 0.35$}}
\end{overpic}
\begin{overpic}[height=.27\textwidth, viewport=105 130 555 577, clip,grid=false]{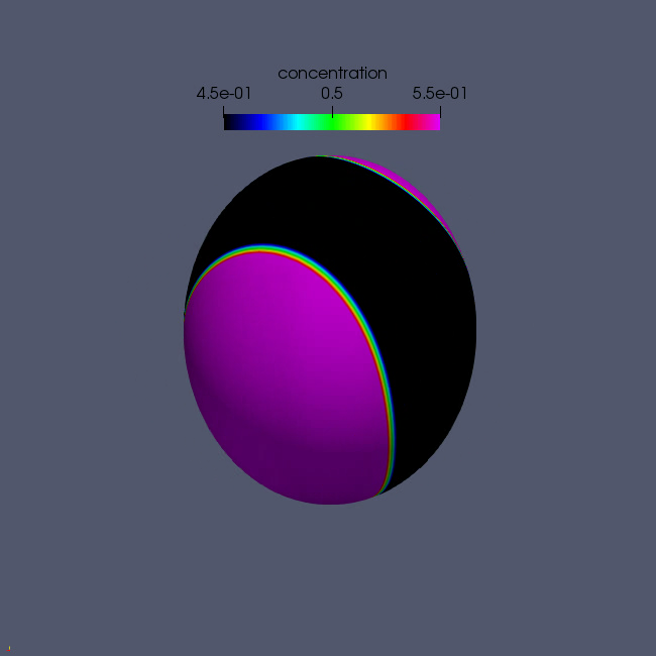}
         \put(40,102){\small{$t = 0.4$}}
\end{overpic}
\vskip .4cm
\begin{overpic}[height=.27\textwidth, viewport=105 130 555 577, clip,grid=false]{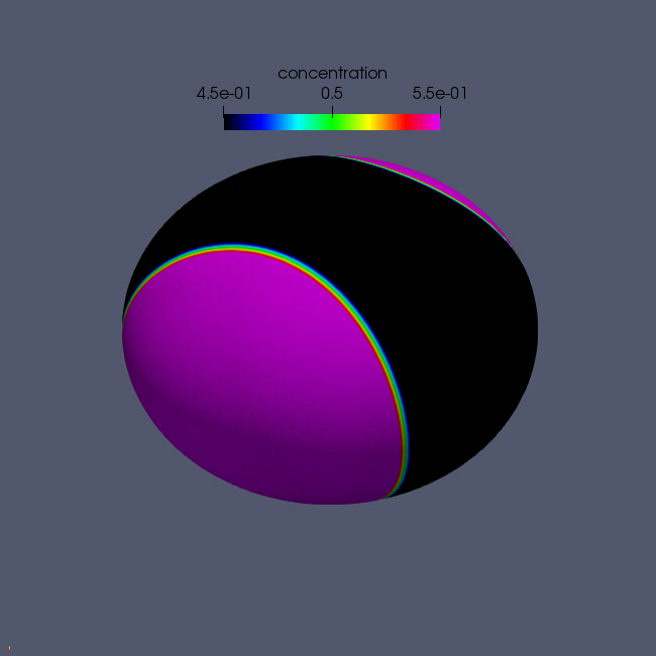}
         \put(40,102){\small{$t = 0.45$}}
\end{overpic}
\begin{overpic}[height=.27\textwidth, viewport=105 130 555 577, clip,grid=false]{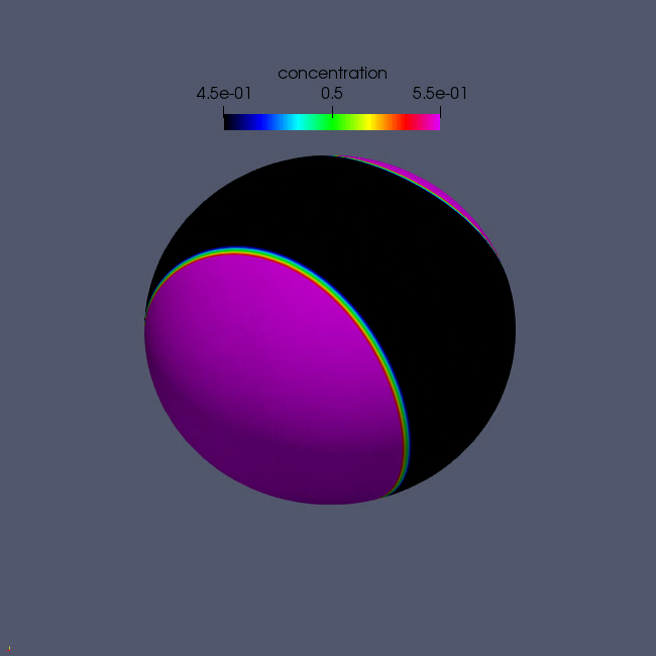}
         \put(40,102){\small{$t = 0.5$}}
\end{overpic}
\begin{overpic}[height=.27\textwidth, viewport=105 130 555 577, clip,grid=false]{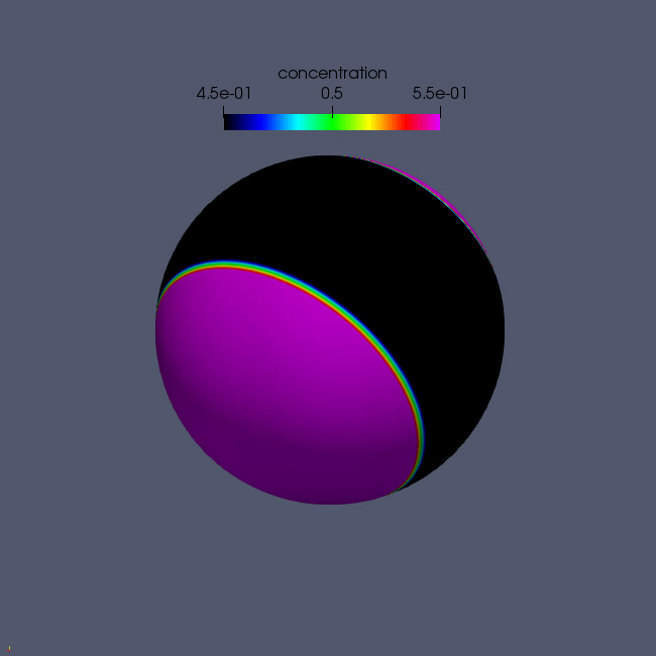}
        \put(40,102){\small{$t = 1$}}
\end{overpic}
\end{center}
\caption{\label{fig:test2c_1}
Test 2c:
Evolution of the numerical solution of the Cahn--Hilliard equation computed with mesh $\ell = 6$
for $t \in [0, 1]$. View: $x_1 x_2$-plane.
The legend is reported in every panel.
}
\end{figure}

\noindent{\bf Test 2c.}
We present the evolution of the discrete Lyapunov energy in Fig.~\ref{fig:E_CH3}~(right).
We see that the results computed with meshes $\ell = 4, 5$ are in excellent
agreement. 
In \cite{elliott2015evolving}, good convergence of the solution is obtained for short times, i.e.~till
around $t = 0.2$. For the choice of parameters associated to test 2c, which
are close to the parameter values in \cite{elliott2015evolving}, we observe
good convergence of the solution for the whole interval of time under consideration.
We notice that the discrete Lyapunov energy computed with meshes $\ell = 4, 5$ decreases (not monotonically)
till around $t = 0.35$. For $t > 0.35$, Fig.~\ref{fig:E_CH3} (right) suggests that the solution converges to a time
periodic solution.

%
%
%
%
%

We show the evolution of the solution computed with mesh $\ell = 6$
in Fig.~\ref{fig:test2c_1}. Note that we use the truncated scale here for $c$ to better illustrate
the solution, which is otherwise rather diffusive (since $\epsilon=0.1$ is used).
From Fig.~\ref{fig:test2c_1}, we see that the solution configuration remains the same for $t > 0.35$
as the surface gets dilated and shrunk. The patterns and the evolution look very similar to those
in Figure 4 of \cite{elliott2015evolving}.

\subsection{Colliding spheres}\label{sec:merging}
We consider an evolving surface $\Gamma(t)$ that undergoes a topological change and experiences a
local singularity. The dynamics of $\Gamma(t)$ is inspired by tests presented in \cite{GOReccomas,olshanskii2017btrace,lehrenfeld2018stabilized}.
The computational domain is $\Omega = [-10/3,10/3]\times[-5/3,5/3]^2$ and
the evolving surface is the zero level set of level set function
\begin{align}\label{eq:merging_level_set}
\phi(\bx, t) = 1 - \frac{1}{|| \bx - \bx_c^+(t)||^3} - \frac{1}{|| \bx - \bx_c^-(t)||^3},
\end{align}
with
\begin{align}\label{eq:centers}
\bx_c^\pm(t) = \pm\, \left(\frac{3}{2} - w\, t, 0, 0 \right),
\end{align}
for $t \in [0, 1.5/|w|]$.
The initial configuration $\Gamma(0)$ is close to two balls of radius 1, centered at
$\bx_c^\pm(0)$. Parameter $w$ is the collision speed.
For $t > 0$, $w>0$ the two spheres approach each other until
time $\tilde{t}\approx 0.235/w$, when they touch at the origin; then, for $t \in (\tilde{t},  1.5/w)$, surface $\Gamma(t)$ is simply connected.

In the vicinity of $\Gamma(t)$, the gradient $\nabla \phi$ and the time derivative $\partial_t \phi$
are well defined and given by simple algebraic expressions.
The normal velocity field of $\Gamma(t)$ can be computed to be
\begin{align}
\bw = -  \frac{\partial_t \phi}{| \nabla \phi |^2} \nabla \phi. \el
\end{align}
We assume the following material velocity:
\begin{align}
\label{eq:merging_velocity}
\bu=\left(1-\tanh{}\left(\frac{ |x_1 |}{0.1}\right)\right)\bw+ \tanh{}\left(\frac{ |x_1 |}{0.1}\right)( -\text{sgn}(x_1) w, 0, 0),
\end{align}
which models the parallel advection away from the merging line and nearly normal motion near the $x_1=0$ plane.

We consider different meshes for $\Omega$ with mesh size $h=2^{-\ell-2}10/3$, where $\ell$ is the refinement level.

The setup is a simple model of two plasma membrane fusion. In all further experiments we take $\epsilon = 0.01$.
We remark that this is a realistic value of interface thickness for applications related to phase separation in lipid bilayers. In fact, if we consider a typical giant vesicle with an average diameter
of 30\,$\mu$m, on which phase separation can be visualized using fluorescence microscopy
\cite{VEATCH20033074} with a resolution of about 300 nm, the thickness of transition region between the
phases is approximately 1\% of vesicle diameter.

\begin{figure}[ht]
\begin{center}
\href{https://youtu.be/I0mZcmFdyEo}{
\begin{overpic}[width=.48\textwidth,
	 clip,grid=false]{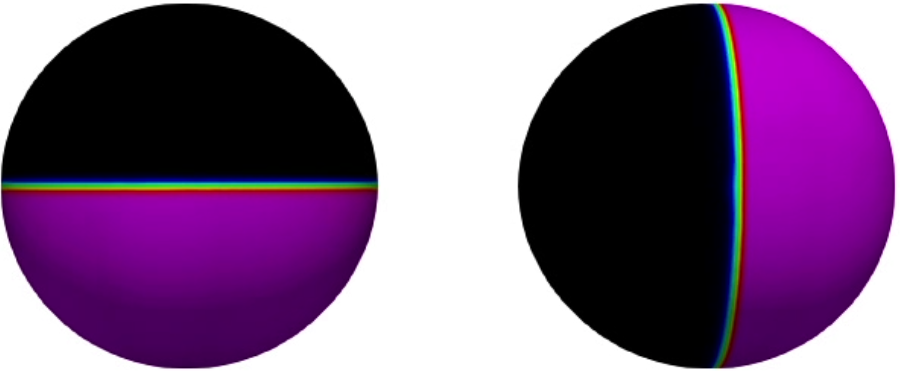}
         \put(40, 42){\small{$t = 0.05$}}
\end{overpic}
\begin{overpic}[width=.48\textwidth,  clip,grid=false]{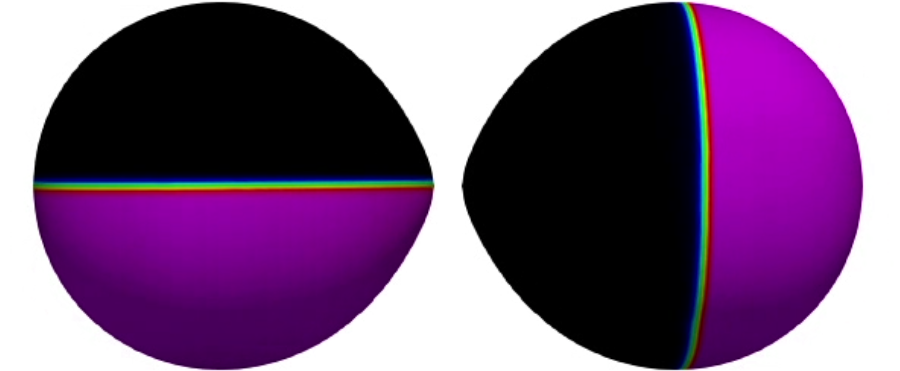}
         \put(40, 42){\small{$t = 0.23$}}
\end{overpic}
}
\vskip .4cm
\href{https://youtu.be/I0mZcmFdyEo}{
\begin{overpic}[width=.48\textwidth,  clip,grid=false]{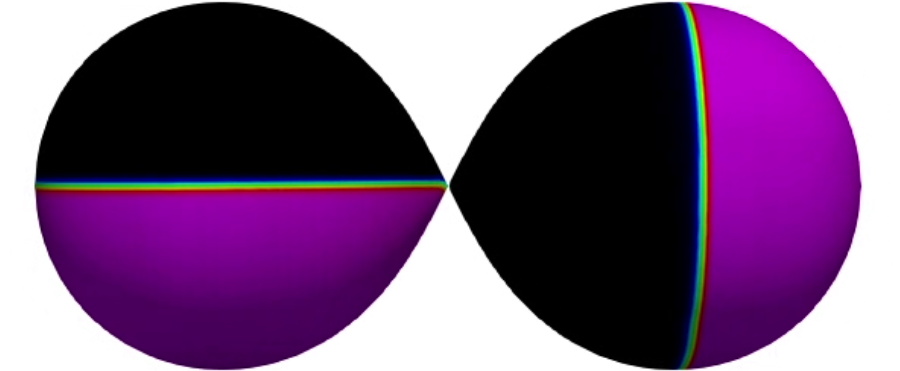}
         \put(40, 42){\small{$t = 0.24$}}
\end{overpic}
\begin{overpic}[width=.48\textwidth,  clip,grid=false]{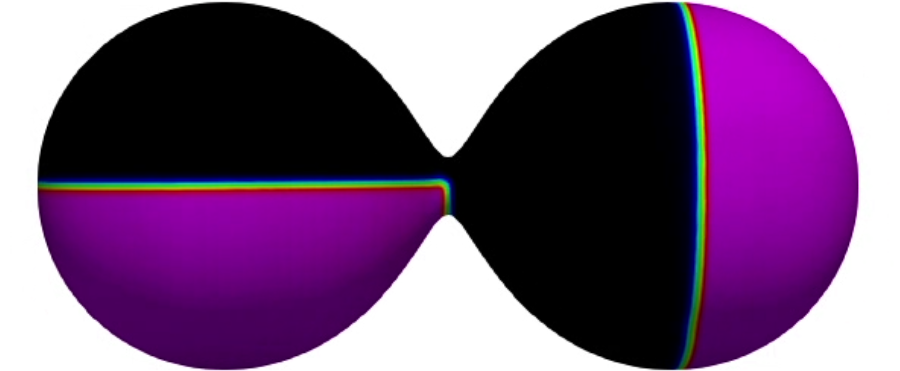}
         \put(40, 43){\small{$t = 0.25$}}
\end{overpic}
}
\vskip .4cm
\href{https://youtu.be/I0mZcmFdyEo}{
\begin{overpic}[width=.48\textwidth,  clip,grid=false]{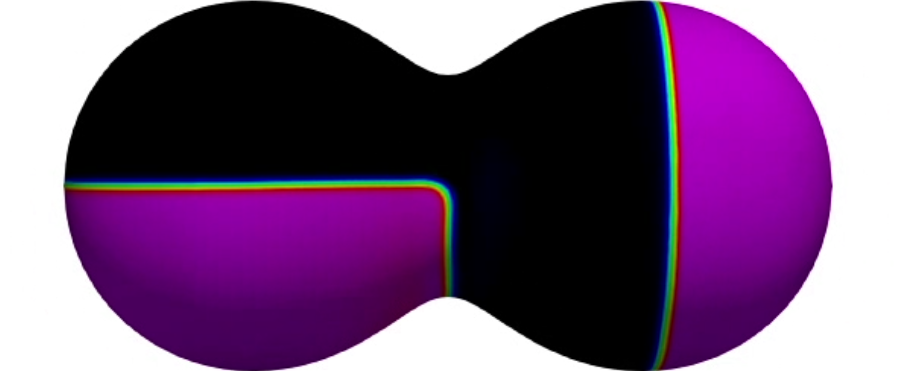}
         \put(42, 50){\small{$t = 0.4$}}
\end{overpic}
\begin{overpic}[width=.48\textwidth, clip,grid=false]{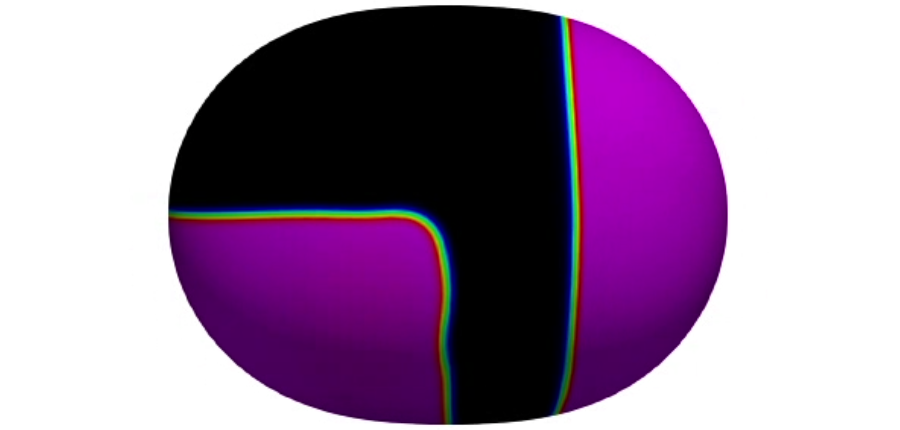}
         \put(44, 50){\small{$t = 1$}}
\end{overpic}
}
\end{center}
\caption{\label{fig:coll_spheres_steady}l
Colliding spheres, pre-separated phases: Evolution of the numerical solution of the Cahn--Hilliard
equation computed with mesh $\ell =6$ for $t \in [0,1]$. View: $x_1 x_2$-plane.
The legend is the same as in Fig.~\ref{fig:test2b}. Click any picture above to run the full animation}
\end{figure}

\subsubsection{Colliding spheres with pre-separated phases}\label{sec:coll_steady}
For this test, we set $w=1$. The configuration at $t = 3/2$ is the ball centered around 0
with radius $2^{1/3}$. We end the simulation at $t = 1$, when the two balls have merged but
they have not evolved into one sphere yet.

As initial solution, we take $c_{x_2}$ with $\delta = 0$ for the ball initially centered at $\bx_c^-(0)$
and $c_{x_1}$  with $\delta = 0$
for the ball initially centered at $ \bx_c^+(0)$. See \eqref{steady_CH}
for the definition of $c_{x_i}$.
We set $\Delta t = 0.001$.

Fig.~\ref{fig:coll_spheres_steady} displays the evolution of the surface and the solution
for $t \in [0.05,1]$. Before the
two balls touch at $t = \tilde{t}$, we see that the ball centered at $\bx_c^-(t)$ (resp., $\bx_c^+(t)$) has a horizontal
(resp., vertical) interface separating a pink and black domain of equal size. For $t >  \tilde{t}$,
the black and pink domain of the ball centered at $\bx_c^-(t)$ get in contact with the black domain
of the ball centered at $\bx_c^+(t)$. Between the pink domain of the ball centered at $\bx_c^-(t)$ and the black domain
of the ball centered at $\bx_c^+(t)$ an interface gets formed. Such interface is virtually located on the curve of minimal length
on the simply connected surface $\Gamma(t)$, for $t >  \tilde{t}$.
This is consistent with the well-known limiting (as $\epsilon\to0$) behavior of
the stationary  Cahn-Hilliard problem that corresponds to solving the isoperimetric problem~\cite{modica1987gradient,sternberg1988effect}.

\subsubsection{Pattern formation on colliding sphere}\label{sec:coll_dyn}


\begin{figure}[ht]
\begin{center}
 \href{https://youtu.be/8JwOQBckUKk}{
\begin{overpic}[width=.48\textwidth, viewport=25 195 656 450, clip,grid=false]{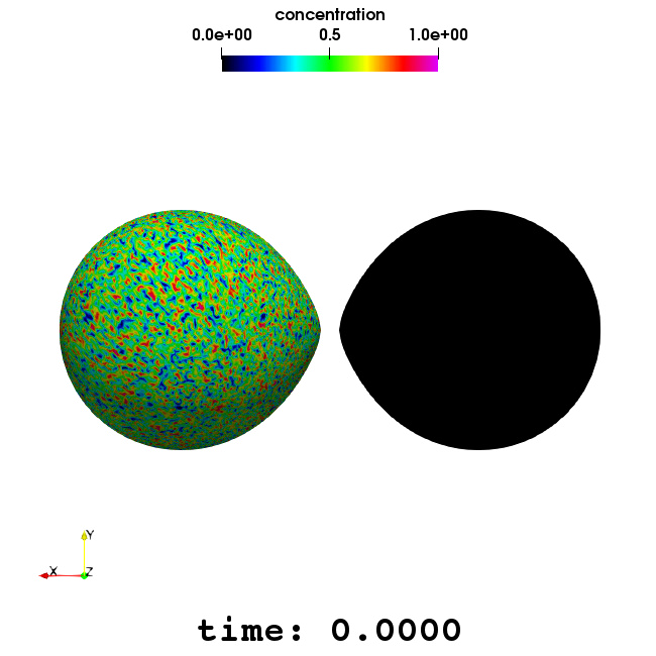}
         \put(44, 42){\small{$t = t_s$}}
\end{overpic}
\begin{overpic}[width=.48\textwidth, viewport=25 195 656 450, clip,grid=false]{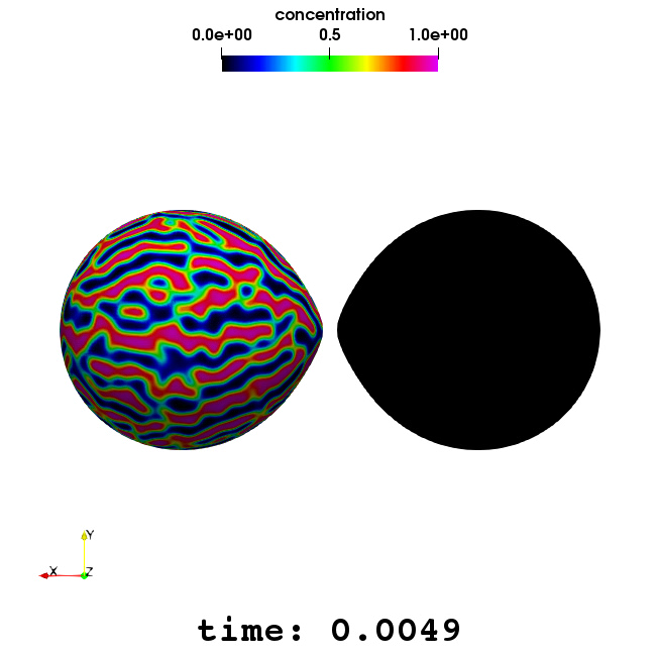}
         \put(35, 42){\small{$t =t_s + 0.005$}}
\end{overpic}}
\vskip .4cm
\href{https://youtu.be/8JwOQBckUKk}{
\begin{overpic}[width=.48\textwidth, viewport=25 195 656 450, clip,grid=false]{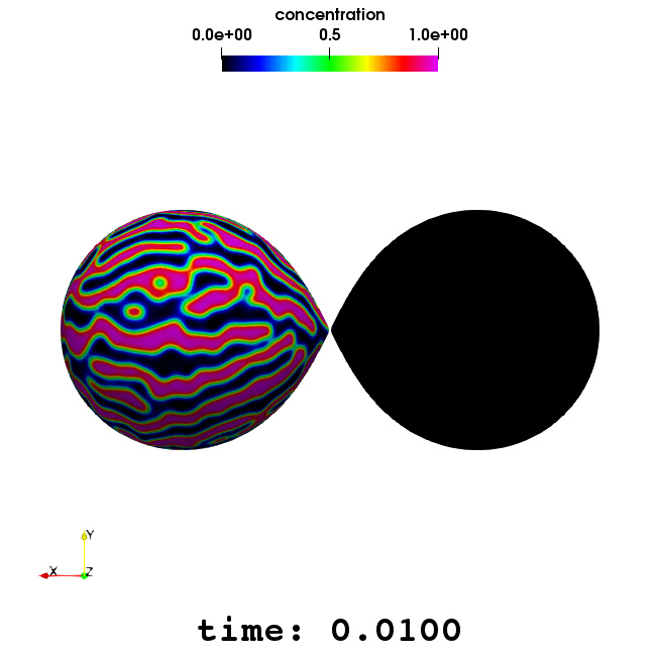}
         \put(37, 42){\small{$t = t_s + 0.01$}}
\end{overpic}
\begin{overpic}[width=.48\textwidth, viewport=25 195 656 450, clip,grid=false]{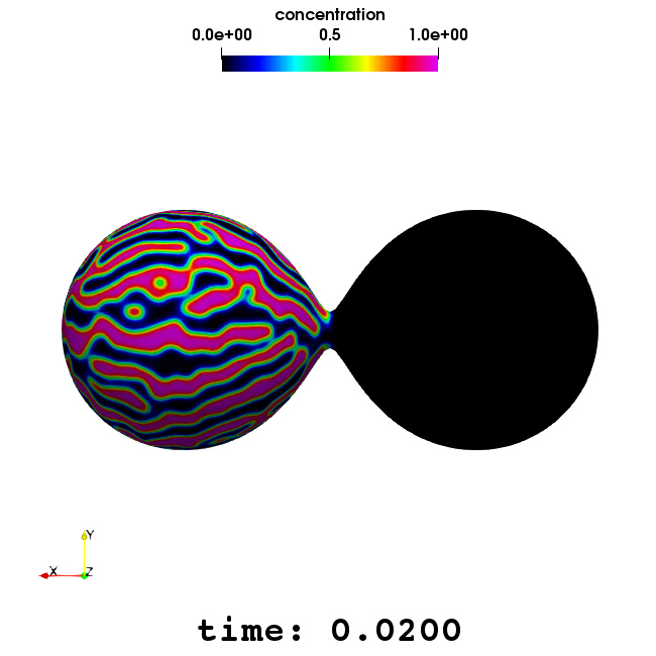}
         \put(37, 42){\small{$t = t_s+ 0.02$}}
\end{overpic}}
\vskip .4cm
\href{https://youtu.be/8JwOQBckUKk}{
\begin{overpic}[width=.48\textwidth, viewport=25 175 656 470, clip,grid=false]{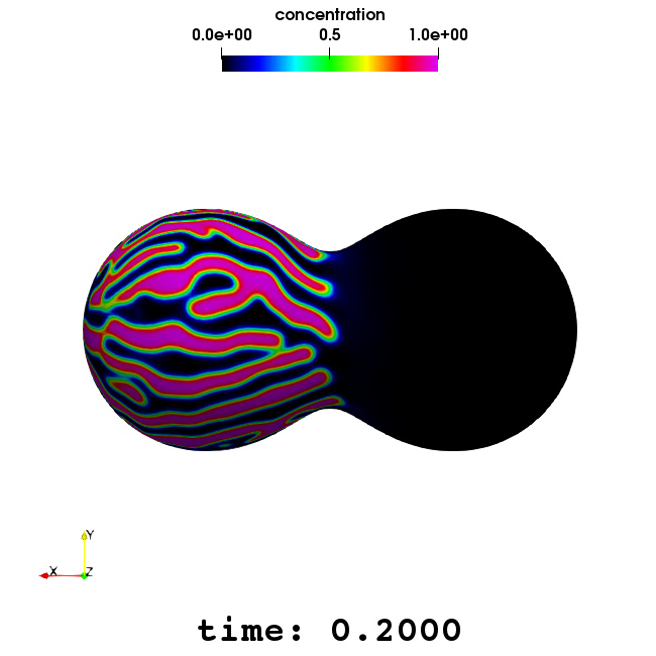}
         \put(38, 46){\small{$t = t_s+ 0.2$}}
\end{overpic}
\begin{overpic}[width=.48\textwidth, viewport=25 175 656 470, clip,grid=false]{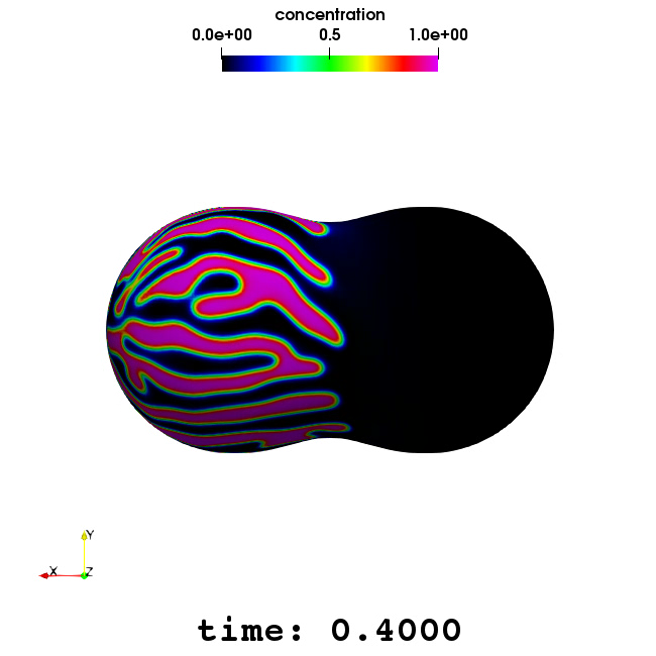}
         \put(38, 46){\small{$t = t_s+0.4$}}
\end{overpic}}
\vskip .4cm
\href{https://youtu.be/8JwOQBckUKk}{
\begin{overpic}[width=.48\textwidth, viewport=25 165 656 500, clip,grid=false]{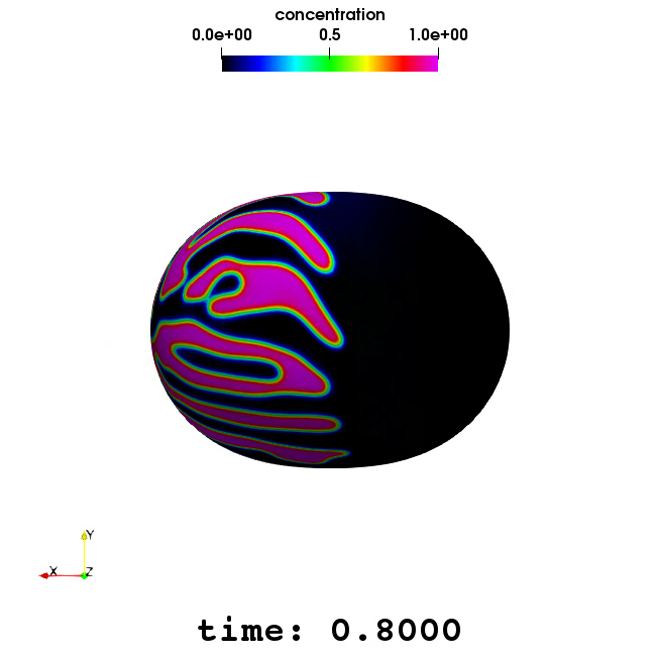}
         \put(38, 52){\small{$t = t_s+0.8$}}
\end{overpic}
\begin{overpic}[width=.48\textwidth, viewport=25 165 656 500, clip,grid=false]{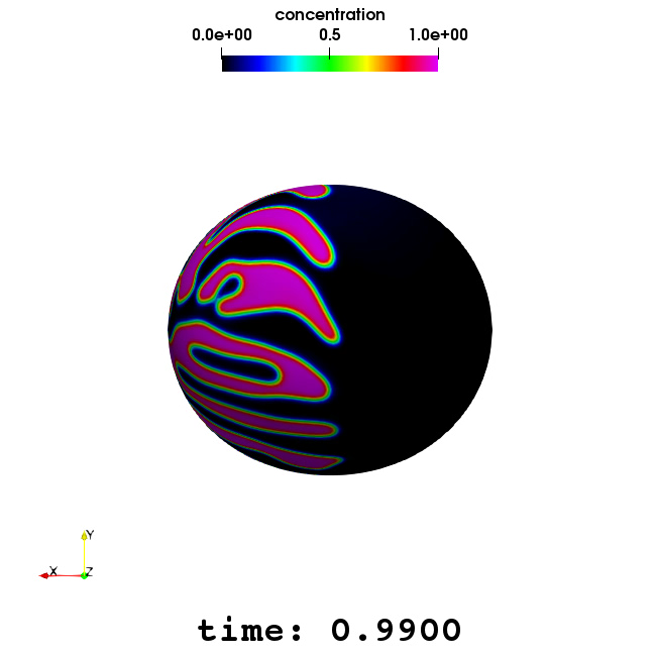}
         \put(38, 52){\small{$t = t_s+1$}}
\end{overpic}}
\end{center}
\caption{\label{fig:coll_spheres_1}
Slow collision: Evolution of the numerical solution of the Cahn--Hilliard
equation computed with mesh $\ell =6$ for time $t \in [t_s, t_s+ 1]$, $t_s=0.23$.
View: $x_1 x_2$-plane.
The legend is the same as in Fig.~\ref{fig:test2b}.
Click any picture above to run the full animation.}
\end{figure}

Let $t_s=0.23/w < \tilde{t}$ be a time shift in order to have phase separation occur close to
the surface collision. All the simulations whose results are reported in this section
start at time $t_s$. We consider two tests:
\begin{itemize}
\item[-] \emph{Slow collision}: $w = 1$, $T = 1+t_s$, $\Delta t = 10^{-4}$ for $t \in (t_s, t_s + 0.01]$ and $\Delta t = 10^{-2}$ for $t \in (t_s+0.01, t_s + 1]$.
\item[-] \emph{Fast collision}: $w = 10$, $T = 0.1+t_s$, $\Delta t = 10^{-4}$ for $t \in (t_s, t_s+ 0.01]$ and $\Delta t = 10^{-3}$ for $t \in (t_s+0.01, t_s+0.1]$.
\end{itemize}
The two initial spheres touch each other at the origin at time $\tilde{t}\approx{}0.235$ for the slow collision test and $\tilde{t} \approx{}0.0235$ for the fast collision test. 
Let $\text{rand}()$ be a uniformly distributed random number between 0 and 1.
For both tests, we take $c_0(\bx) =\text{rand}()$ for the ball initially centered at $\bx_c^-(0)$
and $c_0 = 0$ for the ball initially centered at $ \bx_c^+(0)$.

As we have observed in our previous work \cite{Yushutin_IJNMBE2019}
and is well known (see, e.g.,\cite{guillen2014second,shen_et_al2016}),
the evolution of the solution to the Cahn--Hilliard problem goes through
an initial fast phase, during which a pattern gets formed, followed by a slowdown in the process
of dissipation of the interfacial energy. In order to capture the different time scales,
we prescribe different time steps for the different stages.
It would be less intrusive to use some time-adaptivity strategy, which is not addressed in
this paper.

\begin{figure}[ht]
\begin{center}
\raisebox{.11\textwidth}{\makebox[10ex]{\small{$t = 0$}}}~~
\begin{overpic}[height=.22\textwidth, viewport=105 100 656 550, clip,grid=false]{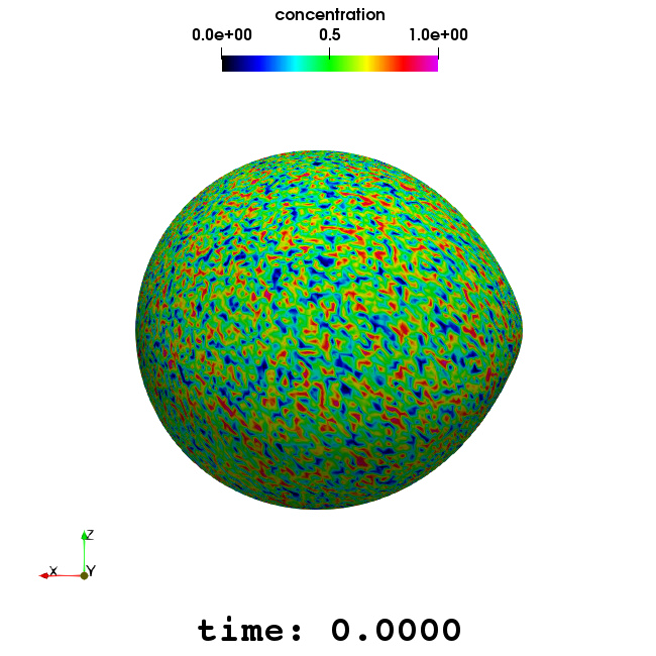}
\end{overpic}\hskip .3cm
\begin{overpic}[width=.48\textwidth, viewport=25 195 656 450, clip,grid=false]{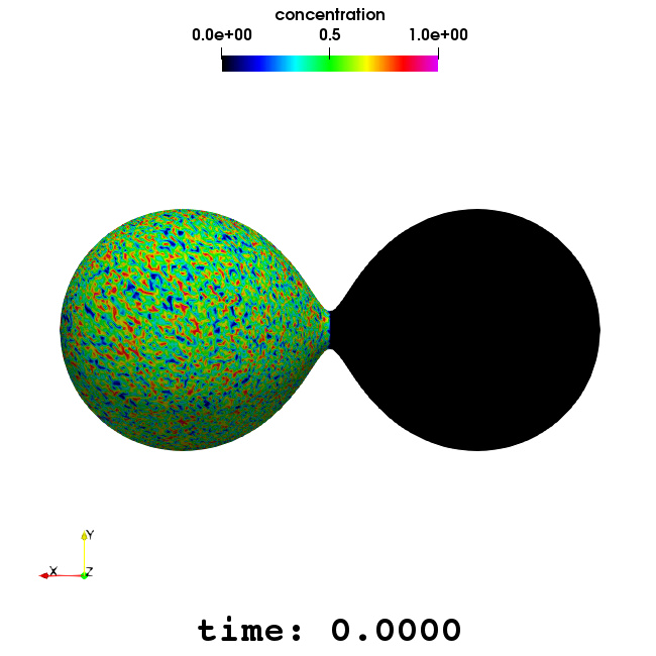}
\end{overpic}
\vskip .1cm
\raisebox{.11\textwidth}{\makebox[10ex]{\small{$t = 0.004$}}}~~
\begin{overpic}[height=.22\textwidth, viewport=105 100 656 550, clip,grid=false]{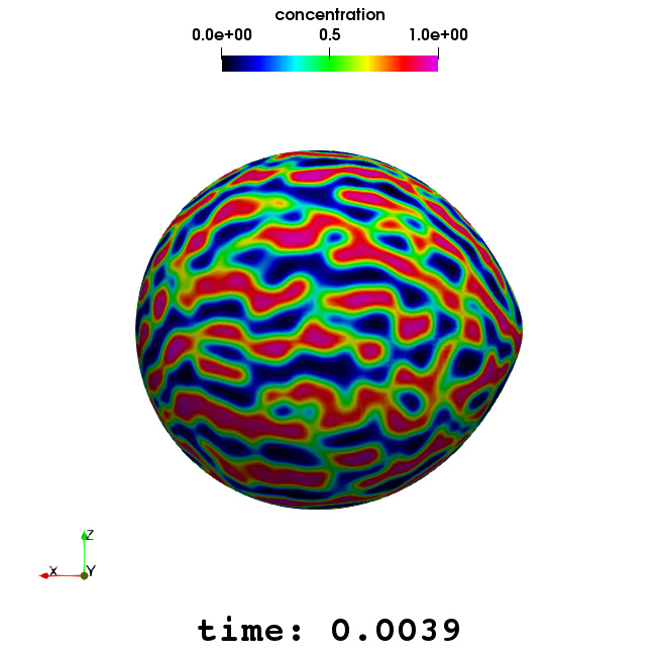}
\end{overpic}\hskip .3cm
\begin{overpic}[width=.48\textwidth, viewport=25 195 656 450, clip,grid=false]{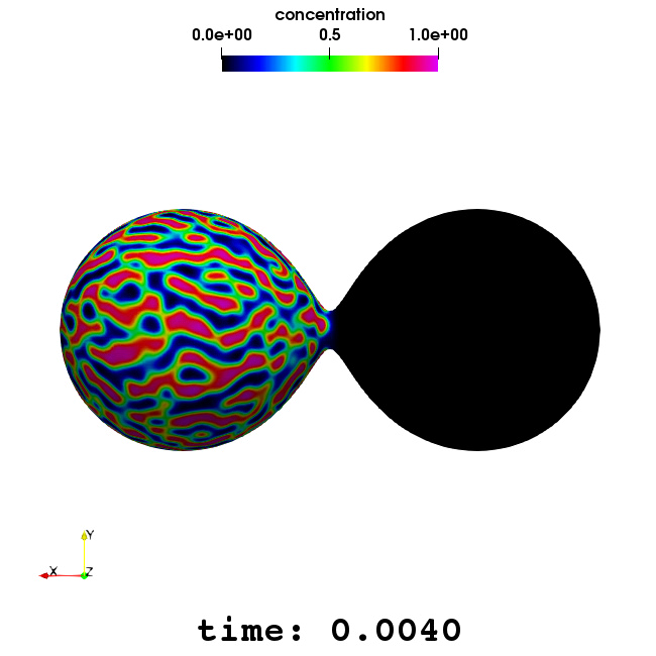}
\end{overpic}
\vskip .1cm
\raisebox{.11\textwidth}{\makebox[10ex]{\small{$t = 0.02$}}}~~
\begin{overpic}[height=.22\textwidth, viewport=105 100 656 520, clip,grid=false]{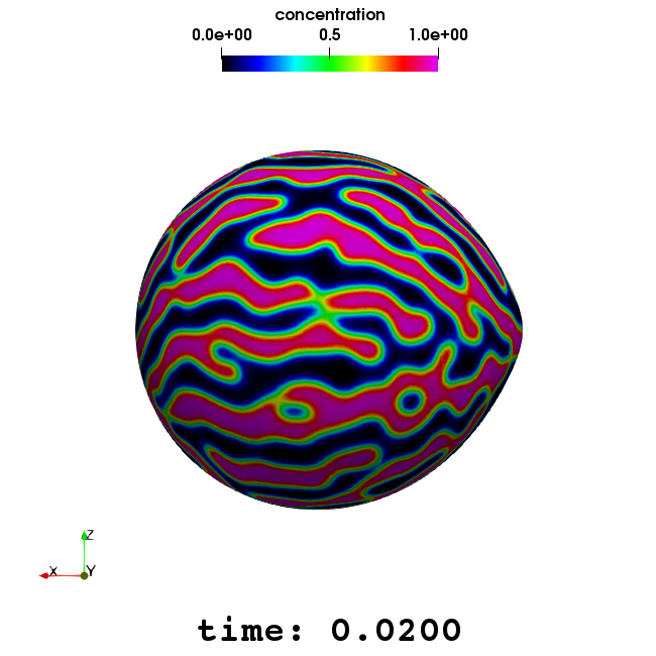}
\end{overpic}\hskip .3cm
\begin{overpic}[width=.48\textwidth, viewport=25 165 656 470, clip,grid=false]{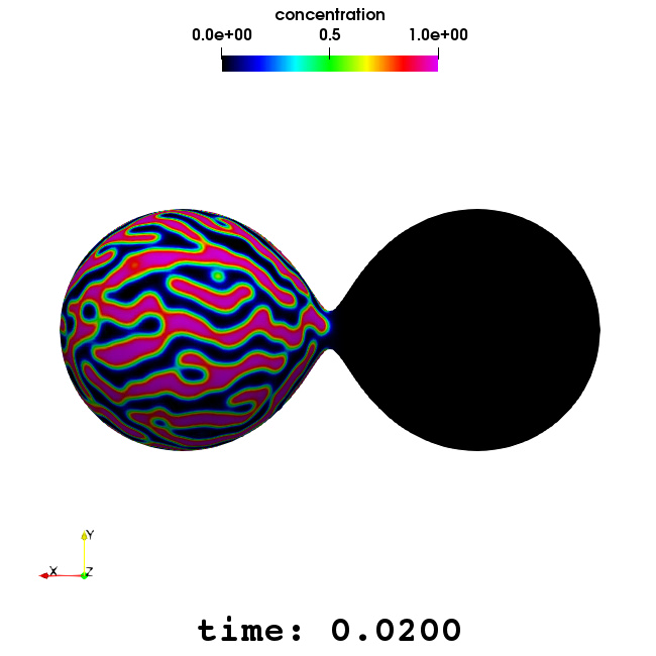}
\end{overpic}
\vskip .1cm
\raisebox{.11\textwidth}{\makebox[10ex]{\small{$t = 0.4$}}}~~
\begin{overpic}[height=.22\textwidth, viewport=105 100 656 520, clip,grid=false]{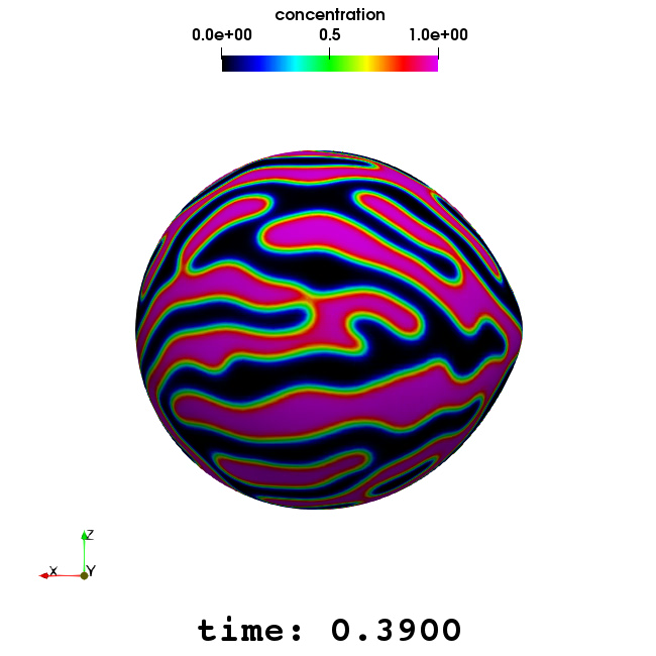}
\end{overpic}\hskip .3cm
\begin{overpic}[width=.48\textwidth, viewport=25 165 656 470, clip,grid=false]{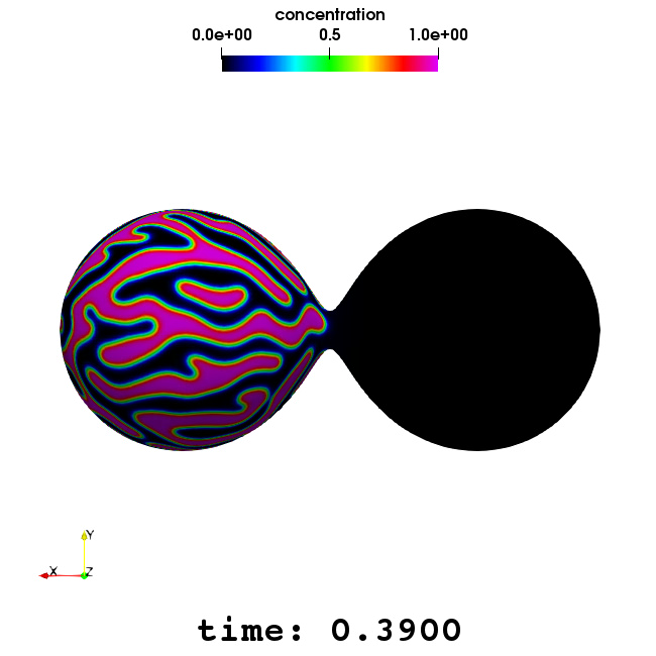}
\end{overpic}
\vskip -.1cm
\raisebox{.11\textwidth}{\makebox[10ex]{\small{$t = 1$}}}~~
\begin{overpic}[height=.22\textwidth, viewport=105 100 656 520, clip,grid=false]{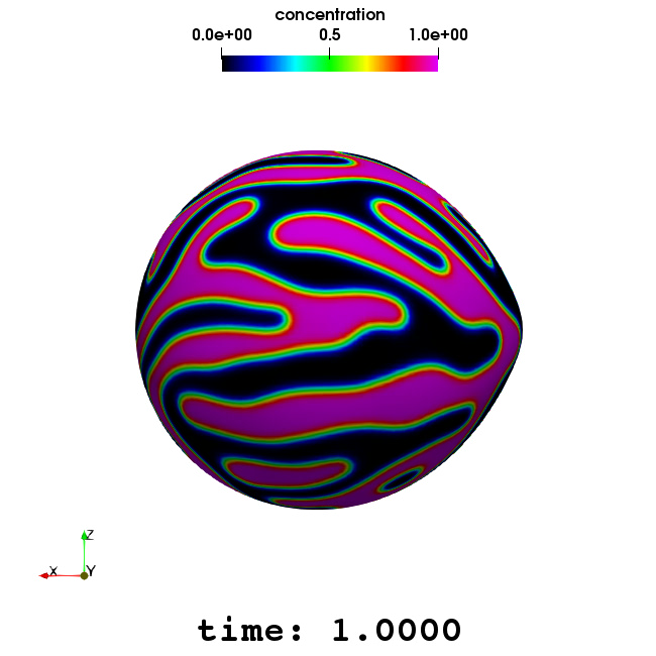}
\end{overpic}\hskip .3cm
\begin{overpic}[width=.48\textwidth, viewport=25 165 656 470, clip,grid=false]{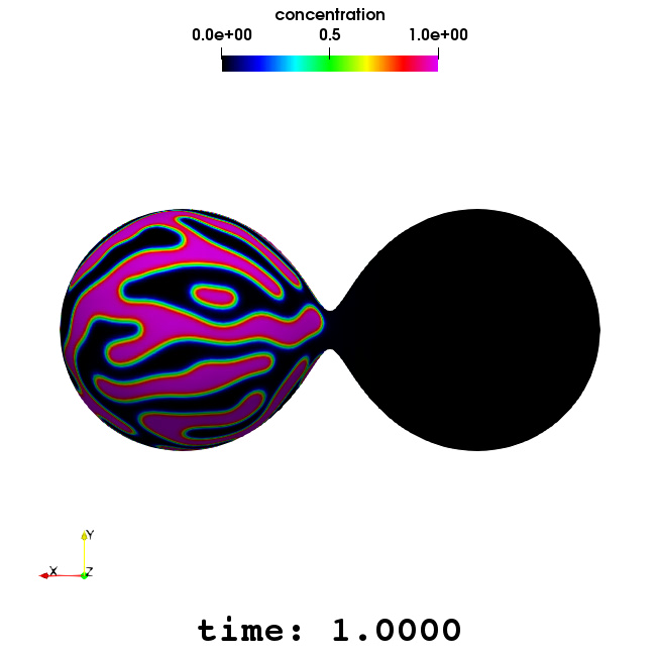}
\end{overpic}
\end{center}
\caption{\label{fig:coll_spheres_still_half}
Comparison between phase separation on two frozen domains for $t \in [0,1]$: single ball (left)
and the hourglass (right). View: $x_1 x_2$-plane.
The legend is the same as in Fig.~\ref{fig:test2b}.}
\end{figure}

Fig.~\ref{fig:coll_spheres_1} shows the  evolution of the surface and the solution
for the slow collision test for $t \in [t_s,t_s + 1]$.
We see the quick formation of a pattern, that is visible already at $t = t_s + 0.005$. Thus,
when the two balls touch at $\tilde{t}  \approx{}0.235$, the already phase separated solution of the ball centered at $\bx_c^-(t)$
enters in contact with the black domain that covers the entire ball centered at $\bx_c^+(t)$.
Similarly to what observed in Sec.~\ref{sec:coll_steady}, for $t >  \tilde{t}$ we see that the interfaces
between the pink domains and the large black domain tend to align with the curve of minimal length
on surface $\Gamma(t)$. Hence, the pink domains remain confined on the portion of the surface in the $x_1 < 0$
semi-space. Possible factors that come into play to determine the appearance of the
pattern on $\Gamma(t)$ are:
\begin{enumerate}
\item The droplet shape of the balls before merging and the bottleneck shape, with a short perimeter cross-section, after merging.
\item The contact with the black domain covering the ball centered at $\bx_c^+(t)$.
\item The colliding dynamics.
\end{enumerate}
To understand the role of factors 1.~and 2., we will study phase separation on
two surfaces frozen in time (not evolving) relative to the slow collision test (i.e., $w = 1$
and $t_s = 0.23$):
half of domain $\Gamma(t_s)$, hereafter called single ball, and surface $\Gamma(t_s + 0.02)$,
hereafter called hourglass. The fast collision test, i.e. $w = 10$, is meant to
understand the role of factor 3.

We start by comparing phase separation on the single ball and on the hourglass.
The initial condition for the single ball is simply $c_0(\bx) =\text{rand}()$ for the entire domain.
The initial condition for simulation on the hourglass mimics
the initial condition for the slow collision test, i.e.~$c_0(\bx) =\text{rand}()$ for $\bx$
with $x_1 \leq 0$ and $c_0 = 0$ for $\bx$ with $x_1>0$.
We consider the same meshes described in Sec.~\ref{sec:merging}.

In Fig.~\ref{fig:coll_spheres_still_half}, we see the results for both simulations for $t \in [0,1]$.
We observe elongated pink domains aligned with the $x_1$-axis, i.e.~the horizontal
axis in the figure, for both cases. Since both surfaces are not evolving,
we believe this patter is solely due to the shapes (droplet and bottleneck). On the other hand,
the interaction with the large black domain on the hourglass
does not seem to significantly influence the pattern for $t < 1$.

\begin{figure}[ht]
\begin{center}
\href{https://youtu.be/BW9ryD1Shug}{
\begin{overpic}[width=.48\textwidth, viewport=25 195 656 450, clip,grid=false]{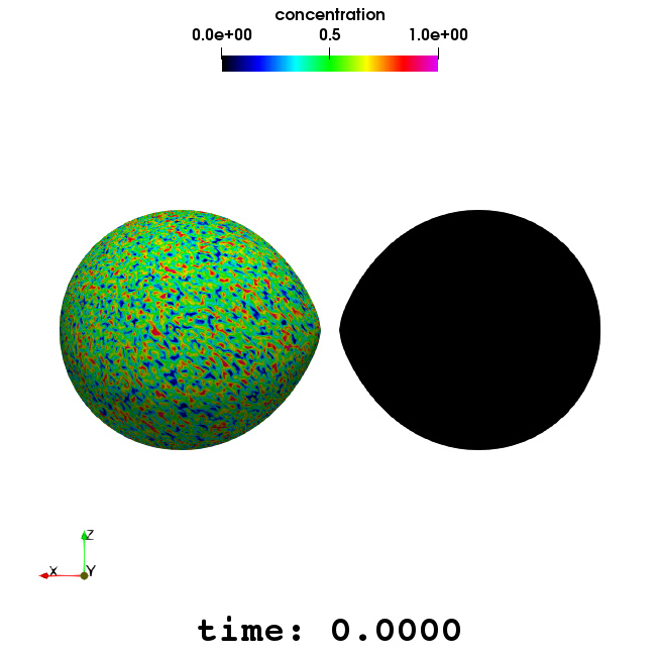}
         \put(43, 42){\small{$t = t_s$}}
\end{overpic}
\begin{overpic}[width=.48\textwidth, viewport=25 175 656 470, clip,grid=false]{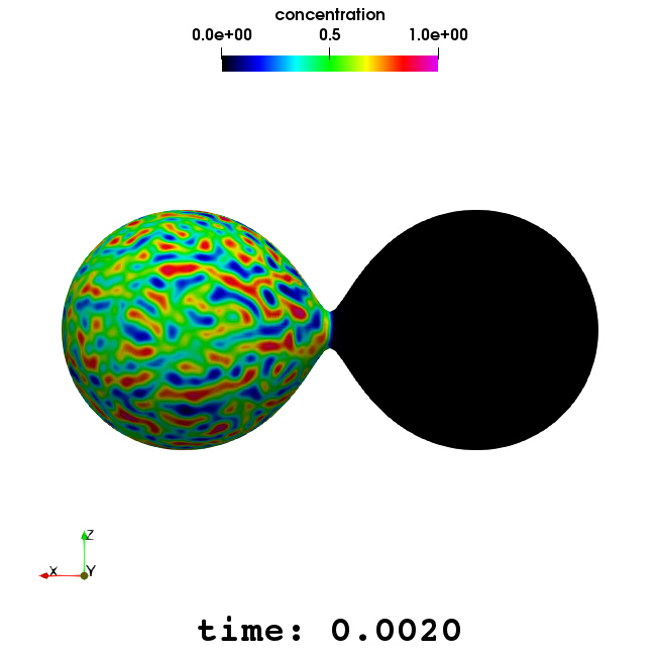}
         \put(36, 44){\small{$t = t_s+0.002$}}
\end{overpic}}
\vskip .4cm
\href{https://youtu.be/BW9ryD1Shug}{
\begin{overpic}[width=.48\textwidth, viewport=25 195 656 450, clip,grid=false]{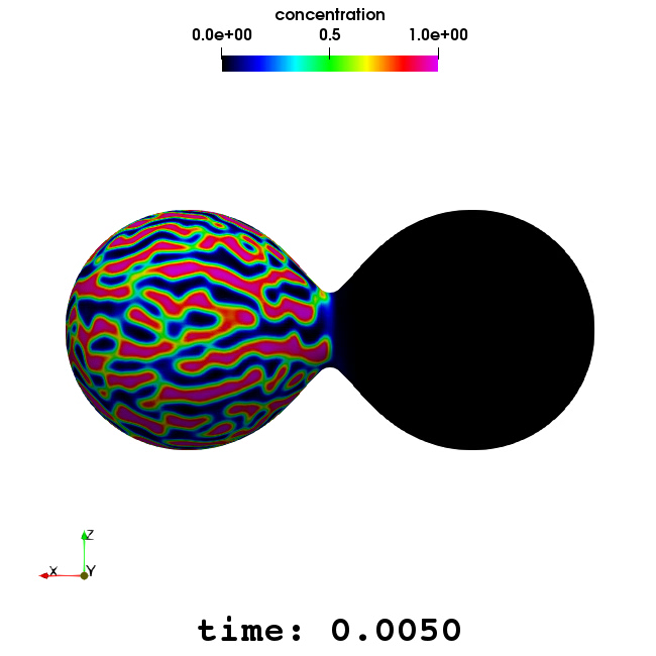}
         \put(36, 44){\small{$t = t_s+0.005$}}
\end{overpic}
\begin{overpic}[width=.48\textwidth, viewport=25 175 656 470, clip,grid=false]{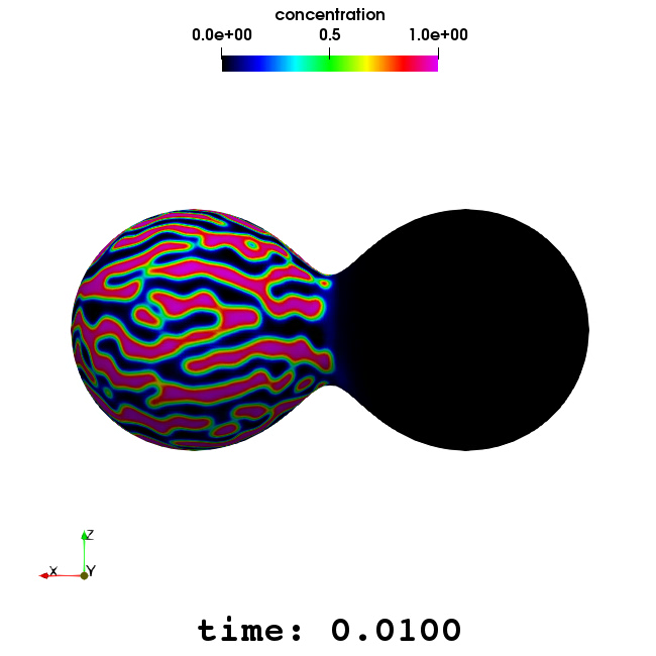}
         \put(38, 44){\small{$t = t_s+0.01$}}
\end{overpic}}
\vskip .4cm
\href{https://youtu.be/BW9ryD1Shug}{
\begin{overpic}[width=.48\textwidth, viewport=25 175 656 470, clip,grid=false]{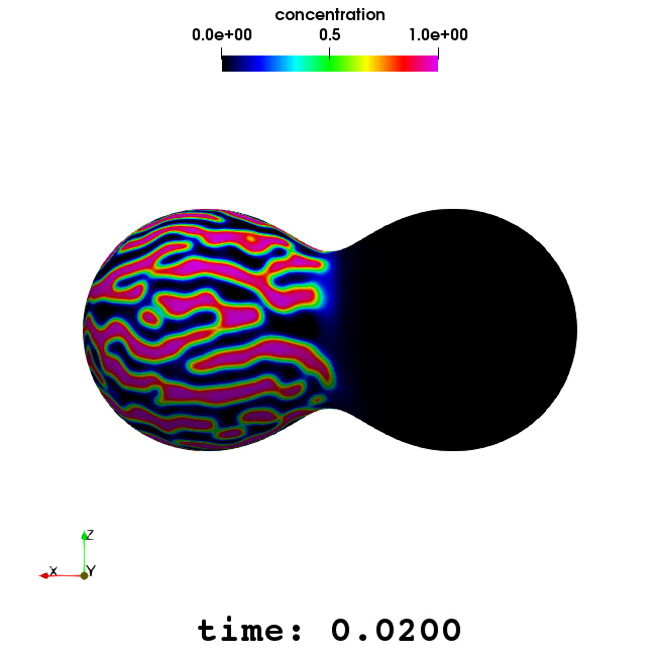}
         \put(38, 46){\small{$t = t_s+0.02$}}
\end{overpic}
\begin{overpic}[width=.48\textwidth, viewport=25 175 656 470, clip,grid=false]{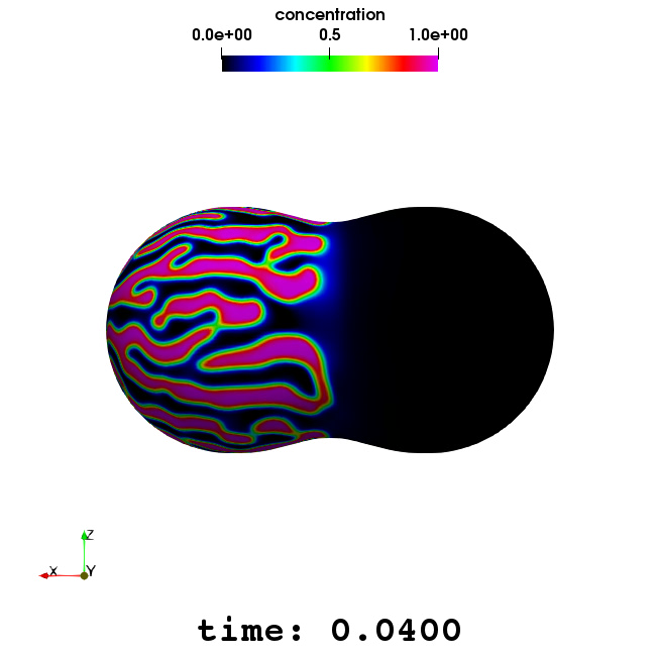}
         \put(38, 46){\small{$t = t_s+0.04$}}
\end{overpic}}
\vskip .4cm
\href{https://youtu.be/BW9ryD1Shug}{
\begin{overpic}[width=.48\textwidth, viewport=25 165 656 500, clip,grid=false]{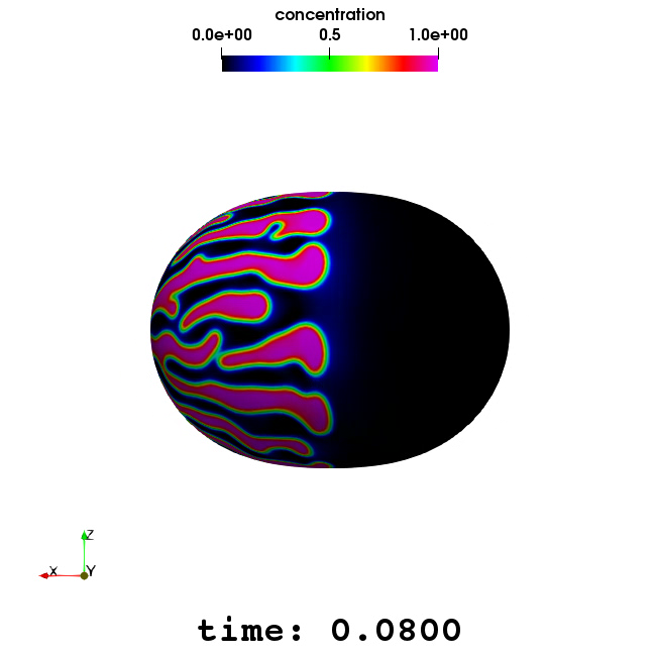}
         \put(38, 52){\small{$t = t_s+0.08$}}
\end{overpic}
\begin{overpic}[width=.48\textwidth, viewport=25 165 656 500, clip,grid=false]{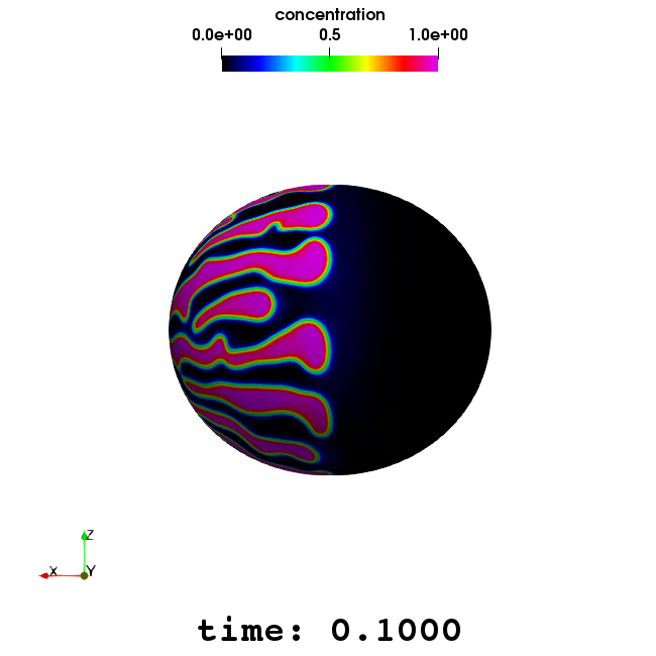}
         \put(38, 52){\small{$t = t_s+0.1$}}
\end{overpic}}
\end{center}
\caption{\label{fig:coll_spheres_10}
Fast collision: Evolution of the numerical solution of the Cahn--Hilliard
equation computed with mesh $\ell =6$ for $t \in [t_s+0, t_s+0.1]$, $t_s=0.023$. View: $x_1 x_2$-plane.
The legend is the same as in Fig.~\ref{fig:test2b}. Click any picture above to run the full animation.}
\end{figure}

Next, we present the results for the fast collision test. We report in
Fig.~\ref{fig:coll_spheres_10} the evolution of the surface and the solution for $t \in [t_s,t_s+0.1]$.
We see that, due to the fast collision dynamics, the two balls touch before
phase separation occurs on the balls centered at $\bx_c^-(t)$.
Upon phase separation, we see the same elongated pink domains
aligned with the $x_1$-axis already observed in Fig.~\ref{fig:coll_spheres_1}
and \ref{fig:coll_spheres_still_half}.
The notable difference between the corresponding solutions
in Fig.~\ref{fig:coll_spheres_1} and \ref{fig:coll_spheres_10}
(taking into account the 10 times faster dynamics in the latter)
are the thicker interfaces along the bottleneck in Fig.~\ref{fig:coll_spheres_10},
where we see larger blue regions in the $t = t_s+0.02, t_s+ 0.04, t_s+0.08$ panels,
and the longer portions of interface on the bottleneck in Fig.~\ref{fig:coll_spheres_10},
$t = t_s+0.1$ panel, which shows a `mushroom' shape of the pink fingers
comparing to a more stretched, oval shape for slow collision motion.

\begin{figure}[ht]
\begin{center}
\begin{overpic}[height=.22\textwidth, viewport=90 128 570 510, clip,grid=false]{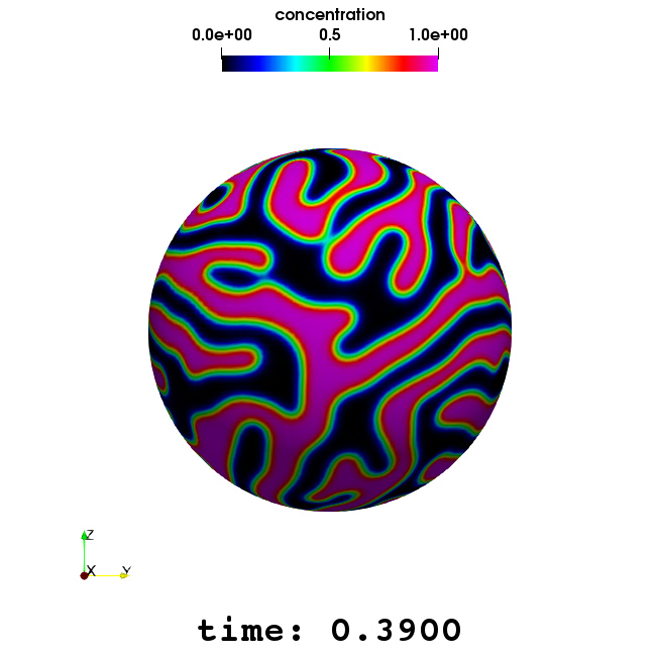}
         \put(42, 85){\small{$t = 0.4$}}
         \put(-57, 40){\small{Single ball}}
         \put(-57, 30){\small{(frozen surface)}}
\end{overpic}\hskip 1cm
\begin{overpic}[height=.23\textwidth, viewport=90 128 570 510, clip,grid=false]{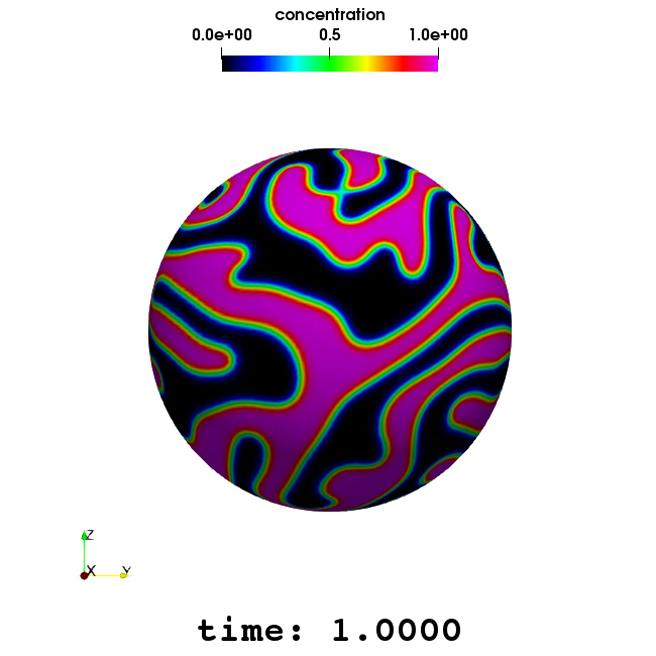}
         \put(42, 85){\small{$t = 1$}}
\end{overpic}
\vskip .4cm
\begin{overpic}[width=.48\textwidth, viewport=25 165 656 500, clip,grid=false]{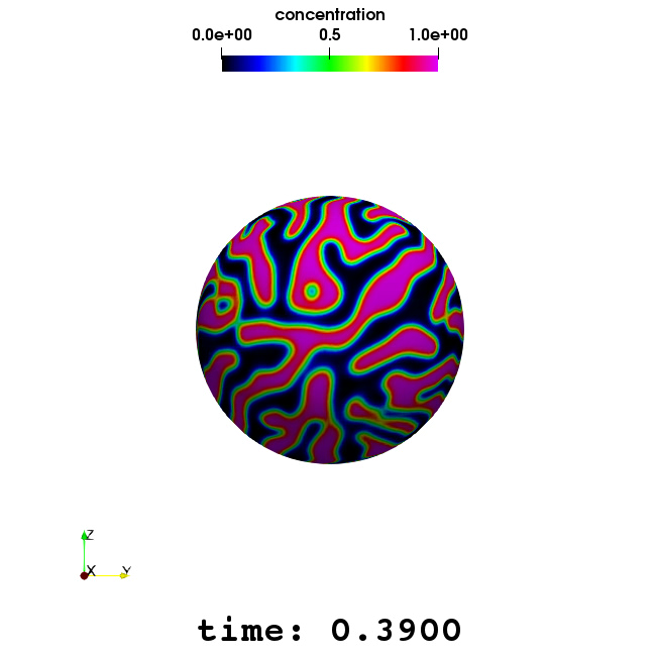}
         \put(42, 52){\small{$t = 0.4$}}
         \put(-12, 25){\small{Slow collision}}
\end{overpic}\hskip -1.5cm
\begin{overpic}[width=.48\textwidth, viewport=25 165 656 500, clip,grid=false]{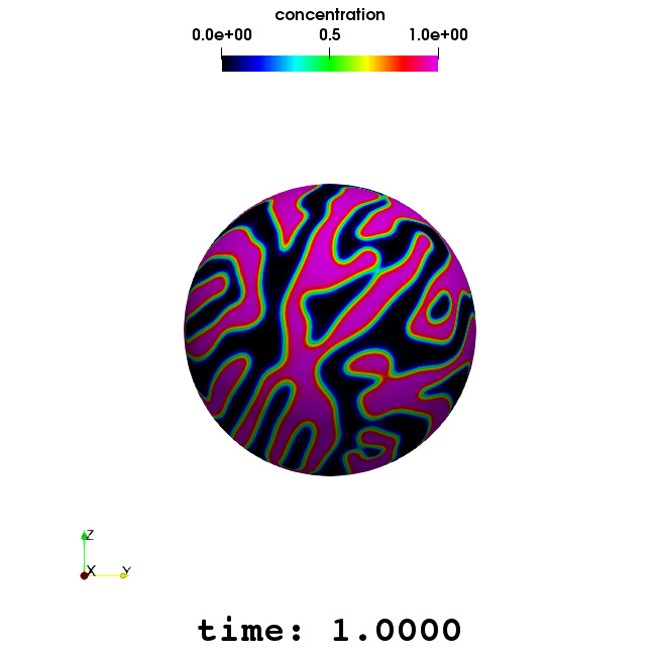}
         \put(42, 52){\small{$t = 1$}}
\end{overpic}
\vskip .4cm
\begin{overpic}[width=.48\textwidth, viewport=25 165 656 500, clip,grid=false]{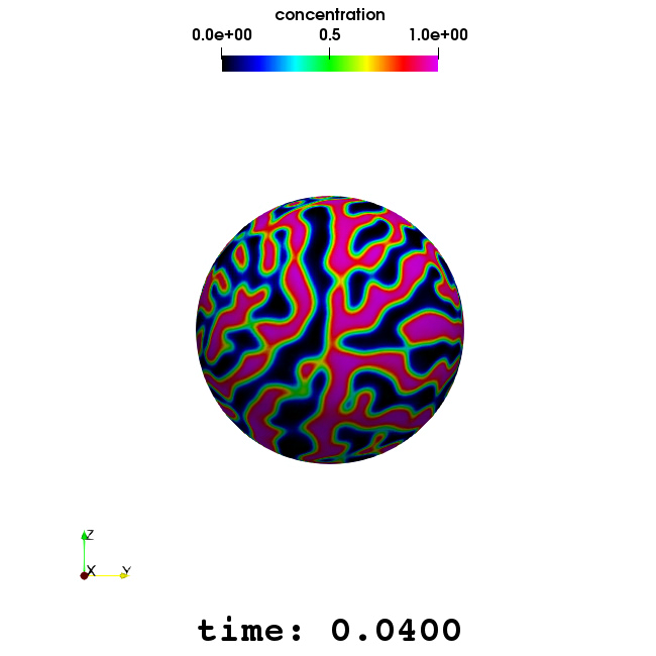}
         \put(42, 52){\small{$t = 0.04$}}
         \put(-12, 25){\small{Fast collision}}
\end{overpic}\hskip -1.5cm
\begin{overpic}[width=.48\textwidth, viewport=25 165 656 500, clip,grid=false]{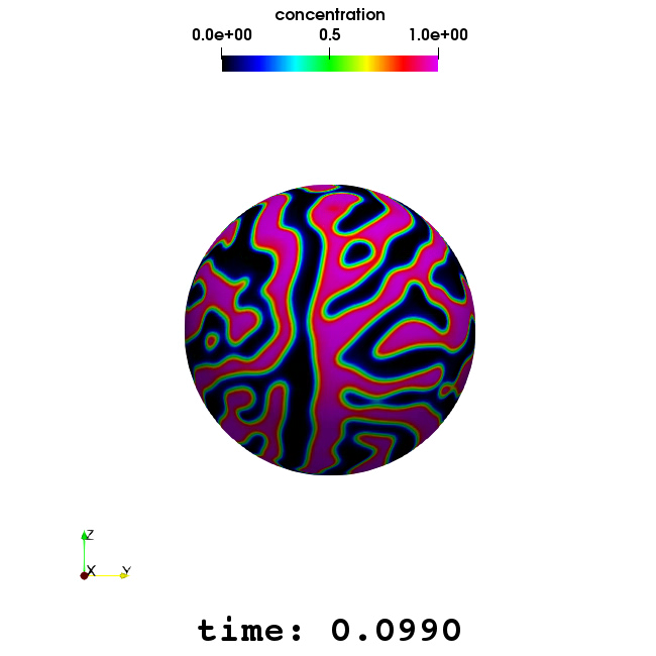}
         \put(42, 52){\small{$t = 0.1$}}
\end{overpic}
\end{center}
\caption{\label{fig:coll_spheres_side}
Colliding spheres, pattern formation: numerical solution of the Cahn--Hilliard at two different times for the
single ball (top), slow collision (center), and fast collision (bottom) tests.
View: $x_2 x_3$-plane.
The legend is the same as in Fig.~\ref{fig:test2b}.}
\end{figure}
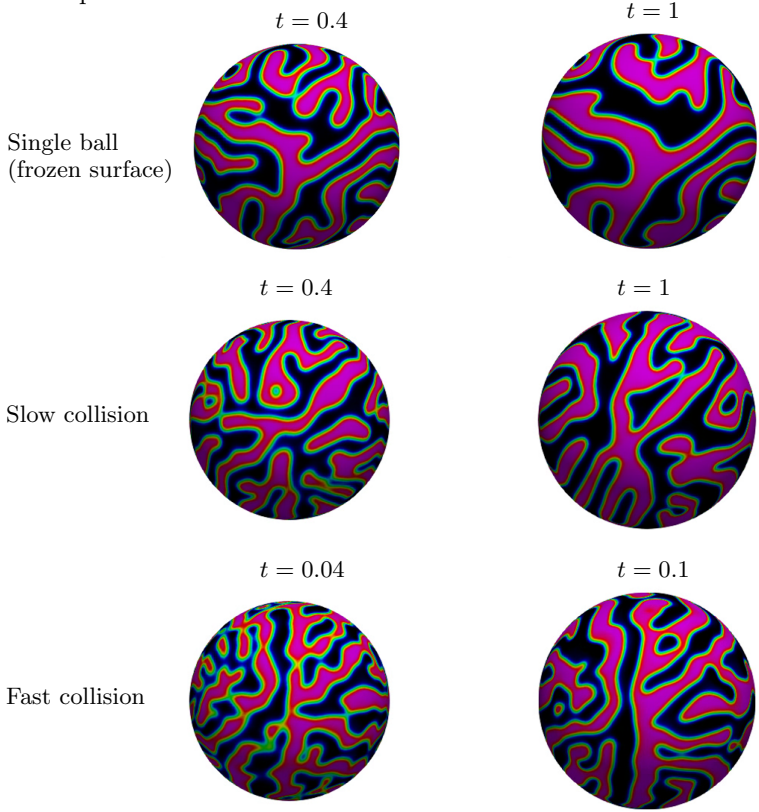

For further comparison, we report in Fig.~\ref{fig:coll_spheres_side} a side view
of the solution at two different times for the
single ball (top), slow collision (center), and fast collision (bottom) tests.
We selected times $t = 0.4,1$ for the single ball and slow collision tests
(same stage of phase separation) and times $t = 0.04,0.1$ for the fast collision test, which feature
the same surface as the slow collision test for  $t = 0.4,1$. Notice that the solution for the fast collision test
is at an earlier stage of phase separation. Given this and the fact that the simulations
are started from a random initial condition, we conclude that the collision speed
does not have a substantial influence on the pattern appearance at the back side of colliding drops.

\subsection{Splitting spheres}\label{sec:splitting}

\begin{figure}[ht]
\begin{center}
\href{https://youtu.be/VTh7UMnMUgM}{
\begin{overpic}[height=.22\textwidth, viewport=130 168 530 470, clip,grid=false]{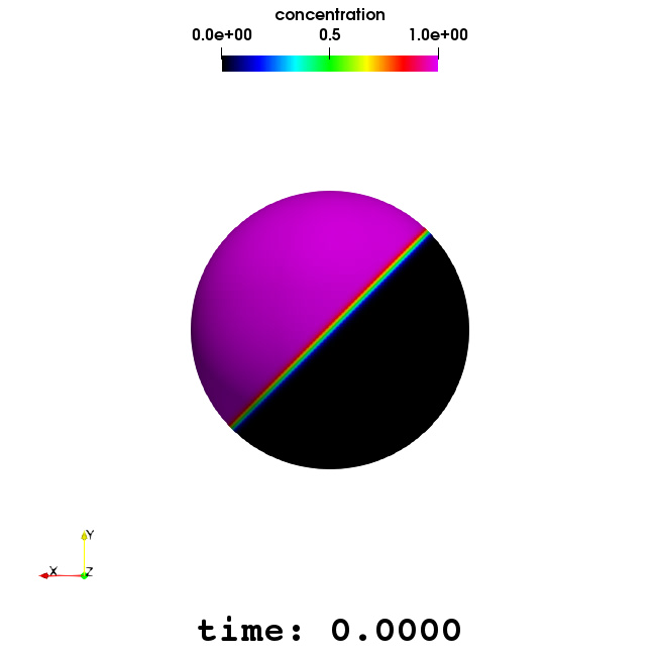}
         \put(42, 85){\small{$t = 0$}}
\end{overpic}
\begin{overpic}[height=.22\textwidth, viewport=130 168 530 470, clip,grid=false]{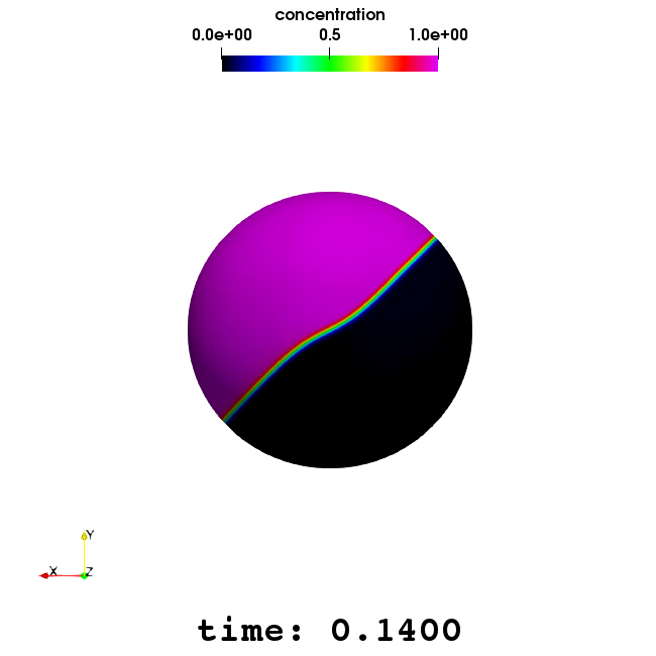}
         \put(42, 85){\small{$t = 0.15$}}
\end{overpic}
\begin{overpic}[height=.22\textwidth, viewport=130 168 530 470, clip,grid=false]{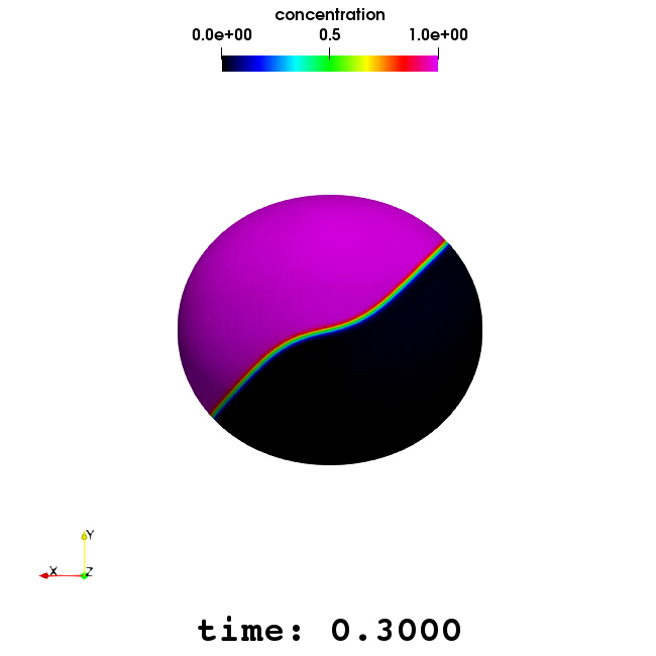}
         \put(42, 85){\small{$t = 0.3$}}
\end{overpic}
}
\vskip .4cm
\href{https://youtu.be/VTh7UMnMUgM}{
\begin{overpic}[height=.22\textwidth, viewport=130 168 530 470, clip,grid=false]{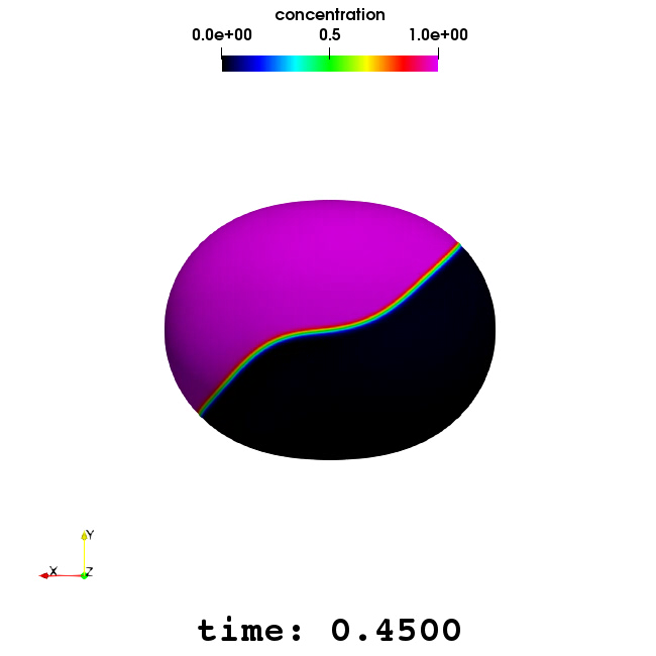}
         \put(40, 77){\small{$t = 0.45$}}
\end{overpic}
\begin{overpic}[height=.22\textwidth, viewport=130 168 530 470, clip,grid=false]{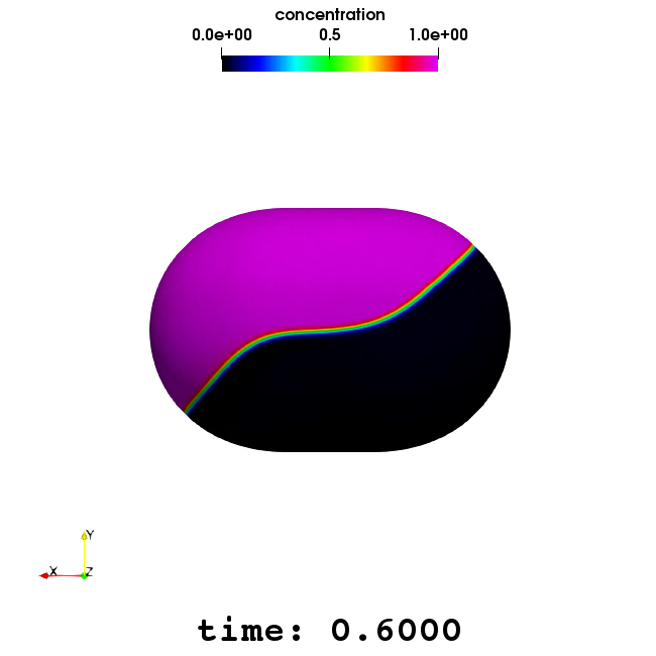}
         \put(40, 77){\small{$t = 0.6$}}
\end{overpic}
\begin{overpic}[height=.22\textwidth, viewport=110 168 550 470, clip,grid=false]{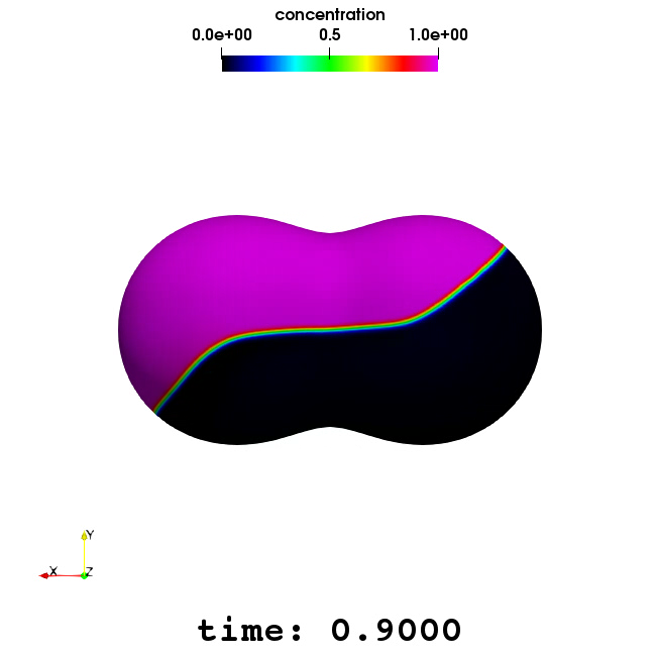}
         \put(40, 70){\small{$t = 0.9$}}
\end{overpic}
}
\vskip .4cm
\href{https://youtu.be/VTh7UMnMUgM}{
\begin{overpic}[width=.48\textwidth, viewport=50 198 616 440, clip,grid=false]{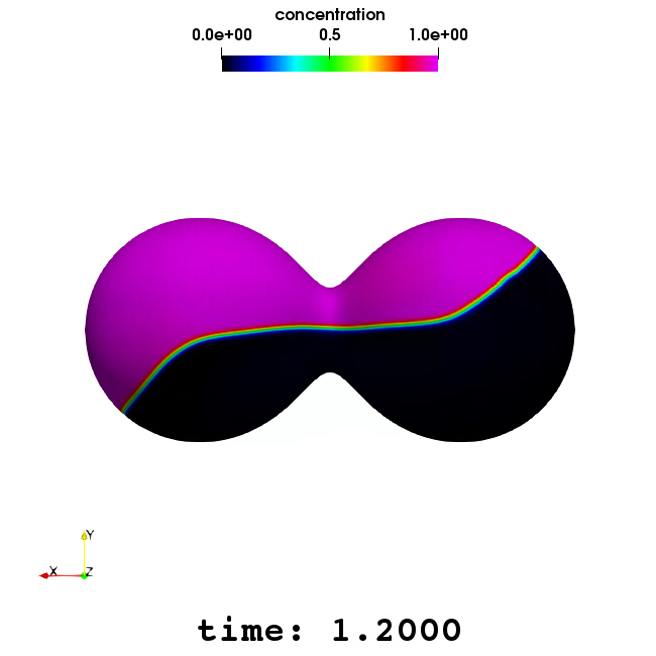}
         \put(40, 44){\small{$t = 1.2$}}
\end{overpic}
\begin{overpic}[width=.48\textwidth, viewport=50 198 616 440, clip,grid=false]{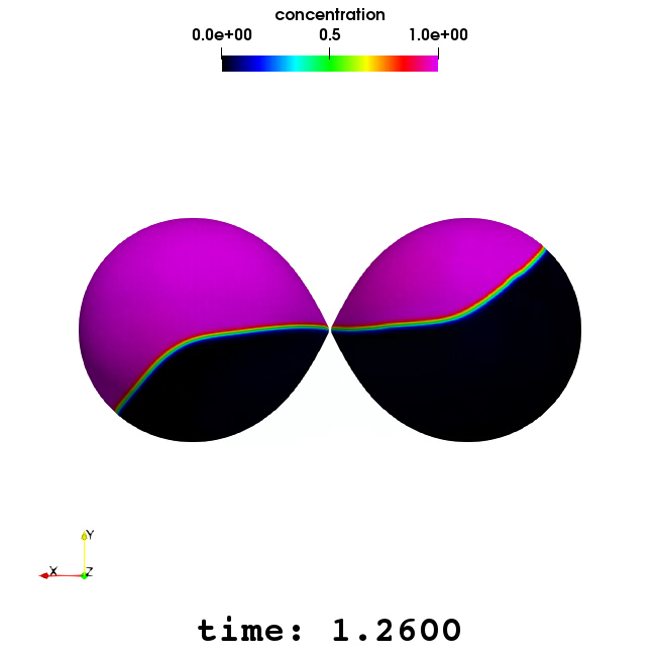}
         \put(40, 44){\small{$t = 1.26$}}
\end{overpic}
}
\vskip .4cm
\href{https://youtu.be/VTh7UMnMUgM}{
\begin{overpic}[width=.48\textwidth, viewport=50 198 616 440, clip,grid=false]{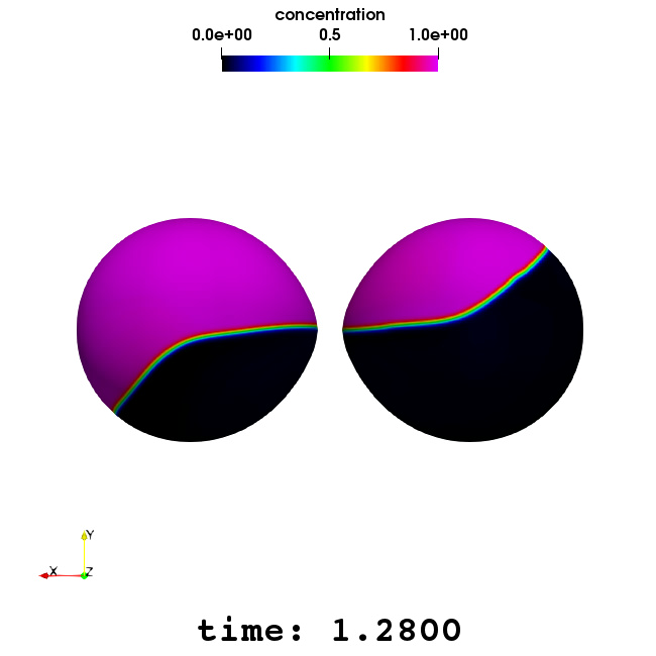}
         \put(40, 44){\small{$t = 1.28$}}
\end{overpic}
\begin{overpic}[width=.48\textwidth, viewport=50 198 616 440, clip,grid=false]{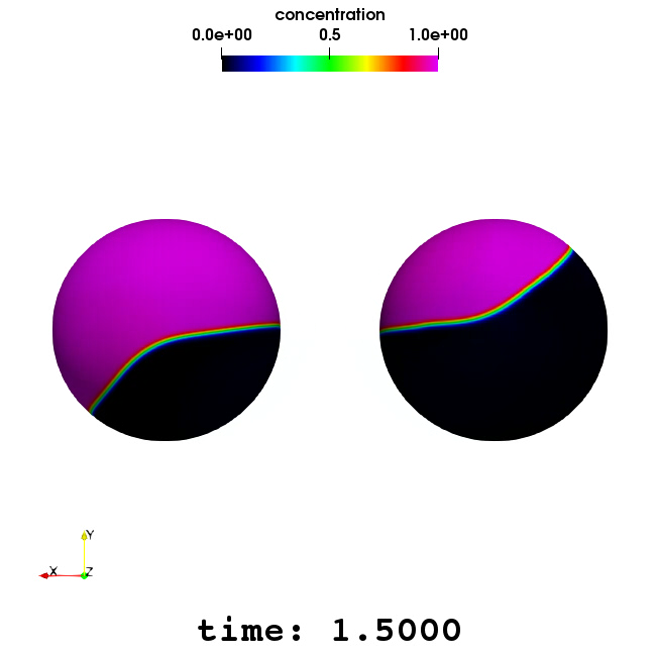}
         \put(40, 44){\small{$t = 1.5$}}
\end{overpic}
}
\end{center}
\caption{\label{fig:split_sphere_rotated}
Splitting spheres, rotated interface: Evolution of the numerical solution of the Cahn--Hilliard
equation computed with mesh $\ell = 6$ for $t \in [0, 1.5]$. View: $x_1 x_2$-plane.
The legend is the same as in Fig.~\ref{fig:test2b}. Click any picture above to run the full animation.}
\end{figure}

\begin{figure}[ht]
\begin{center}
\href{https://youtu.be/MHnn2eliqvA}{
\begin{overpic}[height=.22\textwidth, viewport=130 168 530 470, clip,grid=false]{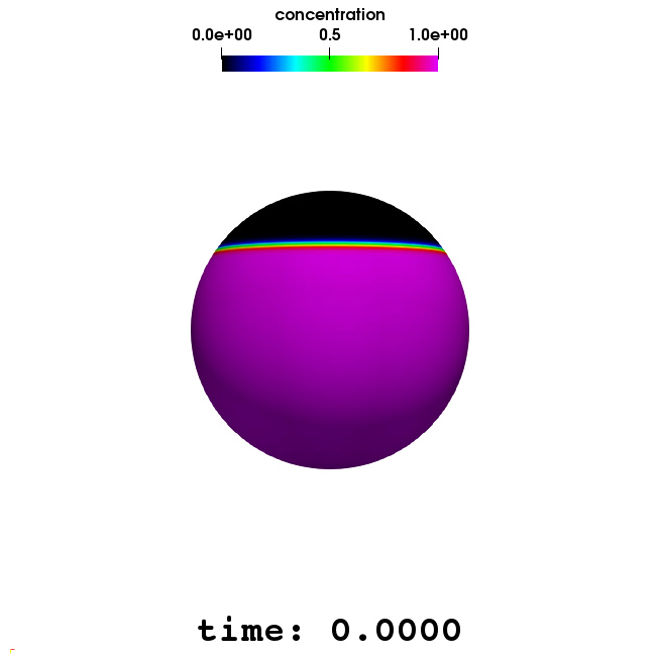}
         \put(42, 85){\small{$t = 0$}}
\end{overpic}
\begin{overpic}[height=.22\textwidth, viewport=130 168 530 470, clip,grid=false]{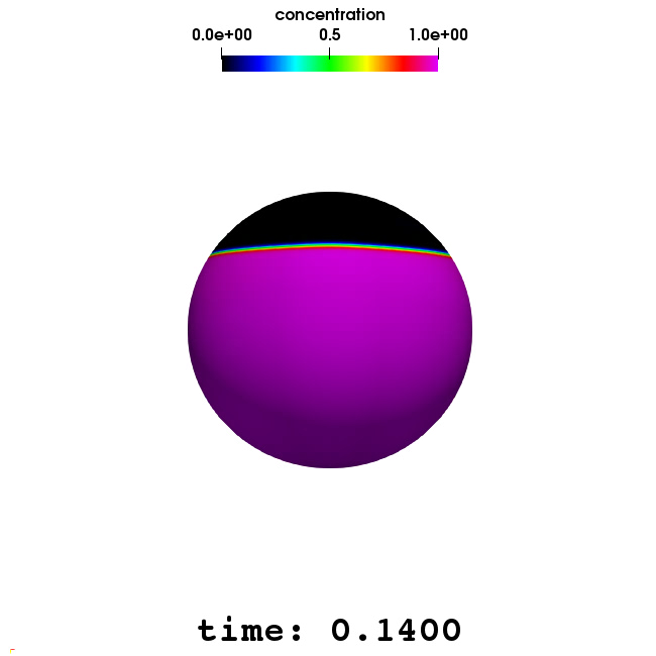}
         \put(42, 85){\small{$t = 0.15$}}
\end{overpic}
\begin{overpic}[height=.22\textwidth, viewport=130 168 530 470, clip,grid=false]{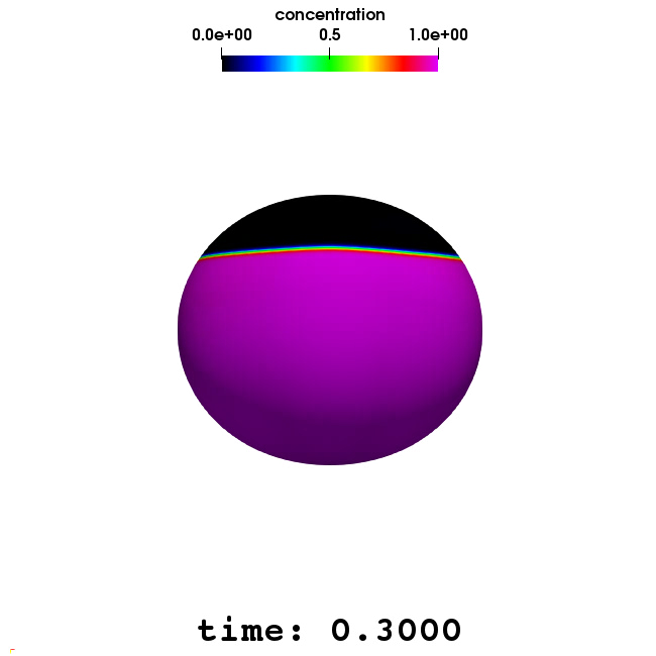}
         \put(42, 85){\small{$t = 0.3$}}
\end{overpic}
}
\vskip .4cm
\href{https://youtu.be/MHnn2eliqvA}{
\begin{overpic}[height=.22\textwidth, viewport=130 168 530 470, clip,grid=false]{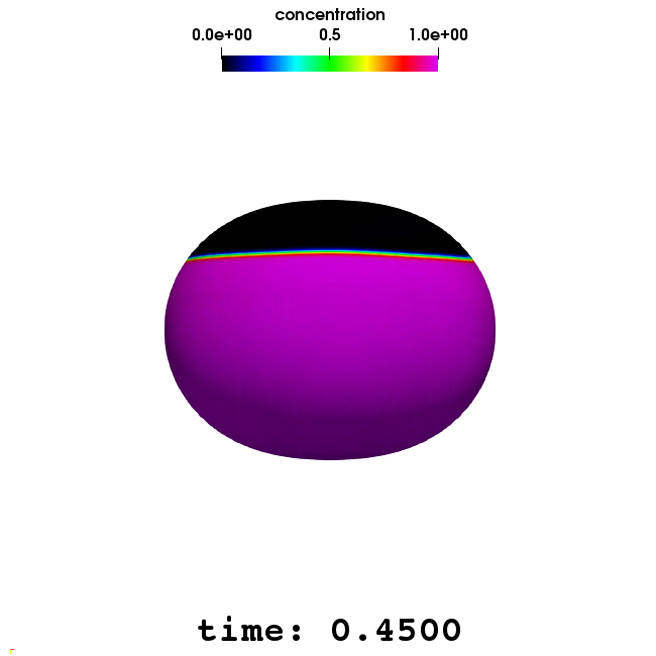}
         \put(40, 80){\small{$t = 0.45$}}
\end{overpic}
\begin{overpic}[height=.22\textwidth, viewport=130 168 530 470, clip,grid=false]{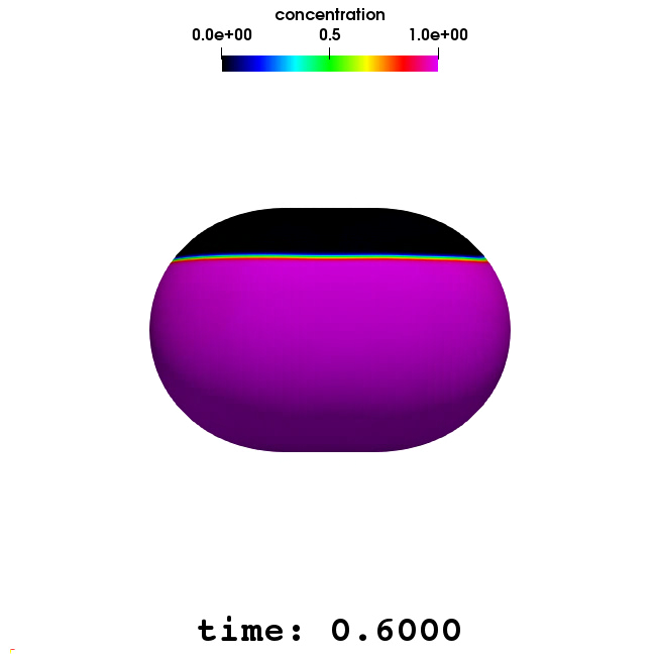}
         \put(40, 80){\small{$t = 0.6$}}
\end{overpic}
\begin{overpic}[height=.22\textwidth, viewport=110 168 550 470, clip,grid=false]{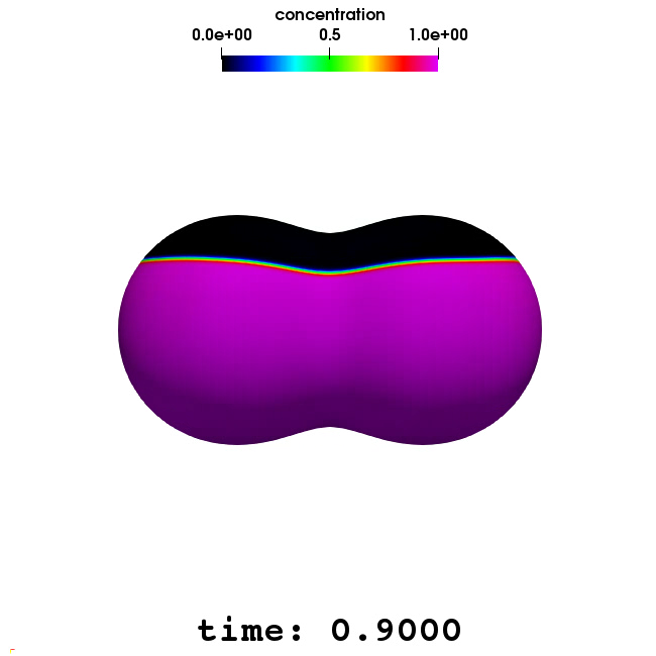}
         \put(40, 67){\small{$t = 0.9$}}
\end{overpic}
}
\vskip .4cm
\href{https://youtu.be/MHnn2eliqvA}{
\begin{overpic}[width=.48\textwidth, viewport=50 198 616 440, clip,grid=false]{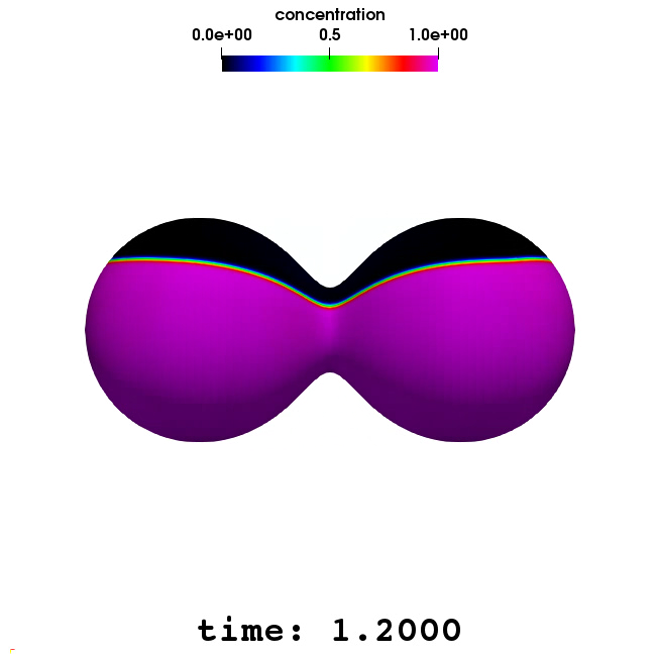}
         \put(40, 44){\small{$t = 1.2$}}
\end{overpic}
\begin{overpic}[width=.48\textwidth, viewport=50 198 616 440, clip,grid=false]{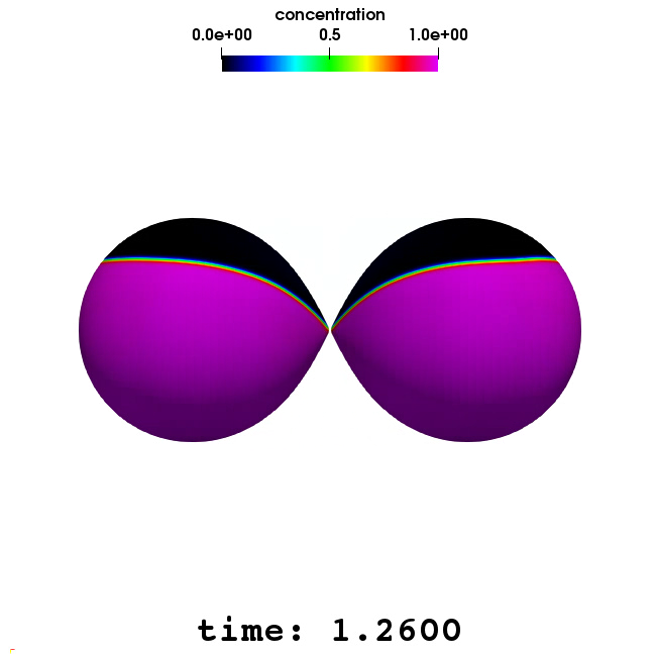}
         \put(40, 44){\small{$t = 1.26$}}
\end{overpic}
}
\vskip .4cm
\href{https://youtu.be/MHnn2eliqvA}{
\begin{overpic}[width=.48\textwidth, viewport=50 198 616 440, clip,grid=false]{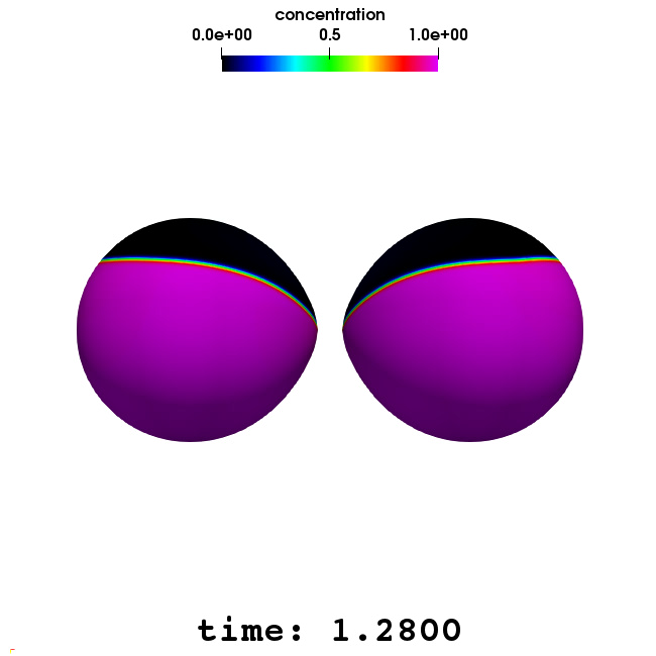}
         \put(40, 44){\small{$t = 1.28$}}
\end{overpic}
\begin{overpic}[width=.48\textwidth, viewport=50 198 616 440, clip,grid=false]{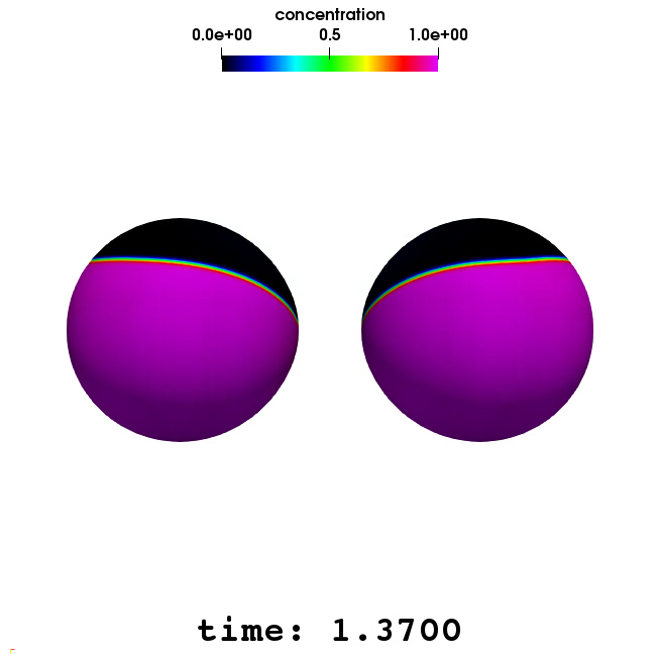}
         \put(40, 44){\small{$t = 1.37$}}
\end{overpic}
}
\end{center}
\caption{\label{fig:split_sphere_horizontal}
Splitting spheres, horizontal interface: Evolution of the numerical solution of the Cahn--Hilliard
equation computed with mesh $\ell = 6$ for $t \in [0, 1.37]$. View: $x_1 x_2$-plane.
The legend is the same as in Fig.~\ref{fig:test2b}. Click any picture above to run the full animation.}
\end{figure}

We consider an evolving surface undergoing a reverse dynamics with respect to the
one considered in Sec.~\ref{sec:merging}: the initial configuration is a sphere centered at the origin,
which splits into two droplets, in turn evolving towards spheres centered at
$\bx_c^\pm(T)$. The computational domain is again $\Omega = [-10/3,10/3]\times[-5/3,5/3]^2$.
The evolving surface $\Gamma(t)$ is the zero level set of $\phi(\bx, 1.5/|w|+t)$ with $w=-1$, where $\phi$ is given by \eqref{eq:merging_level_set}. The velocity vector field is defined by \eqref{eq:merging_velocity} with $w=-1$.
Surface $\Gamma(t)$ is simply connected for $t < 1.5-\tilde{t}$, where
$\tilde{t} = 0.235$.
Obviously, this evolving surface experiences a local singularity too.

We consider the same meshes used for the results reported in Sec.~\ref{sec:merging}.
Below we present the results for 2 numerical experiments:
\begin{itemize}
\item[-] \emph{Rotated interface}: with initial solution is given by $c_0^{x_3}$ from \eqref{steady_CH}, with $\delta=0$ and rotated by $\pi/4$ in $x_1x_2$ plane.

\item[-] \emph{Horizontal interface}: with initial solution $c_0^{x_3}$ from \eqref{steady_CH}, with $\delta=20$.
\end{itemize}

Fig.~\ref{fig:split_sphere_rotated} and \ref{fig:split_sphere_horizontal} show the evolution
of the interface and the solution for $t \in [0,1.5]$ for the rotated and horizontal interface
tests, respectively.
We see that as the two hemispheres are being pulled in opposite directions along the $x_1$-axis,
a horizontal interface forms between the pink domain and the black domain
in both figures.
In the case of the rotated interface test, the horizontal portion of the interface is located symmetrically
with respect to the $x_1$-axis. Thus, it does not get affected by the bottleneck formation.
On the other hand, the horizontal interface in Fig.~\ref{fig:split_sphere_horizontal} is located non-symmetrically
with respect to the $x_1$-axis. Hence, when the bottleneck forms the interface curves.
As a result, the two balls into which the surface has evolved at $t = 1.5$ present
 a curved interface. See Fig.~\ref{fig:split_sphere_horizontal}, bottom right panel.

\section{Conclusions}\label{sec:concl}

We presented a formulation of the Cahn--Hilliard equation
on a time-dependent surface $\Gamma(t)$ that uses tangential calculus
induced by the embedding of $\Gamma(t)$ in $\mathbb{R}^3$.
This makes the formulation particularly suitable for the development
of fully Eulerian numerical techniques for the solution of the problem.
The fully Eulerian technique proposed in this paper is
based on the unfitted trace finite element method for spatial discretization and
the finite difference method combined with an implicit extension procedure
for time discretization.
The degrees of freedom are tailored to a background time-independent mesh and
we used the usual nodal basis to build the systems of algebraic equations at each time step;
space-time integration is avoided. All these ingredients lead to a
rather straightforward implementation of the method in a standard finite element software.
Numerical experiments demonstrated the following main properties of the proposed approach:
(i) optimal second order accuracy for $P_1$ approximation of the concentration, chemical potential,
and the level set function that defines the surface;
(ii) numerical stability even in the case of fast and large deformations;
(iii) the ease of handling phase transition on a surface undergoing topological changes.
The third feature was previously available only for methods based on phase-field  representation of the surface itself.

There are many directions to further develop the proposed method.
Some of them are: employing $P^k$ ($k>1$) trace finite elements together with higher order isoparametric surface recovery; the coupling of lateral phase separation with the surface evolution through the line tension forces; accounting for the surface fluidity and bulk phenomena. We believe that all these extensions fit well within the proposed framework and we plan to address some of them in the future.

\bibliographystyle{siam}
\bibliography{literatur}{}

\end{document}